\documentclass[UTF-8,reqno]{amsart}
\usepackage{enumerate, bbm}
\usepackage{amssymb,url,color, booktabs}
\usepackage{tikz}
\usepackage{mathrsfs}
\usepackage{dutchcal}
\usepackage{color}
\usepackage[colorlinks=true]{hyperref}
\hypersetup{
    linkcolor=blue,          % color of internal links
    citecolor=red,        % color of links to bibliography
    filecolor=blue,      % color of file links
    urlcolor=cyan
}

\usepackage{color}
\usepackage{ulem}

\definecolor{MyDarkBlue}{cmyk}{0.8,0.3,0.8,0.4}
\definecolor{yellow}{rgb}{0.99,0.99,0.70}
\definecolor{white}{rgb}{1.0,1.0,1.0}
\definecolor{black}{rgb}{0.00,0.00,0.00}

\newcommand{\red}{\color{red}}
\numberwithin{equation}{section}

\newcommand{\be}{\begin{eqnarray}}
\newcommand{\ee}{\end{eqnarray}}
\newcommand{\ce}{\begin{eqnarray*}}
\newcommand{\de}{\end{eqnarray*}}
\newtheorem{theorem}{Theorem}[section]
\newtheorem{lemma}[theorem]{Lemma}
\newtheorem{remark}[theorem]{Remark}
\newtheorem{definition}[theorem]{Definition}
\newtheorem{proposition}[theorem]{Proposition}

\newtheorem{corollary}[theorem]{Corollary}

\newcommand\wleqq{\, \tikz[baseline]{\draw (0pt,7pt)--(4pt,5pt) (0pt,6pt)--(4pt,4pt)--(0pt,2pt)--(4pt,0pt) (0pt,1pt)--(4pt,-1pt);} \,}

\newcommand\wle{\, \tikz[baseline]{\draw (0pt,7pt)--(4pt,5pt) (0pt,6pt)--(4pt,4pt)--(0pt,2pt)--(4pt,0pt)--(4pt,0pt);} \,}

\def\eps{\varepsilon}

\def\e{\mathrm{e}}

\def\p{\partial}

\def\[{{\Big[}}
\def\]{{\Big]}}
\def\<{{\langle}}
\def\>{{\rangle}}
\def\({{\big(}}
\def\){{\big)}}

\def\dif{{\mathord{{\rm d}}}}

\def\min{{\mathord{{\rm min}}}}

\def\bbp{{\boldsymbol{p}}}

\def\bbq{{\boldsymbol{q}}}
\def\bba{{\boldsymbol{a}}}
\def\bbw{{\boldsymbol{w}}}
\def\bb2{{\boldsymbol{2}}}
\def\no{\nonumber}
\def\={&\!\!=\!\!&}

\def\bB{{\mathbf B}}

\def\bC{{\mathbf C}}

\def\cB{{\mathcal B}}

\def\cH{{\mathcal H}}

\def\cK{{\mathcal K}}

\def\cM{{\mathcal M}}

\def\cP{{\mathcal P}}

\def\cR{{\mathcal R}}
\def\cS{{\mathcal S}}

\def\cU{{\mathcal U}}

\def\mB{{\mathbb B}}

\def\mE{{\mathbb E}}

\def\mI{{\mathbb I}}

\def\mL{{\mathbb L}}

\def\mN{{\mathbb N}}

\def\mR{{\mathbb R}}
\def\mS{{\mathbb S}}

\def\mX{{\mathbb X}}
\def\mY{{\mathbb Y}}
\def\mZ{{\mathbb Z}}

\def\bP{{\mathbf P}}

\def\bD{{\mathbf H}}
\def\DD{{H}}

\def\b1{{\mathbbm 1}}

\def\sA{{\mathscr A}}
\def\sB{{\mathscr B}}

\def\sF{{\mathscr F}}
\def\sG{{\mathscr G}}

\def\sI{{\mathscr I}}

\def\sQ{{\mathscr Q}}
\def\sR{{\mathscr R}}

\def\sU{{\mathscr U}}

\def\geq{\geqslant}
\def\leq{\leqslant}
\def\ge{\geqslant}
\def\le{\leqslant}

\def\div{\mathord{{\rm div}}}

\def\eps{\varepsilon}

\def\e{\mathrm{e}}

\def\p{\partial}

\def\[{{\Big[}}
\def\]{{\Big]}}
\def\<{{\langle}}
\def\>{{\rangle}}

\def\bu{{\mathbf{u}}}

\def\dif{{\mathord{{\rm d}}}}

\def\min{{\mathord{{\rm min}}}}

\def\no{\nonumber}
\def\={&\!\!=\!\!&}
\def\bt{\begin{theorem}}
\def\et{\end{theorem}}
\def\bl{\begin{lemma}}
\def\el{\end{lemma}}
\def\br{\begin{remark}}
\def\er{\end{remark}}
\def\bd{\begin{definition}}
\def\ed{\end{definition}}
\def\bp{\begin{proposition}}
\def\ep{\end{proposition}}
\def\bc{\begin{corollary}}
\def\ec{\end{corollary}}

\def\geq{\geqslant}
\def\leq{\leqslant}
\def\ge{\geqslant}
\def\le{\leqslant}

\def\div{\mathord{{\rm div}}}

\def\bH{{\mathbf H}}

\def\<{\langle} \def\>{\rangle}

\def\wt{\widetilde}

\def\bbp{{\boldsymbol{p}}}
\def\bbl{{\boldsymbol{\ell}}}
\def\bbg{{\boldsymbol{\gamma}}}

\def\bbq{{\boldsymbol{q}}}
\def\bbd{{\boldsymbol{d}}}
\def\bba{{\boldsymbol{a}}}
\def\bbk{{\boldsymbol{k}}}
\def\bbw{{\boldsymbol{w}}}
\def\bb2{{\boldsymbol{2}}}
\def\bbb1{\boldsymbol{1}}

\def\no{\nonumber}
\def\={&\!\!=\!\!&}

\def\bbrho{\boldsymbol{\varrho}}

\allowdisplaybreaks

\begin{document}

\title[Second order Mean-field  SDEs]
{Second order fractional  mean-field  SDEs with singular kernels and measure initial data}

\author{Zimo Hao, Michael R\"ockner and Xicheng Zhang}

\address{Zimo Hao:
Fakult\"at f\"ur Mathematik, Universit\"at Bielefeld,
33615, Bielefeld, Germany\\ 
and School of Mathematics and Statistics, Wuhan University,
Wuhan, Hubei 430072, P.R.China\\
Email: zhao@math.uni-bielefeld.de}

\address{Michael R\"ockner:
Fakult\"at f\"ur Mathematik, Universit\"at Bielefeld,
33615, Bielefeld, Germany\\
and Academy of Mathematics and Systems Science, CAS, Beijing, China\\
Email: roeckner@math.uni-bielefeld.de
 }

\address{Xicheng Zhang:
School of Mathematics and Statistics, Beijing Institute of Technology, Beijing 100081, China\\
Email: XichengZhang@gmail.com
 }

\thanks{
This work is partially supported by NNSFC grants of China (Nos. 12131019), and the German Research Foundation (DFG) through the 
Collaborative Research Centre(CRC) 1283/2 2021 - 317210226 ``Taming uncertainty and profiting from randomness and low regularity in analysis, stochastics and their applications".}

\begin{abstract}
In this paper we establish the local and global well-posedness of weak and strong solutions to second order fractional mean-field SDEs with singular/distribution interaction kernels and measure initial value, where the kernel can be
Newton or Coulomb potential, Riesz potential, Biot-Savart law, etc. Moreover, we also show the stability, smoothness and 
the short time singularity and large time decay estimates of the distribution density.
Our results reveal a phenomenon that for {\it nonlinear} mean-field equations, the regularity of the initial distribution
could balance the singularity of the kernel. The precise relationship between the singularity of kernels and the regularity of
initial values are calculated, which belongs to the subcritical regime in the scaling sense. 
In particular, our results  provide a microscopic probabilistic explanation 
and establish a unified treatment for 
many physical models such as the fractional Vlasov-Poisson-Fokker-Planck system, the vorticity formulation of 2D-fractal Navier-Stokes equations, 
surface quasi-geostrophic models, 
fractional porous medium equation with viscosity, etc.

\bigskip

%As applications we also discuss the vanishing viscosity and show the existence of weak solutions for kinetic aggregation equations.

\noindent 
\textbf{Keywords}: Mean-field SDEs, Anisotropic Besov spaces, Coulomb potential, Riesz potential.\\

\noindent
 {\bf AMS 2010 Mathematics Subject Classification:} 60H10, 35H10
\end{abstract}

\maketitle \rm

\tableofcontents

\section{Introduction}

In this paper we are concerned with the following second order mean-field stochastic differential equation in $\mR^d$ (also called McKean-Vlasov or distribution-dependent SDE) driven by standard Brownian motion $W_t$,
\begin{align}\label{MV00}
\ddot X_t=(b_t*\mu_{t})(X_t,\dot X_t)+\sqrt2\dot W_t,\ \ (X_0,\dot X_0)\sim \mu_0,
\end{align}
or by an $\alpha$-stable process $L^{(\alpha)}_t$ with $\alpha\in(1,2)$, whose infinitesimal generator is given by the usual fractional Laplacian $\Delta^{\alpha/2}$,
\begin{align}\label{MV0}
\ddot X_t=(b_t*\mu_{t})(X_t,\dot X_t)+\dot L^{(\alpha)}_t,\ \ (X_0,\dot X_0)\sim \mu_0,
\end{align}
where $\mu_0\in\cP(\mR^{2d})$ is the initial probability distribution on $\mR^{2d}$, 
$\mu_t$ stands for the joint law of $(X_t,\dot X_t)$, 
and $b$ is a time-dependent $\mR^d$-valued (Schwartz) distribution over the phase space $\mR^{2d}$,
$b_t*\mu_t$ is the usual convolution of two distributions. Here we require that $b*\mu\in L^1_{loc}(\mR_+\times\mR^{2d})$.
%Note that if $b_t*\mu_t$ is a distribution, we have to make a meaning for the drift $(b_t*\mu_{t})(X_t,\dot X_t)$.
Below for the unification of notations, we shall write
$$
L^{(2)}_t=\sqrt{2}W_t.
$$
Note that if we introduce the velocity variable $V_t:=\dot X_t$, then SDEs \eqref{MV00} and \eqref{MV0} can be written as the following first-order degenerate system: for $\alpha\in(1,2]$, 
\begin{align}\label{MV1}
\left\{
\begin{aligned}
&X_t=X_0+\int_0^t V_s\dif s,\\
&V_t=V_0+\int^t_0(b_s*\mu_{s})(X_s,V_s)\dif s+L^{(\alpha)}_t.
\end{aligned}
\right.
\end{align}
It should be noted that if $b(t,x,v)=b(t,v)$ does not depend on the position variable $x$, then SDE \eqref{MV1} reduces to the following first-order nondegenerate SDE
\begin{align}\label{MV2}
\dif V_t=(b_t*\mu_t)(V_t)+\dif L^{(\alpha)}_t.
\end{align}
Suppose that $\mu_t(\dif x,\dif v)=f(t,x,v)\dif x\dif v$. 
Then by It\^o's formula, $f$ solves the following nonlinear kinetic Fokker-Planck equation in the distributional sense:
\begin{align}\label{FPE00}
\p_t f+v\cdot\nabla_xf=\Delta^{\frac\alpha 2}_vf-\div_v ((b*f)f).
\end{align}
If $b(t,x,v)=b(t,v)$ does not depend on $x$, then $\rho(t,v):=\int_{\mR^d}f(t,x,v)\dif x$ solves the following nonlinear and nondegenerate Fokker-Planck equation:
\begin{align}\label{FPE011}
\p_t \rho=\Delta_v^{\frac\alpha 2}\rho-\div_v ((b*\rho)\rho).
\end{align}

The aim of this work is to show the weak and strong well-posedness of SDEs \eqref{MV1} and
\eqref{MV2} with singular $b$ for lowly regular initial values via studying PDEs \eqref{FPE00} and \eqref{FPE011}.
Especially, we find a precise relationship between the singularity of the kernel $b$ and the regularity of $\mu_0$ in terms of their Besov norms, which is almost sharp and belongs to the subcritical regime in the scaling sense. In particular, 
it allows kernels in \eqref{FPE011} 
like $b(x)=\nabla |x|^{s-d}$ for some $s\in(2-\alpha,d)$ and $\rho_0\in L^1\cap L^{p_0}$ for some $p_0>\frac{d}{s-2+\alpha}$, see Remark \ref{Re10} below.
%We shall explain the singularity of $b$ in the next subsection.

\medskip

\subsection{Background and motivation}

%In this subsection we introduce several concrete physical models that relate to SDEs \eqref{MV1} and \eqref{MV2}
%or Fokker-Planck equations \eqref{FPE00} and \eqref{FPE011}.

%\noindent {\bf Example 1} Let us 
Consider the interacting many-body system  of $N$-particles in $\mR^d$ governed by Newtonian mechanics and stochastic noises, 
\begin{align}\label{Par1}
\left\{
\begin{aligned}
&\dot X^{N,i}_t=V^{N,i}_t,\ \ i=1,2,\cdots,N,\\
& \dot V^{N,i}_t=\frac1N\sum_{j\not=i}K(X^{N,i}_t-X^{N,j}_t)+\sqrt{2\eps}\dot W^{N,i},
\end{aligned}
\right.
\end{align}
and subjected to i.i.d. initial value $(X^{N,i}_0,  V^{N,i}_0)=(\xi^i_0,\eta^i_0)$,  whose law is 
absolutely continuous or singular with respect to the Lebesgue measure,
where $(X^{N,i}_t,V^{N,i}_t)$ stands for the position and velocity of the $i$th-particle at time $t$, 
$K:\mR^d\to\mR^d$ models the pairwise interaction between different particles,
$(\dot W^{N,i})_{i=1,\cdots,N}$ are the generalized derivatives of $N$ independent Brownian motions, which models the collision between particles and
a background medium, $\eps>0$ stands for the intensity of the noise, see \cite{VD90} for the physical background.
In this work we are mainly interested in the interaction kernel $K=A\cdot\nabla U$ (see \cite{Go16}), where $A$ is some $d\times d$-matrix and
\begin{align}\label{UU1}
U(x)=|x|^{-s},\ \ s\in(0,d) \ \mbox{ or } \ U(x)=\ln|x|,\ \ d=1,2.
\end{align}
Note that $U(x)=|x|^{2-d}$ for $d\geq 3$ or $U(x)=\ln|x|$ for $d=2$ is the classical Newton or Coulomb kernel,
$U(x)=|x|^{1-d}$ for $d\geq 2$ is called the Riesz kernel.
Since it is hard to directly study the particle system \eqref{Par1} by observing the motion of each particle
in physical models, one usually considers the empirical measure
$$
\mu^N_t(\dif x,\dif v):=\frac1N\sum_{i=1}^N\delta_{(X^{N,i}_t, V^{N,i}_t)}(\dif x,\dif v).
$$
{ It is assumed that $\mu^N_t$} weakly converges to some measure of type $f(t,x,v)\dif x\dif v$, which is also called
mean-field limit or propagation of chaos in Kac's sense (see \cite{Ja14}, \cite{Go16} and \cite{CD21}), where the density $f$ satisfies the
following Vlasov-Poisson-Fokker-Planck equation (abbreviated as VPFP)
\begin{align}\label{VP1}
\begin{split}
\p_t f+v\cdot\nabla_xf&=\eps\Delta_vf-\div_v(K*\<f\>f )
=\eps\Delta_vf-(K*\<f\>)\cdot\nabla_v f.
\end{split}
\end{align}
Here $\<f\>(t,x):=\int_{\mR^d}f(t,x,v)\dif v$ is the mass or charge density.
%, the asterisk stands for the convolution.
%$E(t,x):=\nabla U*\<f\>(t,x)$ is called the gravitational or Coulombian force field. 
PDE \eqref{VP1} can be regarded as a macroscopic description of the microscopic particle system \eqref{Par1}.

In the three dimensional case, when $\eps=0$, that is, there is no viscosity in \eqref{VP1},  and $U$ is the Newton/Coulomb potential, which corresponds to gravitational/electrostatic interaction, 
PDE \eqref{VP1} is also called Vlasov-Poisson equation and describes the evolution of the distribution density of particles in the phase space $\mR^6$ and appears in galactic dynamics and plasma physics (see \cite{MoVi11} and references therein). 
No viscosity implies that PDE \eqref{VP1} is a time-reversible system in the sense that 
if $f$ is a solution of \eqref{VP1}, then  $\tilde f(t,x,v):=f(-t,x,-v)$ is also a solution of \eqref{VP1}.
The existence of global smooth solutions to the Vlasov-Poisson equation
was shown in \cite{Li-Pe91}, \cite{Pfa92}, \cite{Sch91} under some moment assumptions on the initial value.
The existence of classical solutions to the VPFP equation was obtained in \cite{Bou93}, \cite{CS97}, \cite{OS00} by the same method as in \cite{Li-Pe91}.
All the results require that the initial value $f_0\in L^1\cap L^\infty$ has a finite second order moment in $v$, that is,
\begin{align}\label{20230303}
\int_{\mR^{2d}}|v|^2 f_0(x,v)\dif x\dif v<\infty.
\end{align}
When the initial distribution $f_0$ is a measure,  say in some Morrey space, and satisfies some smallness assumptions,
the well-posedness of VPFP \eqref{VP1} was established in \cite{CS97} (see the earlier work \cite{ZM94} for one dimensional case).
{On the other hand, from the physical viewpoint, the mean-field limit provides an explanation for Vlasov-Poisson equations 
from Newton's second law for the motion of $N$ particles. Such a problem 
has been studied extensively in the context of PDEs (see \cite{Ja14} and \cite{Go16}). }
When $K$ is Lipschitz continuous, the propagation of chaos for particle system \eqref{Par1}
was firstly shown by McKean \cite{McK67}. When the interaction kernel $K$ is bounded measurable,  Jabin and Wang \cite{JW16}
showed the quantitative rate of the propagation of chaos via the entropy method. Based on the BBGKY method, Lacker \cite{La21}
showed the optimal rate of convergence
by the Girsanov transform and deep analysis. For singular potentials $U$, there is also some cutoff version about the mean-field limits,
see \cite{HLP20} and references therein.
A complete review of mean field limits  for Vlasov-Poisson equations is contained in \cite{Ja14} and \cite{Go16}.

From the probabilistic view point, it is natural to expect that the  distributional density solution $f$ of VPFP \eqref{VP1} is the law
of the solution $Z_t=(X_t,V_t)$ to the following mean-field SDE, if a solution exists in It\^o's sense:
\begin{align}\label{MV3}
\dot X_t=V_t,\ \ \dot V_t=\int_{\mR^d}K(X_t-y)\mu_{X_t}(\dif y)+\sqrt{2\eps}\dot W_t,
\end{align}
where $\mu_{X_t}$ stands for the probability distribution measure of $X_t$. As discussed above, although there are a lot of references about
the solvability of PDE \eqref{VP1}, there are few results about SDE \eqref{MV3}, e.g., about its 
weak/strong well-posedness in the probabilistic sense. The difficulty of solving \eqref{MV3} lies in that DDSDE 
\eqref{MV3} is degenerate and the drift is too singular for $K=A\cdot\nabla U$ with $U$ given in \eqref{UU1}.

Besides the second order particle system \eqref{Par1}, we are also interested in the following first order system 
that is related to the point vortex system and other physical models:
\begin{align}\label{SQ0}
\dot X^{N,i}_t=\frac1N\sum_{j\not=i}K(X^{N,i}_t-X^{N,j}_t)+\sqrt{2\eps}\dot W^{N,i}_t.
\end{align}
As above, {it is usually assumed that} the limit of the empirical measure $\mu^N_t(\dif x):=\frac1N\sum_{i=1}^N\delta_{X^{N,i}_t}(\dif x)$ solves the following 
Fokker-Planck equation:
\begin{align}\label{SQ1}
\p_t\rho=\eps\Delta \rho-\div ((K*\rho)\rho),
\end{align}
where $\rho$ is the distribution density of the solution to the following McKean-Vlasov SDE:
\begin{align}\label{SDE19}
X_t=X_0+\int^t_0 (K*\mu_{X_s})(X_s)\dif s+\sqrt{2\eps}W_t.
\end{align}
For the above first order particle system with $K=A\cdot\nabla U$ as described in \eqref{UU1}, 
Serfaty \cite{Ser20} established the conditional mean field convergence 
by modulated free energy. On the two dimensional torus, when $K$ is the Biot-Savart law, Jabin and Wang \cite{JW18}
showed some quantitative rate of convergence of propagation of chaos for the point vortex system by the entropy method.
When the interaction kernel is in $L^p$-space with $p>d$, the strong convergence of propagation of chaos was obtained recently in \cite{RZZ22}. The purpose of the current work is not to prove the propagation of chaos, 
but concentrated on the solvability of SDEs \eqref{MV3} (and hence the corresponding nonlinear Fokker-Planck equations \eqref{VP1} and \eqref{SQ1} respectively) and \eqref{SDE19} with singular interaction kernels.

Before going further, we introduce more examples from physics for equation \eqref{SQ1}.

\noindent {\bf Example 1.} Consider the following vorticity form of the 2D Navier-Stokes equation:
\begin{align}\label{SQ3}
\p_t w=\Delta w+(K_2*w)\cdot\nabla w=\Delta w+\div ((K_2*w) w),
\end{align}
where $K_2$ is the Biot-Savart law given by
$$
K_2(x_1,x_2)=(-x_2,x_1)/|x|^2.
$$
If the initial vorticity $w_0$ is a finite signed measure, the existence of weak solutions to equation \eqref{SQ3} 
was proved by Cottet \cite{Co86} and Giga, Miyakawa and Osada \cite{GMO88}.
The uniqueness is a more difficult problem and was partially proved in \cite{GMO88}, \cite{GW05}and then completely solved by 
Gallagher and Gallay \cite{GG05}. The proofs in \cite{GMO88} and \cite{GG05} strongly depend on the heat kernel estimates
of $\p_t-\Delta-U\cdot\nabla$ established by Osada in \cite{Os87}, 
where $U:\mR_+\times\mR^2\to\mR^2$ is a divergence free vector field which satisfies 
$$
\|U(t)\|_{L^\infty(\mR^2)}\leq C/\sqrt{t},\ \ 
\|\p_2 U_1(t)-\p_1 U_2(t)\|_{L^1(\mR^2)}\leq C.
$$
For the classical Navier-Stokes equation of velocity fields, there are also numerous papers about the initial velocity being in
some homogeneous spaces (see \cite{Ka84}, \cite{GM89}, \cite{Ka94}, \cite{CP96}, \cite{KT01}). It should be noted that 
there is an essential discrepancy between the vorticity formulation and the original two-dimensional
Navier-Stokes equations because if $u\in L^2(\mR^2)$ and $w=\p_2 u_1-\p_1 u_2\in L^1(\mR^2)$, then $\int_{\mR^2}w\dif x=0$.
One of our main results shows that when 
the initial vorticity is in $\cup_{\eps>0}\bB^{-3/2+\eps}_\infty$, then there is a unique smooth solution to equation \eqref{SQ3}.
In particular, it includes the white noise in $\mR^2$ as the initial vorticity, which seems to be new.

\medskip

\noindent {\bf Example 2.}  Consider the following surface quasi-geostrophic equation  in $\mR^2$ (abbreviated by SQG):
for $\alpha\in(0,2)$,
\begin{align}\label{SQG1}
\p_t \theta=\Delta^{\frac\alpha2}\theta+(\sR*\theta)\cdot\nabla \theta,
\end{align}
where $\sR=(-\sR_2,\sR_1)$ and $\sR_i$ is the Riesz transform
$$
\sR_i(x)=x_i/|x|^3,\ x\in\mR^2,\ i=1,2.
$$
The SQG equation appeared in the study of the mathematical theory of meteorology and oceanography \cite{CMT94}, \cite{HPGS95}.
It is classified in subcritical/critical/supercritical cases according to the fractional power 
$\alpha\in(1,2)$/$\alpha=1$/$\alpha\in(0,1)$. In \cite{CW99},
Constantin and Wu showed that in the subcritical case any solution with smooth initial value is smooth for all time. 
In both the critical case and supercritical cases, Chae and Lee \cite{CL03} showed 
the well-posedness of solutions with initial conditions small in Besov spaces (see also \cite{Wu05}). In the critical case, Cafferalli and Vaseur \cite{CV10}
showed the global well-posedness for SQG \eqref{SQG1} for large initial values. All these results require that the initial value is in $L^2$. We are interested in the well-posedness of \eqref{SQG1} when the initial data are measures or certain Schwartz distributions.

\medskip

\noindent {\bf Example 3.}  Let $d\geq 3$, $s\in[0,1]$ and $\eps\geq 0$. 
Consider the following fractional porous medium equation in $\mR^d$ with viscosity:
$$
\p_t u=\eps\Delta u+\div((\nabla\Delta^{-s} u) u ).
$$
For $s\in(0,1]$, it is well known that the nonlocal operator $\Delta^{-s} u$ can be written as (see \cite{St79})
$$
\Delta^{-s} u=K*u,\ \ K(x):=c_{d,s}|x|^{2s-d},
$$
where $c_{d,s}>0$ is a normalizing constant.
Thus, if we take 
$$
U(x)=|x|^{2s-d}\b1_{s\in(0,1]}+\delta_0(\dif x)\b1_{s=0},
$$
then for $s\in[0,1]$,
$$
\nabla\Delta^{-s} u=\nabla U*u.
$$
For $s=0$ and $\eps=0$, it is the classical porous medium equation. For $s=1$ and $\eps=0$, it was investigated in \cite{CRS96} and
\cite{E94}
for the evolution of the vortex density in a superconductor. For $s\in(0,1)$ and $\eps=0$, it is called fractional porous medium equation and
was studied in \cite{BIK15} and \cite{Va17}. 
For $\eps>0$ and $s\in[0,1]$, it was considered as the viscosity approximation of fractional porous medium equation in 
\cite{SV14} and \cite{CHJZ22}. In particular, our results can be applied to the above equation for $\eps>0$. In fact, we shall consider a kinetic version, which seems not physically relevant, but appears interesting for mathematics.

\subsection{Main results}
In this section we provide explicit conditions on $b$ and initial conditions
so that the mean-field SDE \eqref{MV1} admits a unique weak/strong solution in the probabilistic sense.

To state our main results, we first introduce some notations used in this paper.
For a muti-index $\bbp=(p_x,p_v)\in[1,\infty]^2$, we define
\begin{align}\label{LP1}
\|f\|_{\mL^\bbp}:=\|f\|_{\bbp}:=\left(\int_{\mR^d}\|f(\cdot,v)\|_{p_x}^{p_v}\dif v\right)^{1/p_v},
\end{align}
and introduce the scaling parameter $\bba$ by
$$
\bba:=(1+\alpha,1),\ \ \alpha\in(0,2).
$$
For notational convenience, we shall write below
\begin{align*}
\frac{1}{\bbp}:=\Big(\frac{1}{p_x},\frac1{p_v}\Big),\ \ \bba\cdot \frac d\bbp:=\frac{(1+\alpha)d}{p_x}+\frac{d}{p_v},
\end{align*}
and for any $\bbp,\bbq\in[1,\infty]^2$,
$$
\bbp\ge\bbq\Leftrightarrow p_x\ge q_x,\ p_v\ge q_v.
$$
Moreover, we use bold letters to denote constant vectors in $\mR^2$, for example,
\begin{align*}
\bbb1:=(1,1),\quad \bbd:=(d,d).
\end{align*}
We remark that the above mixed $\mL^\bbp$-norm was first introduced in \cite{BP61}. Recently, it was used in an essential way in \cite{RZZ22}
to study the strong convergence of propagation of chaos for particle systems with singular interactions (see also \cite{LX21}, \cite{Z22}).
We emphasize that the order of $(p_x,p_v)$ in definition \eqref{LP1} is important because the above $\mL^\bbp$-norm 
is invariant under the group action $\Gamma_t f(x,v)=f(x-tv,v)$ associated with $\p_t+v\cdot\nabla_v$. More precisely,
for any $\bbp\in[1,\infty]^2$ and $f\in\mL^\bbp$,
\begin{align}\label{AA1}
\|\Gamma_t f\|_{\bbp}=\left(\int_{\mR^d}\|f(\cdot-tv,v)\|_{p_x}^{p_v}\dif v\right)^{1/p_v}=
\left(\int_{\mR^d}\|f(\cdot,v)\|_{p_x}^{p_v}\dif v\right)^{1/p_v}=\|f\|_{\bbp}.
\end{align}

Now if we let $Z_t:=(X_t,V_t)$, then SDE \eqref{MV1} can be written in the following compact form:
\begin{align}\label{MV11}
Z_t=Z_0+\int^t_0\bD_s(Z_s)\dif s+\sigma L^{(\alpha)}_t,
\end{align}
where $\sigma:=(0,\mI)^*\in\mR^{2d}\otimes\mR^d$ and for $z=(x,v)$,
\begin{align}\label{FF}
\bD_t(z):=(v,\DD_t(z))^*\in\mR^{2d},  \  \ \DD_t(z):=(b_t*\mu_t)(z).
\end{align}
We first give the definitions of weak/strong solutions to mean-field SDEs \eqref{MV11}.
\bd\label{Def1}
Let $\alpha\in(1,2]$ and $\sG:=(\Omega,\sF,\bP; (\sF_t)_{t\geq 0})$ be a stochastic basis. Let $Z:=(Z_t)_{t\geq 0}$ and 
$L^{(\alpha)}:=(L^{(\alpha)}_t)_{t\geq 0}$
be two $\mR^{2d}$ and $\mR^d$-valued c\'adl\'ag $\sF_t$-adapted processes. 
\begin{enumerate}[{\bf (A)}]
\item {\bf (Weak solutions)} For given distribution $\mu_0\in\cP(\mR^{2d})$,
we call $(Z,L^{(\alpha)},\sG)$ a weak solution of SDE \eqref{MV11} with initial distribution $\mu_0$ if
\begin{enumerate}[(i)]
\item  $L^{(\alpha)}$ is an $\alpha$-stable process with generator $\Delta^{\alpha/2}$ and $\bP\circ Z^{-1}_0=\mu_0$. 
\item For each $t>0$, $\bP\circ Z^{-1}_t$ admits a {$C^\infty$} density $f(t,z)$ so that
\begin{align}\label{SQ5}
\int^T_0 \sup_{z\in\mR^{2d}}|\DD_s(z)|\dif s<\infty,\ \ \DD_t(z):=(b_t*f(t))(z).
\end{align}
\item For all $t>0$
$$
Z_t=Z_0+\int^t_0\bD_s(Z_s)\dif s+\sigma L^{(\alpha)}_t,\ \ \bP-a.s.
$$
\end{enumerate}
\item {\bf (Strong solutions)} Let   $L^{(\alpha)}$ be an $\alpha$-stable process with generator $\Delta^{\alpha/2}$ and 
$$
\sF^0_t:=\sF_0\vee\sigma\{L^{(\alpha)}_s, s\leq t\}.
$$
For given  $Z_0\in\sF^0_0$,
we call a $\mR^{2d}$-value c\'adl\'ag $\sF^0_t$-adapted process $Z:=(Z_t)_{t\geq 0}$ a strong solution 
of SDE \eqref{MV11} starting from $Z_0$ if (ii) and (iii) above hold.
\end{enumerate}
\ed
\br\rm
The above definitions are standard (for example, see \cite{KS88} and \cite{EK86}), where we used the convention that the $2$-stable process with generator
$\Delta$  is just a Brownian motion with variance $2t$.
We would like to mention that condition \eqref{SQ5} guarantees the well-definedness of
$\int^t_0\bD_s(Z_s)\dif s$. 
\er

Our first result is the following weak/strong well-posedness about the second order mean-field SDE \eqref{MV1} or \eqref{MV11}.
\bt\label{Main0}
Let $\alpha\in(1,2]$,  $q_b\in(\frac\alpha{\alpha-1},\infty]$
and $\bbp_0,\bbrho_0,\bbrho_1\in[1,\infty]^2$ with $\bbb1\leq\tfrac1{\bbp_0}+\tfrac1{\bbrho_0}$ and $\bbrho_0\leq\bbrho_1$. 
Let $\bar\bbrho_1\in[1,\infty]^2$ with  $\frac1{\bbrho_1}+\frac1{\bar\bbrho_1}=\bbb1$ and $\bbp_1:=\bbp_0\wedge\bar\bbrho_1$. Define
$$
    \ \sA_i:=\bba\cdot(\tfrac d{\bbp_i}+\tfrac d{\bbrho_i}-\bbd), \ i=0,1,\quad
    \mZ_{\beta}:=\bB^{\beta}_{\bbp_0;\bba}\cap\bB^{\beta}_{\bbp_1;\bba},\ \beta\in\mR,
$$
where the anisotropic Besov space $\bB^{\beta}_{\bbp;\bba}:=\bB^{\beta,\infty}_{\bbp;\bba}$ is defined in Definition \ref{bs} below.
\begin{enumerate}[{\bf (A)}]
\item {\bf (Weak well-posedness)} Suppose that for  some $\beta_0\in(\frac{\alpha}{q_b}-\alpha,0)$ 
and $\mR\ni\beta_b\not=\sA_i-\beta_0$,
\begin{align}\label{AD067} 
0<\sA_0-\beta_0-\beta_b<\alpha-\tfrac\alpha{q_b}+\beta_0\wedge(-1),
\end{align}
and $b=b_0+b_1$, where for some $\kappa_b>0$,
\begin{align}\label{KB9}
\kappa_b:=\|b_0\|_{L^{q_b}(\mR_+; \bB^{\beta_b}_{\bbrho_0;\bba})}+\|b_1\|_{L^{q_b}(\mR_+; \bB^{\beta_b}_{\bbrho_1;\bba})}<\infty.
\end{align}
For any $\gamma_b>(\alpha-\frac\alpha{q_b}-1-\sA_1)\vee 0$, there is a constant $C_0=C_0(\alpha,d,\beta_0,\bbp_0, \beta_b,q_b,\gamma_b,\bbrho_0,\bbrho_1)>0$ 
such that for any $T\geq 1$ and initial distribution $\mu_0\in\mZ_{\beta_0}$ satisfying
\begin{align}\label{RQ22}
\|\mu_0\|_{\mZ_{\beta_0}}\leq C_0T^{-\frac{\gamma_b}\alpha}/\kappa_{b},
\end{align}
there exists a unique weak solution $Z$ to mean field SDE  \eqref{MV11} on $[0,T]$
in the sense of Definition \ref{Def1} so that $\DD\in\mL^q_T(\bC^\beta_\bba)$ for some $q\in(\frac\alpha{\alpha-1},\infty)$ and $\beta>0$, where $\DD$ is defined  by \eqref{FF}
and $\bC^\beta_\bba:=\bB^{\beta}_{\infty;\bba}$ is the anisotropic H\"older space.
\item {\bf (Strong well-posedness)}  Suppose that for  some $\beta_0\in(\frac{\alpha}{q_b}-\alpha,0)$ and $\mR\ni\beta_b\not=\sA_i-\beta_0$,
\begin{align}\label{JC08}
0<\sA_0-\beta_0-\beta_b<\tfrac32\alpha-\tfrac\alpha{q_b}-2,
\end{align}
and  $b=b_0+b_1$, where
\begin{align}\label{Conb}
b_i\in L^{q_b}(\mR_+;\bB^{s_0,s_1}_{\bbrho_i; x,\bba}\cap \bB^{\beta_b}_{\bbrho_i; \bba}),\ \ i=0,1,
\end{align}
and $s_0=1+\frac\alpha2$, $s_1=\frac\alpha2-1+\beta_b$ 
and the mixed Besov space $\bB^{s_0,s_1}_{\bbrho_i; x,\bba}$ is defined
in Definition \ref{bs} below. Let $\mu_0$ be the law of the initial random variable $Z_0$. Then, under \eqref{RQ22}, for any
$\gamma_b>(\alpha-\frac\alpha{q_b}-1-\sA_1)\vee 0$ and $T>0$, 
 there is a unique strong solution $Z$ to mean field SDE  \eqref{MV11} 
on $[0,T]$ in the sense of Definition \ref{Def1}
 so that $\DD\in\mL^q_T(\bC^\gamma_x\cap\bC^\beta_\bba)$ for some $q\in(\frac\alpha{\alpha-1},\infty)$ and $\gamma>\frac{1+\alpha/2}{1+\alpha}$, $\beta>1-\frac\alpha2$, where $\bC^\gamma_x$ is the H\"older space along the $x$-direction.
 \end{enumerate}
Moreover, we have:
\begin{enumerate}[(i)]
\item {\bf (Smoothness of density and short time singularity estimates)} For any $t\in(0,T]$,
$Z_t$ has a smooth density $f$ with short time regularity
$$
\sup_{t\in(0, T]}\left(t^{\frac{\beta}{\alpha}}\|f(t)\|_{\mZ_{\beta+\beta_0}}\right)<\infty, \quad \text{ for any $\beta\geq 0$}.
$$
\item {\bf (Large time decay estimate)} 
Suppose $\bba\cdot\tfrac{d}{\bbrho_1}>\alpha-\tfrac{\alpha}{q_b}-1$ and $\bbp_0\not=1$. 
Then we can take $\gamma_b=0$ in \eqref{RQ22}, and  there is a global weak solution to mean-field SDE \eqref{MV11} so that
for any $\beta\geq0$ and $i=0,1$,
$$
\sup_{t>1}\Big(t^{\frac{\bba\cdot(\bbd-d/\bbp_i)}{\alpha}}\|f(t)\|_{\bB^{\beta,1}_{\bbp_i;\bba}}\Big)<\infty.
$$
\item {\bf (Weak stability)}  Suppose 
that the assumptions in {\bf (A)} hold. 
Let $\wt b$ be another kernel satisfying \eqref{KB9} and $\wt \mu_0\in\mZ_{\beta_0}$ another initial distribution
that satisfies  \eqref{RQ22}.
Let $\wt f$ be the distribution density of the solution $\wt Z$ corresponding to $\wt b$ and $\wt\mu_0$. Then for any $\beta\geq 0$,
\begin{align}\label{AS629}
\sup_{t\in(0, T]}\left(t^{\frac{\beta}{\alpha}}\|f(t)-\wt f(t)\|_{\mZ_{\beta+\beta_0}}\right)\lesssim_{ C_\beta}\sum_{i=0,1}
\|b_i-\wt b_i\|_{\mL^{q_b}_T(\bB^{\beta_b}_{\bbrho_i;\bba})}+\|\mu_0-\wt \mu_0\|_{\mZ_{\beta_0}},
\end{align}
where $C_\beta=C_\beta(\alpha,d,\beta_0,\bbp_0, \beta_b,q_b,\bbrho_0,\bbrho_1, T,\kappa_b)>0$.
\item {\bf (Strong stability)}
Suppose that the assumptions in {\bf (B)} hold. 
Let $\wt b$ be another kernel satisfying \eqref{Conb} and $\wt \mu_0\in\mZ_{\beta_0}$ another initial distribution
that satisfies  \eqref{RQ22}. Let $\wt Z$ be the strong solution of SDE \eqref{MV11} corresponding to $\wt b$ and $\wt Z_0\sim\wt\mu_0$. Then
\begin{align}\label{AS689}
\mE\left(\sup_{t\in(0, T]}|Z_t-\wt Z_t|^2\right)\lesssim_C\mE|Z_0-\wt Z_0|^2+\sum_{i=0,1}
\|b_i-\wt b_i\|^2_{\mL^{q_b}_T(\bB^{\beta_b}_{\bbrho_i;\bba})}+
\|\mu_0-\wt \mu_0\|_{\mZ_{\beta_0}}^2,
\end{align}
where $C=C(\alpha,d,\beta_0,\bbp_0, \beta_b,q_b,\bbrho_0,\bbrho_1, T,\kappa_b)>0$.
\item {\bf (Global existence)} If $\bbp_0=\bbb1$ or $(\div_v b)^-=0$ and $\bbp_0=(p_0,p_0)$,
$\bbrho_1=(\varrho_1,\varrho_1)$, then the smallness condition \eqref{RQ22} can be dropped.
\end{enumerate}
\et

\br\rm
If $\bbrho_0=\bbrho_1$, then $\bbp_0=\bbp_1$ since $\bbp_0\leq\bar\bbrho_1$ by $\bbb1\leq\tfrac1{\bbp_0}+\tfrac1{\bbrho_0}$.
Moreover,
$(\div_v b)^-=0$ means that for any nonnegative $g\in C^\infty_c(\mR^{2d})$,
$$
\int_{\mR^{2d}}\<b(t,z),\nabla g(z)\>\dif z\le0,\quad \text{for a.e. $t\in\mR_+$}.
$$
\er
\br\rm
Let $\mu_0\in \cP(\mR^{2d})\cap\mZ_{\beta_0}$. Define for $m\in\mN$,
$$
\mu_0^m(\dif x,\dif v):=\mu_0(2^{-(1+\alpha)m}\dif x,2^{-m}\dif v).
$$
Then  $\mu^m_0\in \cP(\mR^{2d})\cap\mZ_{\beta_0}$ and by \eqref{SC8} below, we have
$$
\|\mu_0^m\|_{\bB^{\beta_0}_{\bbp_i;\bba}}\le \tfrac{2^{m(\bba\cdot(\frac d{\bbp_i}-\bbd)-\beta_0)}}{1-2^{\beta_0}}
\|\mu_0\|_{\bB^{\beta_0}_{\bbp_i;\bba}},\ \ i=0,1.
$$
From this, one sees that \eqref{RQ22} holds for $\mu_0^m$ as long as $0>\beta_0>\bba\cdot(\frac d{\bbp_i}-\bbd)$ and $m$ is large enough.

\er

\br\rm
The condition $\beta_0\in\{0,\sA_i-\beta_b,i=0,1\}$ are purely technical and can be dropped 
since we can choose different $\beta_0'<\beta_0$
so that \eqref{AD067} still holds. This condition comes from the following Besov estimate:
for $\bbp\leq\bbp'$ and $\beta'\not=\beta-\sA$ with $\sA:=\bba\cdot(\tfrac d{\bbp}-\tfrac d{\bbp'})$,
$$
\|P_tf\|_{\bB^{\beta',1}_{\bbp';\bba}}\lesssim_C(1\wedge t)^{-\frac{(\beta'-\beta+\sA)\vee0}\alpha}\|f\|_{\bB^{\beta,\infty}_{\bbp;\bba}},\ \ t>0,
$$
where $P_t$ is the kinetic semigroup (see Lemma \ref{Le215} below).
Note that this estimate is stronger than the usual estimate
$\|\nabla_v P_tf\|_{\bbp}\leq C(1\wedge t)^{-\frac1\alpha}\|f\|_{\bbp}$. Indeed,
$$
\|\nabla_v P_tf\|_{\bbp}\leq C\|P_tf\|_{\bB^{1,1}_{\bbp;\bba}}\leq C(1\wedge t)^{-\frac1\alpha}\|f\|_{\bB^{0,\infty}_{\bbp;\bba}}
\leq C(1\wedge t)^{-\frac1\alpha}\|f\|_{\bbp}.
$$
The above estimate plays a crucial role in our proofs.
\er

\br\rm
Let $b:\mR_+\times\mR^{2d}\times\mR^{2d}$ be a measurable function with $|b(t,z,z')|\leq K_t(z-z')$, where $K$ satisfies 
\eqref{KB9} with $-1<\beta_b<0$. If we replace $b*\mu$ by 
$$
(b\circledast\mu)(z):=\int_{\mR^{2d}}b(z,z')\mu(\dif z'),
$$
then the weak well-posedness still holds (see Remark \ref{Re48} below). In this case, it is of course not expected that \eqref{MV11} 
admits a smooth density except for the case where that $z\mapsto b(t,z,z')$ is smooth.
\er

\br\rm
The stability estimates \eqref{AS629} and \eqref{AS689} provide us the possibility of showing the propagation of chaos
for the cutoff kernel $b$. For instance, consider the kernel $K(x)=\nabla |x|^{s-d}$ and the cutoff kernel
$K_N(x):=\nabla (|x|\vee N^{-\delta})^{s-d}$, where $s\in(0,d)$ and $\delta>0$. For any $\varrho\in(\frac{d}{d-s},\infty]$ and $\beta<\frac dp-d+s$,
by \eqref{AW3} in Appendix A, one sees that
$$
\|K-K_N\|_{\bB^{\beta}_\varrho}\leq N^{-\delta(\frac d p-d-\beta+s)}. 
$$
Thus one can consider the regularized particle system \eqref{Par1} with $K=K_N$ (see \cite{HLP20}).
\er

Now we sketch the main ingredience of the proof of Theorem \ref{Main0}.
Let $P_t$ be the kinetic semigroup associated with $\Delta^{\alpha/2}_v-v\cdot\nabla_x$ (see \eqref{SCL11} for a definition).
First of all we need to solve FPE \eqref{FPE00} and obtain some precise asymptotic estimates in anisotropic Besov spaces.
More precisely, by Duhamel's principle,
$$
f(t)=P_t f_0-\cH^b_t(f,f),
$$
where 
$$
\cH^b_t(f,g)=\int^t_0 P_{t-s}\div_v((b_s*f_s)g_s)\dif s.
$$
The key point is to obtain a priori estimates for the bilinear form $(f,g)\mapsto\cH^b(f,g)$ in suitable spaces, and then use the contraction mapping theorem. 
To achieve this goal, we first establish Schauder estimates
for the non-homogenous kinetic equation $\p_t u=(\Delta^{\alpha/2}_v-v\cdot\nabla_x)u+f$ 
 in time-weighted spaces (see Theorem \ref{Th26} below).  To show the smoothness and large time decay estimate for the solution, 
 by the semigroup property of $P_t$ and shifting the time, the 
 following simple observation is important (see Lemma \ref{Le38} below):
 $$
 \cH^b_{t+s}(f,g)=P_s\cH^b_{t}(f,g)+\cH^{b^t}_{s}(f^t,g^t),\ \ s,t>0,
 $$
 where $f^t(\cdot,x,v)=f(t+\cdot,x,v)$.
 After establishing the solvability of \eqref{FPE00}, to prove
 the weak well-posedness of McKean-Vlasov SDE \eqref{MV11}, 
 we shall freeze the distribution-dependent drift along the smooth solution of the nonlinear Fokker-Planck equation, 
 and then utilize the uniqueness of linear Fokker-Planck equations. For the strong well-posedness, we use Zvonkin's transformation to kill
 the singular drift, where the freezed drift is H\"older continuous in the spatial variable, but has a singularity around the starting point 
 with respect to the time variable.

As a corollary, we also consider the following first order non-degenerate mean-field SDE 
\begin{align}\label{MV22}
\dif X_t=(b_t*\mu_{X_t})(X_t)\dif t+\dif L^{(\alpha)}_t,\ \alpha\in(1,2].
\end{align}
We have  the following weak/strong well-posedness result.
\bt\label{Main1}
Let $\alpha\in(1,2]$,  $q_b\in(\frac\alpha{\alpha-1},\infty]$
and $p_0,\varrho_0,\varrho_1\in[1,\infty]$ with $1\leq\tfrac1{p_0}+\tfrac1{\varrho_0}$ and $\varrho_0\leq \varrho_1$. 
Let $\bar\varrho_1\in[1,\infty]$ with  $\frac1{\varrho_1}+\frac1{\bar\varrho_1}=1$ and 
$p_1:=p_0\wedge\bar\varrho_1$. Define
$$
    \ \sA_i:=\tfrac d{p_i}+\tfrac d{\varrho_i}-d, \ i=0,1,\quad
    \mZ_{\beta}:=\bB^{\beta}_{p_0}\cap\bB^{\beta}_{p_1},\ \beta\in\mR,
$$
where $\bB^{\beta}_{p}=\bB^{\beta,\infty}_p$ is the isotropic Besov space defined in Definition \ref{bs} below.
\begin{enumerate}[{\bf (A)}]
\item {\bf (Weak well-posedness)} Suppose that for  some $\beta_0\in(\frac{\alpha}{q_b}-\alpha,0)$ and  $\mR\ni\beta_b\not=\sA_i-\beta_0$,
\begin{align}\label{AD0167} 
0<\sA_0-\beta_0-\beta_b<\alpha-\tfrac\alpha{q_b}+\beta_0\wedge(-1),
\end{align}
and $b=b_0+b_1$, where
\begin{align}\label{KB19}
\kappa_b:=\|b_0\|_{L^{q_b}(\mR_+; \bB^{\beta_b}_{\varrho_0})}+\|b_1\|_{L^{q_b}(\mR_+; \bB^{\beta_b}_{\varrho_1})}<\infty.
\end{align}
For any $\gamma_b>(\alpha-\frac\alpha{q_b}-1-\sA_1)\vee 0$, there is a constant $C_0=(\alpha,d,\beta_0, p_0, \beta_b,q_b, \varrho_0, \varrho_1)>0$ 
such that for any $T\geq 1$ and initial distribution $\mu_0\in\mZ_{\beta_0}$ satisfying
\begin{align}\label{RQ202}
\|\mu_0\|_{\mZ_{\beta_0}}\leq C_0T^{-\frac{\gamma_b}\alpha}/\kappa_{b},
\end{align}
there exists a unique weak solution $X$ to mean field SDE  \eqref{MV22} on $[0,T]$ in the sense of Definition \ref{Def1} so that $\DD\in\mL^q_T(\bC^\beta)$ for some $q\in(\frac\alpha{\alpha-1},\infty)$ and $\beta>0$, where $\DD_t(x):=(b_t*\mu_{X_t})(x)$
and $\bC^\beta:=\bB^{\beta}_{\infty}$ is the isotropic H\"older space. 
\item {\bf (Strong well-posedness)}  Suppose that for  some $\beta_0\in(\frac{\alpha}{q_b}-\alpha,0)$ and  $\mR\ni\beta_b\not=\sA_i-\beta_0$,
\begin{align}\label{JC008}
0<\sA_0-\beta_0-\beta_b<\tfrac32\alpha-\tfrac\alpha{q_b}-2.
\end{align}
For any $\gamma_b>(\alpha-\frac\alpha{q_b}-1-\sA_1)\vee 0$ and $T\geq 1$, 
under \eqref{KB19} and \eqref{RQ202}, 
there is a unique strong solution $X$ to mean field SDE  \eqref{MV22} on $[0,T]$ in the sense of Definition \ref{Def1}
so that $\DD\in\mL^q_T(\bC^\beta)$ for some $q\in(\frac\alpha{\alpha-1},\infty)$ and $\beta>1-\frac\alpha2$, where $\DD_t(x):=(b_t*\mu_{X_t})(x)$
\end{enumerate}
Moreover, we have the following:
\begin{enumerate}[(i)]
\item {\bf (Smoothness of density and short time singularity estimates)} For any $t\in(0,T]$,
$Z_t$ has a smooth density $\rho(t)$ with regularity that for any $\beta\geq 0$,
$$
\sup_{t\in(0, T]}\left(t^{\frac{\beta}{\alpha}}\|\rho(t)\|_{\mZ_{\beta+\beta_0}}\right)<\infty.
$$
\item {\bf (Large time decay estimate)} 
Suppose $\tfrac{d}{\varrho_1}>\alpha-\tfrac{\alpha}{q_b}-1$ and $p_0\not=1$. Then we can take $\gamma_b=0$ in \eqref{RQ202}, and  there is a global weak solution to mean-field SDE \eqref{MV22} so that
for any $\beta\geq0$ and $i=0,1$,
$$
\sup_{t>1}\Big(t^{\frac{d-d/p_i}{\alpha}}\|\rho(t)\|_{\bB^{\beta,1}_{p_i}}\Big)<\infty.
$$

\item {\bf (Two-sided estimates)} For any $T>0$, there are constants $C_1,C_2,C_3>0$ 
such that for all $t\in(0,T)\times\mR^d$,
$$
\frac{C^{-1}_1}{t^{d/2}}\int_{\mR^d} \e^{-\frac{|x-y|^2}{C_2t}}\mu_0(\dif x)\leq \rho(t,y)\leq 
\frac{C_1}{t^{d/2}}\int_{\mR^d} \e^{-\frac{C_2|x-y|^2}{t}}\mu_0(\dif x),\ \ \alpha=2,
$$
and
$$
\int_{\mR^d} \frac{C_3^{-1}t\mu_0(\dif x)}{(t^{1/\alpha}+|x-y|)^{d+\alpha}}\leq \rho(t,y)\leq \int_{\mR^d} \frac{C_3t\mu_0(\dif x)}{(t^{1/\alpha}+|x-y|)^{d+\alpha}},\ \ \alpha\in(1,2).
$$
\item {\bf (Global existence)} If $p_0=1$ or $(\div b)^-=0$, then we can drop the smallness condition \eqref{RQ202}.
\end{enumerate}
\et

\br\rm\label{Re10}
When $p_0=1$, conditions \eqref{AD0167} and \eqref{JC008} reduce to
\begin{align}\label{dF1}
\tfrac{d}{\varrho_0}+\tfrac\alpha{q_b}-\beta_b<\alpha-1,\ \ \tfrac{d}{\varrho_0}+\tfrac\alpha{q_b}-\beta_b<\tfrac{3\alpha}2-2,
\end{align}
which coincide with {\bf (C0)} and {\bf (C0$_S$)} in \cite{CJM22}. In this case, under the same Besov regularity 
on $b$, {\bf (A)} and {\bf (B)} were established by Chaudru de Raynal, Jabir and Menozzi in \cite{CJM22}. We assume $p_0\ge1$, since we were interested to cover more examples from physics for which typical $p_0$ has to be taken in $(1,\infty]$.
For example, if $b(x)=\nabla |x|^{s-d}$ for some $s\in(2-\alpha,d)$ and $\mu_0\in\cap_{\eps>0}\bB^{-\eps}_{p_0}$, where $p_0>\frac{d}{s+\alpha-2}$, then the assumptions in Theorem \ref{Main1} hold by \eqref{HH2} in Appendix A below.
Moreover, it is well known that for each $x\in\mR^d$, the Dirac measure $\delta_x\in\bB^{0}_1$. 
Thus one sees that for $b=(\delta_{x_0},\cdots,\delta_{x_d})$ and $\mu_0\in\cap_{\eps>0}\bB^{-\eps}_{p_0}$, where $p_0>\frac{d}{\alpha-1}$, the assumptions in Theorem \ref{Main1} hold. In this case, SDE \eqref{MV22} is a density-dependent or Nemyskii's type SDEs.
\er

\br\label{Rk111H}
\rm In this remark we make simple regularity analysis to illustrate the sharpness of condition \eqref{AD0167} for $\beta_0\in[-1,0)$, 
which belongs to
the subcritical regime. If the strict inequality in \eqref{AD0167} is an equality, then it is called critical.
First of all, we consider the following  SDE with linear drift
\begin{align}\label{SDE1}
\dif X_t=H_t(X_t)\dif t+\dif L^{(\alpha)}_t,\ \ X_0\sim\mu_0,
\end{align}
where $H\in L^q_tL^p_x$ for some $q,p\in[1,\infty]$. 
By simple scaling analysis, it is well known that for SDE \eqref{SDE1} or the associated Kolmogorov equations to be uniquely solvable, 
even though $\mu_0$ is smooth, we have to require (see \cite{XZ21})
\begin{align}\label{AZ90}
\tfrac d p+\tfrac\alpha q<\alpha-1.
\end{align}
In other words, for linear SDE \eqref{SDE1}, the regularity of initial values does not affect the singularity of $H$.
Now we consider the following first-order DDSDE with convolution drift:
\begin{align}\label{AZ66}
    \dif X_t=b_t*\mu_t(X_t)\dif t+\dif L^{(\alpha)}_t,\ \ X_0\sim\mu_0,
\end{align}
where $\mu_t$ is the law of $X_t$. 
%Let $\dot\bB^\beta_p=\dot\bB^{\beta,\infty}_p$ be the usual homogeneous Besov space. 
Suppose that for some $\beta_0,\beta_b\in\mR$ and $q_b,p_b,p_0\in[1,\infty]$,
\begin{align}\label{AZ98}
b\in L^{q_b}(\mR_+; \bB^{\beta_b}_{p_b}),\ \ \mu_0\in\bB^{\beta_0}_{p_0}.
\end{align}
We want to seek reasonable conditions on $q_b,\beta_b,p_b$ and $\beta_0,p_0$ so that SDE \eqref{AZ66} is uniquely solvable.
By It\^o's formula, $\mu_t$ solves the Fokker-Planck equation
in the distributional sense:
$$
\p_t\mu=\Delta^{\frac\alpha 2}\mu-\div ((b*\mu)\mu).
$$
Let $P_t$ be the semigroup associated with $\Delta^{\frac\alpha 2}$.
By Duhamel's formula,
\begin{align}\label{JC101}
\mu_t=P_t\mu_0-\int^t_0P_{t-s}\div ((b_s*\mu_s)\mu_s)\dif s.
\end{align}
For any $p_1\ge p_0$ and $\gamma\ge \frac d{p_1}-\frac d {p_0}$, by the regularization effect of $P_t$ (see \eqref{AD0446} below), we have
$$
\|P_t\mu_0\|_{\bB^{\gamma+\beta_0}_{p_1}} \leq 
C(1\wedge t)^{-(\gamma+\frac{d}{p_0}-\frac{d}{p_1})/\alpha}\|\mu_0\|_{\bB^{\beta_0}_{p_0}}.
$$
In particular,  for any $T>0$,  without considering the second term in \eqref{JC101},
\begin{align}\label{JC10}
\mbox{$\mu$ at most belongs to $L^{q_1}_T(\bB^{\gamma+\beta_0}_{p_1})$ for any $(q_1,p_1,\gamma)\in\sI$,}
\end{align}
where
\begin{align}\label{JC19}
\sI:=\Big\{(q_1,p_1,\gamma): p_1\geq p_0,\ 0\leq\gamma+\tfrac{d}{p_0}-\tfrac{d}{p_1}<\tfrac{\alpha}{q_1}\Big\}.
\end{align}
Moreover, under \eqref{AZ98}, for any $(q_1,p_1,\gamma)\in\sI$, by Young's convolution inequality,
$$
\DD_t(x):=b_t*\mu_t(x)\in L^q_T(\bB^\beta_p),
$$
where $(\beta,q,p)$ is defined by
\begin{align}\label{AZ99}
\beta=\beta_b+\gamma+\beta_0,\quad\tfrac{1}{q}=\tfrac{1}{q_b}+\tfrac{1}{q_1},\quad 1+\tfrac1{p}=\tfrac1{p_b}+\tfrac1{p_1}.
\end{align}
After this analysis, 
in order to make sense for $\int^T_0H_s(X_s)\dif s$ and to solve SDE \eqref{SDE1}, 
it is reasonable to assume $\DD_t(x)\in L^q_T(L^p)$
for some $q,p$ satisfying \eqref{AZ90}. Thus by \eqref{AZ99} with $\beta>0$, \eqref{JC19} and \eqref{AZ90}, we have
$$
-\beta_0-\beta_b<\gamma<\tfrac{\alpha}{q_1}+\tfrac{d}{p_1}-\tfrac{d}{p_0}
=\tfrac\alpha q-\tfrac\alpha{q_b}+d+\tfrac{d}{p}-\tfrac{d}{p_b}-\tfrac{d}{p_0}
<\alpha-1-\tfrac\alpha{q_b}+d-\tfrac{d}{p_b}-\tfrac{d}{p_0}.
$$
This is the case of \eqref{AD0167} for $\beta_0\in[-1,0)$.
However, for $\beta_0\in(\frac\alpha{q_b}-\alpha,-1)$, at this moment, we do not know whether the following condition that comes from \eqref{AD0167} 
$$
2\beta_0-\beta_b<\alpha-\tfrac{\alpha}{q_b}+d-\tfrac{d}{p_b}-\tfrac{d}{p_0}
$$
is sharp. Of course, for a smaller differentiability index $\beta_0$, the initial value could be more singular.
We will study this in a forthcoming work.
\iffalse
In addition, in order to make the product $(b*\mu)\mu=H\mu$ of two distributions well-defined, 
by \eqref{JC10}, we need to at least assume
$\mu\in L^{\tilde q_1}_T(\bB^{\tilde \gamma+\beta_0}_{p_0})$ for some $(\tilde q_1,p_0,\tilde\gamma)\in\sI$,
$$
\tfrac1{\tilde q_1}+\tfrac1q=1,\ \ -\beta_0<\tilde\gamma<\tfrac\alpha{\tilde q_1},
$$
which together with \eqref{AZ99} yields
$$
-\beta_0<\tfrac{\alpha}{\tilde q_1}
=\alpha-\tfrac\alpha{q_b}-\tfrac\alpha{q_1}\stackrel{\eqref{AZ90}}{<}\alpha-\tfrac\alpha{q_b}
+\tfrac{d}{p_1}-\tfrac{d}{p_0}+\beta_b+\beta_0.
$$
By $\frac1{p_1}+\frac1{p_b}=1$ we obtain
$$
2\beta_0-\beta_b<\alpha-\tfrac{\alpha}{q_b}+d-\tfrac{d}{p_b}-\tfrac{d}{p_0}.
$$
This together with \eqref{AZ94}  implies \eqref{AD0167}  with $\varrho_0=p_b$. 
\fi
\er

\subsection{Related works}

The study of kinetic equations has a long history. These are used to describe the evolution of the distribution density
of a single particle that is characterized by Newton's second law in gas theory, aerodynamics, plasma physics and 
non-equilibrium statistical physics (see \cite{C69}, \cite{Ja14},\cite{ Go16}).
In 1872, Boltzmann put forward the fundamental kinetic equation of gas theory, which is a non-linear integro-differential equation, 
now named  Boltzmann equation, that is,
$$
\p_t f+v\cdot\nabla_x f=Q(f,f),
$$
where $Q(f,f)$ is the collision operator (see \cite{C88}).
It describes the motion of the molecules as a certain random process determined by collisions between pairs of molecules (see \cite{C88}). 

The two-sided estimates of the heat kernel of kinetic operator $\cK:=\Delta_v-v\cdot\nabla_x+b\cdot\nabla_x$ was studied in \cite{DM10} for Lipschitz drift $b$ and in \cite{CMPZ22} for H\"older drift $b$.
The $L^p$-regularity estimate of kinetic operators $\cK$ was shown in Bouchut \cite{Bo02}. 
For the fractional kinetic operator $\cK^{(\alpha)}:=\Delta^{\alpha/2}_v-v\cdot\nabla_x+b\cdot\nabla_x$,
the $L^p$-maximal regularity estimates were established in \cite{CZ18} and \cite{HMP19}
(see also \cite{DY22}), the Schauder estimates were obtained in \cite{HWZ20} and \cite{IS21}. In fact, our proofs about the solvability of FPE \eqref{FPE00} 
strongly depend on the method developed in \cite{HWZ20}.

Consider now the following McKean-Vlasov or distribution-dependent SDE
\begin{align}\label{SJ1}
\dif X_t =b(t,X_t,\mu_{X_t})\dif t +\sigma(t,X_t,\mu_{X_t})\dif W_t.
\end{align}
When $b$ and $\sigma$ are Lipschitz continuous in $x$ and in the unknown distribution $\mu$ with respect to the Euclidean distance and
Wasserstein metric, it is by now standard to show the existence and uniqueness of strong solutions to DDSDE \eqref{SJ1} (cf. \cite{Szn91} and \cite{CD21}).
Moreover, under some one-sided Lipschitz assumptions, Wang  \cite{Wa18}
showed the strong well-posedness and some functional inequalities for DDSDE \eqref{SJ1}. For coefficients $b,\sigma$ only measurable in the measure variable, hence also containing the case of distribution density dependent (also called Nemytskii-type) coefficients, existence  and uniqueness of weak solutions was proved in a series of papers (see \cite{BR18, BR20, BR21, BR23, BR22}) and in the case of multiplicative L\'evy noise \cite{RXZ20}.  
In \cite{HSS21}, Hammersley, Siska and Szpruch proved the existence of weak solutions to SDE \eqref{SJ1} on a domain 
 with continuous and unbounded coefficients under Lyapunov-type conditions. Moreover, uniqueness is also obtained under some functional Lyapunov conditions. 
 When the diffusion matrix is uniformly non-degenerate and are only measurable and of at most linear growth, by using the classical Krylov estimates, Mishura and Veretennikov \cite{MV16} showed the existence of weak solutions. 
 The uniqueness is also proved when the diffusion coefficient does not depend on the unknown distribution $\mu$
 and is Lipschitz continuous in $x$ and $b$ is Lipschitz continuous with respect to  $\mu$.  
 When $b(x,\mu)=\int_{\mR^d}\bar b(x,y)\mu(\dif x)$ with $|\bar b(x,y)|\leq h(x-y)$, where $h\in L^p(\mR^d)$ for some $p>d$,
and $\sigma(x,\mu)=\sigma(x)$ does not depend on the distribution and is uniformly elliptic and H\"older continuous,
the well-posedness of DDSDE \eqref{SJ1} was shown in \cite{RZ21}.
More recently, Chaudru de Raynal, Jabir and  Menozzi \cite{CJM22} studied the well-posedness to the following DDSDE
with convolution drift
$$
\dif X_t=(b_t*\mu_{X_t})(X_t)\dif t+\dif L^{(\alpha)}_t.
$$
Under \eqref{dF1}, for $b\in L^{q_b}_T(\bB^{\beta_b}_{p_b})$ and any initial distribution $\mu_0$, 
by smoothing $b$ and a duality method, they showed the weak/strong 
well-posedness to the above SDE. {Therein they also obtained some H\"older regularity of the solution. Later, Han in \cite{Ha22} improved it to the smooth solution under some conditions of $b$.} 

Next consider the following linear second order kinetic SDE:
\begin{align}\label{Kn1}
\ddot X_t=H_t(X_t,\dot X_t)+\dot L^{(\alpha)}_t,\ \ (X_0,\dot X_0)\sim \mu_0,
\end{align}
where $H$ is only a function of $(t,x,v)$. When $\alpha=2$ and $H$ is $2/3+\eps$-H\"older continuous in $x$ and 
$\eps$-H\"older continuous in $v$, the strong well-posedness was established in \cite{Ch17} and \cite{WZ16}.
The $L^p$-type drift was considered in \cite{Zh18}. In \cite{CM22}, Chaudru de Raynal and Menozzi showed the weak well-posedness when
$H\in L^q_t L^p_{x,v}$, where $\frac{4d}{p}+\frac 2q<1$. For $\alpha\in(1,2)$, the strong well-posedness for kinetic SDE \eqref{Kn1}
was established in \cite{HWZ20} when $H\in L^\infty_t(C^{\gamma,\beta}_{x,v})$, where $\gamma\in(\frac{2+\alpha}{2(1+\alpha)},1)$ and $\beta>1-\frac\alpha 2$. It seems that there are few results about the well-posedness of kinetic DDSDEs with singular interaction kernels. The aim of this paper is to fill this gap.

After the contents of this paper was presented in a talk of ``The 
Non-local operators, probability and singularities''
online seminar, Stephane Menozzi informed us in a very kind e-mail
that he with P.E. Chaudru de Raynal and J.F. Jabir have a paper in
preparation with related results and sent us a current (draft) version of it
(\cite{CJM23}). We would like to thank Stephane very much for
this. After a detailed discussion with Stephane, we emphasize the following differences between our and their work.
\begin{itemize}
\item The first and foremost difference concerns the motivation and the models in applications. In their paper, they continue the research from their very interesting previous work \cite{CJM22} and consider non-degenerate McKean-Vlasov equations where the initial data is a real function, i.e. $\beta_0\ge0$, while in our present paper, we consider degenerate kinetic equations with singular kernels and measure initial data, which can be deduced naturally from the second order $N$-particle system \eqref{Par1}. As a consequence, our results can be applied to non-degenerate equations (see Theorem \ref{Main1}).  

\item The key condition \eqref{AD0167} is not the same {as ({\bf C$_1$}) in \cite{CJM23}}. More precisely, there is only one $\beta_b$ in our condition \eqref{AD0167}. But the condition {({\bf C$_1$})} in \cite{CJM23} includes $2\beta_b$. Moreover, in \cite{CJM23} $\beta_0\ge0$ is required, whereas $\beta_0<0$ is allowed in our paper. It should be noted that our condition is sharp in the scaling sense (see Remark \ref{Rk111H}). 
Furthermore, in the PDE part, we cover the critical case in Theorem \ref{Main}, where the key point is that in our Schauder estimate Theorem \ref{Th26} the norms in the left and right hand sides are $\bB^{\beta',1}_{p'}$ and $\bB^{\beta,\infty}_{p}$, respectively.

\item Another important difference is in the methods of showing the well-posedness of the mean-field SDEs. In \cite{CJM23}, the authors mollify the singular interaction kernel and show some uniform estimate for the solution to the nonlinear Fokker-Planck equation. Then, by the weak convergence method
they obtain the existence of a solution in the sense of a martingale problem. 
In our work, we use the contraction mapping theorem to show the well-posedness of the nonlinear Fokker-Planck equation and obtain
the well-posedness of the linearized SDE. Then, based on the uniqueness of the linear Fokker-Planck equation,  
the solution to the linearized SDE is just the solution to the mean-field SDE (see Section \ref{Sec4.2H} for more details).   
\item
For concrete physical models like the Navier-Stokes equations, in \cite{CJM23}, in order to avoid the non-integrability at infinity 
the authors cut off the kernel {by an indicator function of a bounded domain}, whereas, in our work, we consider the kernel being the sum of two singular kernels so that we can cover 
the physical interaction kernel in the whole space.
\item In the current work, we also pay much attention to the smoothness, stability, and large-time behavior of solutions to nonlinear Fokker-Planck equations, which are not considered in \cite{CJM23}. In fact, there are many techniques developed to study these aspects. For instance, we use a time shift trick to obtain the smoothness and large time decay estimates of solutions to nonlinear Fokker-Planck equations. Moreover, by choosing suitable parameters, based on an interpolation technique, we can drop the smallness assumption of $\|u_0\|_{\cM}$ in the PDE's literature, 
which is important to study the related mean-field SDEs, see Remark \ref{Rk53H} below.   
\item Finally, there are also some differences in the techniques. For instance, for the treatment of the convolution $b*u$,
it is regarded as an inner product $\<b(x-\cdot),u(\cdot)\>$ in \cite{CJM23}. 
In our paper, we use Young's inequality to estimate it. Moreover, we work with the kinetic operator and in
anisotropic Besov spaces, which are more difficult to handle, since for instance $\cR_j^aP_t\ne P_t \cR_j^a$. To obtain the
well-posedness of the equations with large initial data locally in time in the critical case, we show an asymptotic estimate for small time
for the semigroup $P_t$, see Lemma \ref{Lem315H}.
\end{itemize} 

Summarizing the above, in our opinion \cite{CJM23} is a very valuable contribution to the subject containing various interesting and different focal points in comparison to our paper.

\subsection{Structure of the paper}
This paper is organized as follows. In Section 2 we introduce the anisotropic Besov spaces as well as some of their important properties.
Most of the results are well known and can be found in \cite{Tri06}. Moreover, we also prepare the basic Schauder estimates
for the model kinetic equation in time-weighted spaces (see Theorem \ref{Th26}). 

In Section 3, we study the global and local well-posedness for nonlinear kinetic Fokker-Planck equation \eqref{FPE00}
 with small or large initial data.
Moreover, we also investigate the smoothness, stability and large time decay estimates of the solution. 
More precisely, in Subsection 3.1, we prove our main result Theorem \ref{Main} about the global well-posedness for small initial values
for the nonlinear kinetic Fokker-Planck equation \eqref{FPE00}. In Lemma \ref{lem240}, we establish a priori bilinear 
estimate for $\cH^b_t(f,g)$, where Young's inequality for convolutions provides all the regularity induced by the kernel. 
In particular, the dependence on parameters is calculated carefully.
Then in Lemma  \ref{Le38}, we show the core estimate in the solution space.
%Especially, explicit conditions about the initial conditions and the kernel $b$ are given in {\bf (H$_0$)}.
In Subsection 3.2, we consider large time decay estimates.
In Subsection 3.3, we show the local well-posedness for large initial values in Theorem \ref{Cor35}. 
In Subsection 3.4, we prove two
results about the global well-posedness for \eqref{FPE00} for arbitrary large initial value.
The corresponding results for non-degenerate fractional Fokker-Planck equations are also presented.

In Section 4 we devote to the study of mean-field kinetic SDEs and prove Theorems \ref{Main0} and \ref{Main1}.
In Subsection 4.1, we first show the weak and strong well-posedness for linear kinetic SDEs with H\"older drifts and
driven by $\alpha$-stable processes by Zvonkin's transformation. In particular, for strong well-posedness, we need to 
assume the drift to be in $L^q_T(\bC^\frac{\alpha+\beta}{1+\alpha}_x\cap \bC^\beta_v)$, where 
 $q\in(\frac\alpha{\alpha-1},\infty]$ and $\beta>1-\frac\alpha 2$.
Note that for $q=\infty$, this has already been obtained in \cite{HWZ20}. 
In Subsection 4.2, we prove our main theorems by freezing the drifts.

In Section 5 we study the well-posedness of fractional Vlasov-Poisson-Fokker-Planck equation \eqref{VP1} 
with measures as initial data and the associated kinetic McKean-Vlasov SDEs. We also obtain an optimal large time decay estimate when the initial data are small. 
For Wiener noise in three dimensional case, when the initial data have finite second order moment,
using well-known results, we prove global existence of smooth solutions for large initial values.
In the fractional case, the global well-posedness of  \eqref{VP1} for large initial values seems to be open since there is no finite second order moment (see \eqref{20230303}), {which is crucial to show global well-posedness in \cite{Bou93, CS97, OS00}}.

In Section 6, we apply our main results to the vorticity formulation of fractional Navier-Stokes and subcritical SQG equations.
For two (and three) dimensional fractional Navier-Stokes equation, we show the global existence of smooth solutions when the initial vorticity is in $\bB^{-(1+\alpha)/2+}_\infty$ (and small). Moreover, in the two dimensional case,
%when the initial vorticity is a probability measure, 
we show the weak and strong well-posedness for the associated mean-field SDEs.
For the three dimensional case, we also obtain an optimal large time decay estimate. To the best of our knowledge, they are new, at least for fractional Navier-Stokes equations.

In Section 7 we present an application to fractional kinetic porous medium equations with viscosity and show the global existence 
of smooth solutions for large initial values.

We conclude this introduction by introducing the following convention: 
Throughout this paper, we use $C$ with or without subscripts to denote constants, whose values
may change from line to line. We also use $:=$ to indicate a definition and set
$$
a\wedge b:=\max(a,b),\ \ a\vee b:=\min(a,b),\ \ a^+:=0\vee a.
$$ 
By $A\lesssim_C B$ and $A\asymp_C B$
or simply $A\lesssim B$ and $A\asymp B$, we respectively mean that for some constant $C\geq 1$,
$$
A\leq C B,\ \ C^{-1} B\leq A\leq CB.
$$
Below we collect some frequently used notations for the readers' convenience.
\begin{itemize}
\item $\cP(\mR^N):=$ set of all probability measures on $\mR^N$.
\item $\cM(\mR^N):=$ set of all finite Radon measures on $\mR^N$.
\item $\|\cdot\|_{\cM}:=$ total variation norm.
\item $\bB^{s,q}_{\bbp;\bba}:=$ anisotropic Besov space.
\item $\bB^{s,q}_{\bbp}:=$ isotropic Besov space.
\item $\bC^s_{\bba}:=$ anisotropic H\"older-Zygmund space.
\item $\cR^\bba_j:=$ anisotropic block operator.
\item $P_t:=$ semigroup of $\Delta^{\alpha/2}_v-v\cdot\nabla_x$ or $\Delta^{\alpha/2}$ (see \eqref{SCL11} below).
\item $\mX_\beta,\mY_\beta,\mZ_\beta:=$ {some subspaces of Besov} spaces  (see \eqref{MY} below).
\item $\mS_\beta:=$ space of solutions (see \eqref{AS59} below).
\end{itemize}

 \section{Anisotropic Besov space and kinetic semigroup}
 
 In this section we first recall the definition of anisotropic Besov space with mixed integrability indices as well as its 
 basic properties (see \cite{Tri06} and \cite{ZZ21}).
Then we study the fractional kinetic semigroup and 
prove time-weighted Schauder's type estimates in anisotropic Besov spaces
for fractional kinetic operators. In order to study the small time singularity and large time decay estimates, 
we introduce a time-weighted $L^p$-space by considering different weights for small time and large time.

\subsection{Anisotropic Besov spaces}

For an $L^1$-integrable function $f$ in $\mR^{2d}$, let $\hat f$ be the Fourier transform of $f$ defined by
$$
\hat f(\xi):=(2\pi)^{-d}\int_{\mR^{2d}} \e^{-{\rm i}\xi\cdot z}f(z)\dif z, \quad\xi\in\mR^{2d},
$$
and $\check f$ the Fourier  inverse transform of $f$ defined by
$$
\check f(z):=(2\pi)^{-d}\int_{\mR^{2d}} \e^{{\rm i}\xi\cdot z}f(\xi)\dif\xi, \quad z\in\mR^{2d}.
$$
Recall $\bba=(1+\alpha,1)$. For $z=(x,v)$ and $z'=(x',v')$ in $\mR^{2d}$, we introduce the anisotropic distance 
$$
|z-z'|_\bba:=|x- x'|^{1/(1+\alpha)}+|v-v'|.
$$
Note that $z\mapsto|z|_{\bba}$ is not smooth.
For $r>0$ and $z\in\mR^{2d}$, we also introduce the ball centered at $z$ and with radius $r$ with respect to the above distance
as follows:
$$
B^\bba_r(z):=\{z'\in\mR^{2d}:|z'-z|_\bba\leq r\},\ \ B^\bba_r:=B^\bba_r(0).
$$
Let $\chi^\bba_0$ be  a symmetric $C^{\infty}$-function  on $\mR^{2d}$ with
$$
\chi^\bba_0(\xi)=1\ \mathrm{for}\ \xi\in B^\bba_1\ \mathrm{and}\ \chi^\bba_0(\xi)=0\ \mathrm{for}\ \xi\notin B^\bba_2.
$$
We define
$$
\phi^\bba_j(\xi):=
\left\{
\begin{aligned}
&\chi^\bba_0(2^{-j\bba}\xi)-\chi^\bba_0(2^{-(j-1)\bba}\xi),\ \ &j\geq 1,\\
&\chi^\bba_0(\xi),\ \ &j=0,
\end{aligned}
\right.
$$
where for $s\in\mR$ and $\xi=(\xi_1,\xi_2)$,
$$
2^{s\bba }\xi=(2^{s(1+\alpha)}\xi_1, 2^{s}\xi_2).
$$
Note that
\begin{align}\label{Cx8}
{\rm supp}(\phi^\bba_j)\subset\big\{\xi: 2^{j-1}\leq|\xi|_\bba\leq 2^{j+1}\big\},\ j\geq 1,\ {\rm supp}(\phi^\bba_0)\subset B^\bba_2,
\end{align}
and
\begin{align}\label{AA13}
\sum_{j\geq 0}\phi^\bba_j(\xi)=1,\ \ \forall\xi\in\mR^{2d}.
\end{align}

Let $\cS$ be the space of all Schwartz functions on $\mR^{2d}$ and $\cS'$ the dual space of $\cS$, called the tempered distribution space.
For given $j\geq 0$, the  dyadic anisotropic block operator  $\mathcal{R}^\bba_j$ is defined on $\cS'$ by
\begin{align}\label{Ph0}
\mathcal{R}^\bba_jf(z):=(\phi^\bba_j\hat{f})\check{\ }(z)=\check{\phi}^\bba_j*f(z),
\end{align}
where the convolution is understood in the distributional sense and by scaling,
\begin{align}\label{SX4}
\check{\phi}^\bba_j(z)=2^{(j-1)(2+\alpha)d}\check{\phi}^\bba_1(2^{(j-1)\bba}z),\ \ j\geq 1.
\end{align}
For $j\in\mN$, by definition it is easy to see that
\begin{align}\label{KJ2}
\cR^\bba_j=\cR^\bba_j\widetilde\cR^\bba_j,\ \mbox{ where }\ \widetilde\cR^\bba_j:=\sum_{|i-j|\leq 2}\cR^\bba_i,
\end{align}
where we used the convention that $\cR^\bba_i:=0$ for $i<0$. Moreover,
by the symmetry of $\phi^\bba_j$,
$$
\<\cR^\bba_j f,g\>=\< f,\cR^\bba_jg\>,\ \ f\in\cS', \ g\in\cS.
$$
Similarly, we can define the isotropic block operator $\cR_jf=\check\phi_j*f$ in $\mR^d$, where
\begin{align}\label{Cx9}
{\rm supp}(\phi_j)\subset\big\{\xi: 2^{j-1}\leq|\xi|\leq 2^{j+1}\big\},\ j\geq 1,\ {\rm supp}(\phi_0)\subset B_2.
\end{align}

The following Bernstein inequality is standard and proven in \cite{ZZ21}. We omit the details.
\bl[Bernstein's inequality]\label{BI00}
For any $\bbk:=(k_1,k_2)\in\mN^2_0$, $\bbp\le\bbp'\in[1,\infty]^2$, there is a constant $C=C(\bbk,\bbp,\bbp',\alpha,d)>0$ such that for all $j\geq 0$,
\begin{align}\label{Ber}
\|\nabla^{k_1}_x\nabla^{k_2}_v\cR^\bba_j f\|_{\bbp'}\lesssim_C 2^{j \bba\cdot(\bbk+\frac{d}{\bbp}-\frac{d}{\bbp'})}\|\cR^\bba_j  f\|_{\bbp},
\end{align}
where $\|\cdot\|_\bbp$ is defined in \eqref{LP1}.
\el

Now we introduce the following anisotropic Besov spaces (cf. \cite[Chapter 5]{Tri06}).
\begin{definition}\label{bs}
Let $s\in\mR$, $q\in[1,\infty]$ and $\bbp\in[1,\infty]^2$. The  anisotropic Besov space is defined by
$$
\mathbf{B}^{s,q}_{\bbp;\bba}:=\left\{f\in \cS': \|f\|_{\mathbf{B}^{s,q}_{\bbp;\bba}}
:= \left(\sum_{j\geq0}\big(2^{ js}\|\cR^\bba_{j}f\|_{\bbp}\big)^q\right)^{1/q}<\infty\right\},
$$
where $\|\cdot\|_\bbp$ is defined in \eqref{LP1}.
Let $(s_0,s_1)\in\mR^2$. The mixed Besov space is defined by
$$
\bB^{s_0,s_1}_{\bbp;x,\bba}:=\left\{f\in \cS': \|f\|_{\mathbf{B}^{s_0,s_1}_{\bbp;x,\bba}}
:= \sup_{k,j\geq 0}2^{\frac{ks_0}{1+\alpha}}2^{ js_1}\|\cR^x_k\cR^\bba_{j}f\|_\bbp<\infty\right\},
$$
where for a function $f:\mR^{2d}\to\mR$, 
$$
\cR^x_kf(x,v):=\cR_k f(\cdot,v)(x).
$$
Moreover, we also define
$$
\bB^{s}_{\bbp;x}:=\bB^{s,\infty}_{\bbp;x}:=\left\{f\in \cS': \|f\|_{\mathbf{B}^{s}_{\bbp;x}}
:= \sup_{k\geq 0}2^{k s}\|\cR^x_kf\|_\bbp<\infty\right\}.
$$
Similarly, one defines the usual  isotropic Besov spaces  $\bB^{s,q}_{p}$ 
in $\mR^d$ in terms of isotropic block operators $\cR_j$. If there is no confusion, we shall write
$$
\bB^s_{\bbp;\bba}:=\bB^{s,\infty}_{\bbp;\bba},\ \ \bB^s_{p}:=\bB^{s,\infty}_{p}.
$$
\end{definition}
\br\rm
The mixed Besov space $\bB^{s_0,s_1}_{\bbp;x,\bba}$ shall be used in the study of strong solutions of kinetic SDEs.
Note that for any $f\in \bB^{s,q}_{\bbp;\bba}$, by \eqref{AA13} we have
\begin{align}\label{SE1}
f=\sum_{j\geq 0}\cR^\bba_j f\  \mbox{ in $\bB^{s,q}_{\bbp;\bba}$.}
\end{align}
\er
\br\rm
Any finite Radon measure $\mu$ belongs to $\bB^{0,\infty}_{\bbb1;\bba}$. Indeed,
by definition and \eqref{SX4},
\begin{align}\label{Cx3}
\|\mu\|_{\bB^{0,\infty}_{\bbb1;\bba}}=\sup_{j\geq0}\|\cR^\bba_j\mu\|_{\bbb1}
\leq\sup_{j\geq0}\|\check\phi^\bba_j\|_{\bf 1} \|\mu\|_\cM
=\Big(\|\check\phi^\bba_0\|_1\vee\|\check\phi^\bba_1\|_1\Big) \|\mu\|_\cM<\infty,
\end{align}
where $\|\mu\|_\cM$ stands for the norm of total variation. 
Moreover, for any $m\in\mN$, if we define
$$
\mu_m(\dif z):=\mu(\dif(2^{-m\bba} z)),
$$
then for $j\geq 1$, by scaling $\check\phi^\bba_j(2^{m\bba} z)=2^{-m(2+\alpha)d}\check\phi^\bba_{m+j}(z)$,  it is easy to see that
$$
\cR^\bba_j\mu_m(z)=2^{-m(2+\alpha)d}(\cR^\bba_{j+m}\mu)(2^{-m\bba} z).
$$
Hence,  for any $\bbp\in[1,\infty]^2$,
$$
\|\cR^\bba_j\mu_m\|_{\bbp}=2^{m\bba\cdot(\frac d{\bbp}-\bbd)}\|\cR^\bba_{j+m}\mu\|_{\bbp},\ j\geq 1.
$$
Moreover, for $j=0$, noting that
$$
\sF(\check\phi^\bba_0(2^{m\bba}\cdot))(\xi)=2^{-m(2+\alpha)d} \phi^\bba_0(2^{-m\bba} \xi)
$$
and
$$
\phi^\bba_0(2^{-m\bba} \xi)=\phi^\bba_0(2^{-m\bba} \xi)\sum_{j=0}^m \phi^\bba_j(\xi),
$$
we have
$$
\check\phi^\bba_0(2^{m\bba}\cdot)=\sum_{j=0}^m\check\phi^\bba_0(2^{m\bba}\cdot)*\check\phi^\bba_j.
$$
Therefore, 
\begin{align*}
\|\cR^\bba_0\mu_m\|_{\bbp}&=2^{m\bba\cdot\frac d{\bbp}}\|\check\phi^\bba_0(2^{m\bba}\cdot)*\mu\|_\bbp
\leq 2^{m\bba\cdot\frac d{\bbp}}\sum_{j=0}^m\|\check\phi^\bba_0(2^{m\bba}\cdot)*\cR^\bba_j\mu\|_\bbp\\
&\le 2^{m\bba\cdot \frac d{\bbp}}\|\check\phi^\bba_0(2^{m\bba}\cdot)\|_{\bbb1}\sum_{j=0}^m\|\cR^\bba_{j}\mu\|_{\bbp}
\le 2^{m\bba\cdot(\frac d{\bbp}-\bbd)}\left[2\b1_{s>0}+\tfrac{2^{-(m+1)s}-1}{2^{-s}-1}\b1_{s<0}\right]\|\mu\|_{\bB^{s,\infty}_{\bbp;\bba}}.
\end{align*}
Combining the above calculations, we obtain that for any $s\in\mR\backslash\{0\}$ and $\bbp\in[1,\infty]^2$,
\begin{align}\label{SC8}
\|\mu_m\|_{\bB^{s,\infty}_{\bbp;\bba}}\le2^{m\bba\cdot(\frac d{\bbp}-\bbd)}\left[2\b1_{s>0}+\tfrac{2^{-ms}-2^s}{1-2^{s}}\b1_{s<0}\right]\|\mu\|_{\bB^{s,\infty}_{\bbp;\bba}},\ \ m\in\mN.
\end{align}
\er
For a function $f:\mR^{2d}\to\mR$,  the first-order difference operator is defined by
$$
\delta^{(1)}_hf(z):=\delta_hf(z):=f(z+h)-f(z),\ \ z, h\in\mR^{2d}.
$$
For $M\in\mN$, the $M$-order difference operator  is defined recursively by
$$
\delta^{(M+1)}_hf(z)=\delta_h\circ\delta^{(M)}_hf(z).
$$
The following characterization is well-known (cf. \cite{ZZ21} and \cite[Theorem 2.7]{HZZZ22}).
\bp
For $s>0$, $q\in[1,\infty]$ and $\bbp\in[1,\infty]^2$, an equivalent norm of $\bB^{s,q}_{\bbp;\bba}$ is given by
\begin{align}\label{CH1}
\|f\|_{\bB^{s,q}_{\bbp;\bba}}\asymp \left(\int_{|h|_\bba\leq 1}\left(\frac{\|\delta^{([s]+1)}f\|_{\bbp}}{|h|^s_\bba}\right)^q\frac{\dif h}{|h|^{(2+\alpha)d}_\bba}\right)^{1/q}+\|f\|_{\bbp},
\end{align}
where $[s]$ is the integer part of $s$. In particular, $\bC^s_{\bba}:=\mathbf{B}^{s,\infty}_{\infty;\bba}$ is the anisotropic 
H\"older-Zygmund space, and for $s\in(0,1)$, there is a constant $C=C(\alpha,d,s)>0$ such that
$$
\|f\|_{\bC^s_{\bba}}\asymp_C\|f\|_\infty+\sup_{z\not= z'}|f(z)-f(z')|/|z-z'|^{s}_\bba.
$$
Similarly,  the mixed and $x$-direction H\"older-Zygmund spaces are defined by
$$
\bC^{s_0,s_1}_{x,\bba}:=\bB^{s_0,s_1}_{\infty;  x,\bba},\ \ \bC^{s}_{x}:=\bB^{s}_{\infty;  x}.
$$
\ep

The following lemma provides important properties about Besov spaces. Although the proofs are standard, 
for the readers' convenience, we provide detailed proofs in the appendix.
\bl\label{lemB1}
\begin{enumerate}[(i)]
 \item For any $\bbp\in[1,\infty]^2$,  $s'>s$ and $q\in[1,\infty]$, it holds that
\begin{align}\label{AB2}
\bB^{0,1}_{\bbp;\bba}\hookrightarrow\mL^\bbp\hookrightarrow\bB^{0,\infty}_{\bbp;\bba},\ \ \bB^{s',\infty}_{\bbp;\bba}\hookrightarrow \bB^{s,1}_{\bbp;\bba}\hookrightarrow \bB^{s,q}_{\bbp;\bba}.
\end{align}
\item For any  $\beta,\beta_1,\beta_2\in\mR$, $q,q_1,q_2\in[1,\infty]$ and $\bbp,\bbp_1,\bbp_2\in[1,\infty]^2$ with
$$
\beta=\beta_1+\beta_2,\ \ 1+\tfrac{1}{\bbp}=\tfrac1{\bbp_1}+\tfrac1{\bbp_2},\ \ \tfrac{1}{q}=\tfrac1{q_1}+\tfrac1{q_2},
$$
it holds that
\begin{align}\label{Con}
\|f*g\|_{\bB^{\beta,q}_{\bbp;\bba}}\leq 5\|f\|_{\bB^{\beta_1,q_1}_{\bbp_1;\bba}}\|g\|_{\bB^{\beta_2,q_2}_{\bbp_2;\bba}}.
\end{align}
\item For any $s>0$, $q\in[1,\infty]$  and $\bbp, \bbp_1,\bbp_2\in [1,\infty]^2$ with $\tfrac1\bbp=\tfrac1{\bbp_1}+\tfrac1{\bbp_2}$,
there is a constant $C=C(\alpha, d,s,q,\bbp_1,\bbp_2)>0$ such that for any $(f,g)\in \bB^{s,q}_{\bbp_1;\bba}\times \bB^{s, q}_{\bbp_2;\bba}$,
\begin{align}\label{EB07}
\|fg\|_{\bB^{s,q}_{\bbp;\bba}}\lesssim_C
\|f\|_{\bB^{s,q}_{\bbp_1;\bba}}\|g\|_{\bbp_2}+\|f\|_{\bbp_1}\|g\|_{\bB^{s,q}_{\bbp_2;\bba}}.
\end{align}
\item Let $\bbp,\bbp_1,\bbp_2\in[1,\infty]^2$ and $s,s_0,s_1\in\mR$, $\theta\in[0,1]$, $q,q_1,q_2\in[1,\infty]$. 
Suppose that
\begin{align}\label{DD01}
\tfrac1{q}=\tfrac{1-\theta}{q_1}+\tfrac{\theta}{q_2},\ 
\tfrac1{\bbp}\leq\tfrac{1-\theta}{\bbp_1}+\tfrac{\theta}{\bbp_2},\ \ s-\bba\cdot\tfrac{d}{\bbp}=(1-\theta)\big(s_1-\bba\cdot\tfrac{d}{\bbp_1}\big)+\theta\big(s_2-\bba\cdot\tfrac{d}{\bbp_2}\big).
\end{align}
Then there is a constant $C=C(\bbp,\bbp_1,\bbp_2,s,s_1,s_2,\theta)>0$ such that
\begin{align}\label{Sob}
\|f\|_{\bB^{s,q}_{\bbp;\bba}}\lesssim_C\|f\|^{1-\theta}_{\bB^{s_1,q_1}_{\bbp_1;\bba}}\|f\|^\theta_{\bB^{s_2,q_2}_{\bbp_2;\bba}}.
\end{align}
In particular, for $\bbb1\leq\bbp_1\leq\bbp\leq\infty$, $q\in[1,\infty]$ and $s=s_1+\bba\cdot\big(\tfrac{d}{\bbp}-\tfrac{d}{\bbp_1}\big)$,
\begin{align}\label{Sob1}
\|f\|_{\bB^{s,q}_{\bbp;\bba}}\lesssim_C\|f\|_{\bB^{s_1,q}_{\bbp_1;\bba}}.
\end{align}
\item For any $\eps>0$ and $\bbp\in[1,\infty]^2$, $(s_0,s_1)\in\mR^2$, we have the following embeddings:
\begin{align}\label{IA2}
\bB^{{s_0-\eps,s_1+\eps}}_{\bbp; x,\bba}\hookrightarrow\bB^{{s_0,s_1}}_{\bbp; x,\bba}
\hookrightarrow\bB^{{s_0+\eps,s_1-\eps}}_{\bbp; x,\bba},
\end{align}
and for $s_1>0$,
\begin{align}\label{IA3}
\bB^{{s_0,s_1}}_{\bbp; x,\bba}\hookrightarrow
\bB^{{s_0+s_1}}_{\bbp; x}.
\end{align}
\item For $\theta\in[0,1]$, $\gamma,\beta\in\mR$ and $\bbp\in[1,\infty]^2$, there is a constant $C=C(d,\theta,\bbp)>0$
such that
\begin{align}\label{IA1}
\|f\|_{\bB^{\theta \gamma, (1-\theta)\beta}_{\bbp; x,\bba}}\lesssim\|f\|^\theta_{\bB^{ \gamma/(1+\alpha)}_{\bbp;x}}
\|f\|^{1-\theta}_{\bB^\beta_{\bbp;\bba}}.
\end{align}
\end{enumerate}
\el
\br\rm
Inequality \eqref{Sob} is a version of the classical Gagliardo-Nirenberg inequality in Besov spaces.
Inequality \eqref{Sob1} is the classical Sobolev inequality.
\er
\subsection{Time-weighted Schauder estimates}
In this subsection we prove the basic time-weighted Schauder estimates for the model kinetic equation:
$\p_t u=(\Delta_v-v\cdot\nabla_x)u+f$. For this goal, we need to make a detailed study for the kinetic semigroup.
Here much attention is paid on the large time  estimate of kinetic semigroup in Besov spaces.

For $\alpha\in(0,2]$,
let $p_t(z)=p_t(x,v)$ be the distributional density of $Z_t:=(\int_0^tL^{(\alpha)}_s\dif s,L^{(\alpha)}_t)$, where $L^{(\alpha)}_t$ is the standard $\alpha$-stable process as in the introduction. Since 
$(L^{(\alpha)}_{\lambda t})_{t\geq 0}\stackrel{(d)}{=}\lambda^{\frac1\alpha}(L^{(\alpha)}_t)_{t\geq 0}$ for $\lambda>0$, 
it is easy to see that $p_t(x,v)$ 
enjoys the following scaling property
\begin{align}\label{SCL1}
p_t(x,v)=t^{-\frac{(2+\alpha)d}{\alpha}}p_1(t^{-\frac{1+\alpha}\alpha}x,t^{-\frac1\alpha}v),\ \ t>0.
\end{align}
The kinetic semigroup of operator $\Delta^{\alpha/2}_v-v\cdot\nabla_x$ is given by
\begin{align}\label{SCL11}
P_t f(z):=\mE f(\Gamma_t z+Z_t)=\Gamma_t(p_t*f)(z)=(\Gamma_tp_t*\Gamma_tf)(z)=\int_{\mR^d} \Gamma_t{p_t}(z-z') \Gamma_t{f}(z')\dif z',
\end{align}
where 
$$
\Gamma_tf(z):=f(\Gamma_tz):=f(x-tv,v).
$$ 

The following lemma is originally proven in \cite{HWZ20} (see also \cite{ZZ21}).
\bl\label{lem001}
For any $l\geq 0$, there is a constant $C=C(l,\alpha,d)>0$ such that
$$
\|\cR^\bba_j p_t \|_{1}+\|\cR^\bba_j \Gamma_t p_t \|_{1}\lesssim_C (t 2^{j\alpha})^{-l}\wedge 1,\ \forall j\geq 1, t>0.
$$
\el
\br\rm
The above estimate only holds for high frequency operator $\cR^\bba_j$, $j\geq 1$.
\er
As an easy consequence we have
\bl\label{Le10}
For any $\bbp\in[1,\infty]^2$, $k\in\mN_0$ and $\beta\in\mR$ , there is a constant $C>0$ such that
$$
\|\nabla^k_v p_t\|_{\bB^{\beta}_{\bbp;\bba}}\lesssim_C t^{\frac{\bba\cdot(d/\bbp-\bbd)-k}{\alpha}}(1+t^{-\frac{\beta\vee 0}\alpha}),\ t>0.
$$
\el
\begin{proof}
By scaling \eqref{SCL1} and the change of variables, it is easy to see that for any $m\in\mN_0$,
\begin{align}\label{SCL12}
\|\nabla^m_vp_t\|_\bbp=t^{\frac{\bba\cdot(d/\bbp-\bbd)-m}{\alpha}}\|\nabla^m_vp_1\|_\bbp\lesssim_C t^{\frac{\bba\cdot(d/\bbp-\bbd)-m}{\alpha}}.
\end{align}
Without loss of generality, we assume $\beta\geq 0$. By definition,  Bernstein's inequality \eqref{Ber} and 
Lemma \ref{lem001} with $l=\frac{\beta+k+\bba\cdot(\bbd-d/\bbp)}{\alpha}$, we have
\begin{align*}
\|\nabla^k_vp_t\|_{\bB^{\beta}_{\bbp;\bba}}
&\leq\|\cR^\bba_0\nabla^k_vp_t\|_\bbp+ \sup_{j\geq1}\big(2^{ j\beta}\|\cR^\bba_{j}\nabla^k_v p_t\|_{\bbp}\big)\\
&\lesssim\|\nabla^k_vp_t\|_\bbp+ \sup_{j\geq1}2^{j(\beta+k+\bba\cdot(\bbd-d/\bbp))}
\|\cR^\bba_{j}p_t\|_1\\
&\lesssim t^{\frac{\bba\cdot(d/\bbp-\bbd)-k}{\alpha}}+t^{\frac{\bba\cdot(d/\bbp-\bbd)-k-\beta}{\alpha}}.
\end{align*}
The proof is complete.
\end{proof}

The following lemma is proved in \cite[Lemma 6.7]{HWZ20}. For the readers' convenience, we reproduce the proof here.
\bl\label{Le29}
For $t>0$ and $j\geq 1$, we define
$$
\Theta^t_j:=\Big\{\ell\in\mN_0\,\big{|}\,2^\ell\leq2^{\alpha+2}(2^j+t2^{j(1+\alpha)}),\ 2^j\leq2^{\alpha+2}(2^\ell+t2^{\ell(1+\alpha)})\Big\}.
$$
\begin{enumerate}[(i)]
\item For any $\ell\notin\Theta^t_j$, it holds that
$$
\<\cR^\bba_j f,\Gamma_t\cR^\bba_\ell g\>=0.
$$
\item For any $0\not=\beta\in\mR$, it holds that
\begin{align}\label{AS88}
\sum_{\ell\in\Theta^t_j}2^{-\ell\beta}\lesssim 2^{-j\beta}(1+(t2^{j\alpha}))^{|\beta|}.
\end{align}
\end{enumerate}
\el
\begin{proof}
(i) By Parsavel's identity of Fourier's transform, we have
\begin{align*}
&\<\cR^\bba_j f,\Gamma_t\cR^\bba_\ell g\>
=\int_{\mR^{2d}}\widehat{\cR^\bba_j f}(\xi_1,\xi_2)\widehat{\Gamma_t\cR^\bba_\ell g}(\xi_1,\xi_2)\dif\xi_1\dif\xi_2\\
&\qquad=\int_{\mR^{2d}}\phi^\bba_j(\xi_1,\xi_2)f(\xi_1,\xi_2)\phi^\bba_\ell(\xi_1,\xi_2-t\xi_1)g(\xi_1,\xi_2-t\xi_1)\dif\xi_1\dif\xi_2.
\end{align*}
Note that
$$
{\rm supp}\phi^\bba_j\subset\{\xi: 2^{j-1}\leq|\xi|_\bba\leq 2^{j+1}\}.
$$
For fixed $j\in\mN$ and $t>0$, it is easy to see that $\<\cR^\bba_j f,\Gamma_t\cR^\bba_\ell g\>=0$ if and only if
$$
\{\xi: 2^{j-1}\leq|\xi|_\bba\leq 2^{j+1}\}\cap\{\xi=(\xi_1,\xi_2): 2^{\ell-1}\leq|(\xi_1, \xi_2-t\xi_1)|_\bba\leq 2^{\ell+1}\}=\emptyset,
$$
which is equivalent to say $\ell\notin\Theta^t_j$.

(ii) Fix $\beta>0$. Note that for $\ell\in\Theta^t_j$ with $\ell\leq j$, 
$$
2^{-\ell}\leq 2^{\alpha+2} 2^{-j}(1+t2^{\ell\alpha})\leq 2^{\alpha+2}2^{-j} (1+t2^{j\alpha})=:D\Rightarrow \ell\geq-\ln D/\ln 2.
$$
Hence,
\begin{align*}
\sum_{\ell\in\Theta^t_j}2^{-\ell\beta}&=\sum_{\ell\in\Theta^t_j;\ell\geq j}2^{-\ell\beta}+\sum_{\ell\in\Theta^t_j;\ell<j}2^{-\ell\beta}
\leq\sum_{\ell\geq j}2^{-\ell\beta}+\sum_{-\ln D/\ln 2\leq \ell}2^{-\ell\beta}\\
&\leq \frac{2^{-j\beta}}{1-2^{-\beta}}+ \frac{D^\beta}{1-2^{-\beta}}\lesssim 2^{-j\beta}(1+t2^{j\alpha})^\beta.
\end{align*}
Similarly, for $D:=2^{\alpha+2}2^j(1+t2^{j\alpha})$, we have
\begin{align*}
\sum_{\ell\in\Theta^t_j}2^{\ell\beta}&=\sum_{\ell\in\Theta^t_j;\ell\geq j}2^{\ell\beta}+\sum_{\ell\in\Theta^t_j;\ell<j}2^{\ell\beta}
\leq\sum_{\ell\leq\ln D/\ln 2}2^{\ell\beta}+\sum_{\ell<j}2^{\ell\beta}\\
&\lesssim D^\beta+2^{j\beta}\lesssim 2^{j\beta}(1+t2^{j\alpha})^\beta.
\end{align*}
The proof is complete.
\end{proof}
Using the above two lemmas, we can show the following crucial lemma.
\bl
Let $\bbp,\bbp'\in[1,\infty]^2$ with $\bbp\leq\bbp'$. Write $\sA:=\bba\cdot(\tfrac d{\bbp}-\tfrac d{\bbp'})$.  
Recall $\bB^{\beta}_{\bbp;\bba}:=\bB^{\beta,\infty}_{\bbp;\bba}$.
\begin{enumerate}[(i)]
\item {\bf (High frequency part of $P_tf$)} For any $\beta\in\mR$ and $\gamma\geq 0$, there is a constant $C=C(\alpha, d,\bbp,\bbp',\beta,\gamma)>0$ such that for any $f\in \bB^{\beta}_{\bbp;\bba}$ and all  
$j\geq 1$ and $t>0$,
\begin{align}
\|\cR^\bba_jP_tf\|_{\bbp'}
\lesssim_C 2^{j(\sA-\beta)}
((2^{j\alpha}t)^{-\gamma}\wedge 1)\|f\|_{\bB^{\beta}_{\bbp;\bba}}.\label{AD0306}
\end{align}
\item {\bf (Low frequency part of $P_tf$)} For any  $\beta\in\mR$, there is a constant $C=C(\alpha, d,\bbp,\bbp',\beta)>0$ such that for  any $f\in \bB^{\beta}_{\bbp;\bba}$ and all $t>0$,
\begin{align}
\|\cR^\bba_0 P_tf\|_{\bbp'}\lesssim_C (1\vee t)^{-\frac{\sA}{\alpha}} \|f\|_{\bB^{\beta}_{\bbp;\bba}},\label{AD0396}
\end{align}
and for any $F\in ( \bB^{\beta}_{\bbp;\bba})^d$ and all $t>0$,
\begin{align}
\|\cR^\bba_0 P_t\div_vF\|_{\bbp'}\lesssim_C (1\vee t)^{-\frac{\sA+1}{\alpha}} \|F\|_{\bB^{\beta}_{\bbp;\bba}}.\label{AD0316}
\end{align}
\end{enumerate}
\el
\begin{proof}
(i) By the interpolation theorem, we may assume $\beta\not=0$.
 Note that by (i) of Lemma \ref{Le29},
\begin{align}
\cR^\bba_jP_tf=\cR^\bba_j\Big(\Gamma_tp_t*\Gamma_tf\Big)
\stackrel{\eqref{SE1}}{=}\sum_{\ell\geq 0}\cR^\bba_j\Gamma_tp_t*\Gamma_t\cR^\bba_\ell f
=\sum_{\ell\in\Theta^t_j}\cR^\bba_j\Gamma_tp_t*\Gamma_t\cR^\bba_\ell f.\label{Mm1}
\end{align}
By Bernstein's inequality \eqref{Ber} and Young's inequality,  we have
\begin{align*}
\|\cR^\bba_jP_t f\|_{\bbp'}&\lesssim 2^{j\sA}\|\cR^\bba_jP_t f\|_{\bbp}\leq 2^{j\sA}\|\cR^\bba_j\Gamma_tp_t\|_1\sum_{\ell\in\Theta^t_j}\|\Gamma_t\cR^\bba_\ell f\|_{\bbp}\stackrel{\eqref{AA1}}{=}2^{j\sA}\|\cR^\bba_j\Gamma_tp_t\|_1\sum_{\ell\in\Theta^t_j}\|\cR^\bba_\ell f\|_{\bbp}\no\\
&\leq 2^{j\sA}\|\cR^\bba_j\Gamma_tp_t\|_1\sum_{\ell\in\Theta^t_j}2^{-\ell\beta}\|f\|_{\bB^{\beta}_{\bbp;\bba}}
\stackrel{\eqref{AS88}}{\lesssim} \|\cR^\bba_j\Gamma_tp_t\|_12^{j(\sA-\beta)}(1+(t2^{j\alpha}))^{|\beta|}\|f\|_{\bB^{\beta}_{\bbp;\bba}},
\end{align*}
which together with Lemma \ref{lem001} with $l=\gamma$ and $l=\gamma+|\beta|$ yields \eqref{AD0306}.

(ii) We only prove \eqref{AD0316}.
Noting that by definitions,
$$
\cR^\bba_0P_t\div_v F=\check\phi^\bba_0*\Gamma_t(p_t*\div_v F)
=\Gamma_t(\Gamma_{-t}\check\phi^\bba_0*(p_t*\div_vF)),
$$
we have
\begin{align}\label{AB1}
\|\cR^\bba_0P_t\div_v F\|_{\bbp'}\stackrel{\eqref{AA1}}{=}\|\Gamma_{-t}\check\phi^\bba_0*(p_t*\div_vF)\|_{\bbp'}.
\end{align}
By Young's inequality and integration by parts, we have
for $\eps>0$, 
$$
\|\cR^\bba_0P_t\div_v F\|_{\bbp'}\leq\|\Gamma_{-t}\check\phi^\bba_0\|_1\|\nabla_v p_t*F\|_{\bbp'}
\lesssim\|\check\phi^\bba_0\|_1\|\nabla_v p_t*F\|_{\bB^\eps_{\bbp';\bba}},
$$
which, by  \eqref{Con} and Lemma \ref{Le10},  implies  that for $\bbq\in[1,\infty]$ being defined by $\bbb1+\tfrac1{\bbp'}=\tfrac 1{\bbq}+\tfrac1{\bbp}$ and
$t\geq 1$, 
\begin{align}\label{AA3}
\|\cR^\bba_0P_t\div_v F\|_{\bbp'}\lesssim\|\nabla_v p_t\|_{\bB^{\eps-\beta}_{\bbq;\bba}}\|F\|_{\bB^{\beta}_{\bbp;\bba}}\lesssim 
t^{\frac{\bba\cdot(d/\bbq-\bbd)-1}{\alpha}} \|F\|_{\bB^{\beta}_{\bbp;\bba}}=t^{-\frac{\sA+1}{\alpha}} \|F\|_{\bB^{\beta}_{\bbp;\bba}}.
\end{align}
Moreover, for $t\in(0,1]$ and $\eps>0$, by \eqref{AB1} and \eqref{Con},
\begin{align*}
\|\cR^\bba_0P_t\div_v F\|_{\bbp'}
&=\|p_t*(\Gamma_{-t}\check\phi^\bba_0*\div_vF)\|_{\bbp'}\leq \|p_t\|_1\|\Gamma_{-t}\check\phi^\bba_0*\div_vF\|_{\bbp'}\\
&=\|(\nabla_v\Gamma_{-t}\check\phi^\bba_0)*F\|_{\bbp'}\lesssim\|(\nabla_v\Gamma_{-t}\check\phi^\bba_0)*F\|_{\bB^\eps_{\bbp';\bba}}\\
&\lesssim\|\nabla_v\Gamma_{-t}\check\phi^\bba_0\|_{\bB^{\eps-\beta}_{\bbq;\bba}}\|F\|_{\bB^{\beta}_{\bbp;\bba}}.
\end{align*}
Since $\check\phi^\bba_0\in\cS$, by the chain rule, we have
$$
\sup_{t\in[0,1]}\|\nabla_v\Gamma_{-t}\check\phi^\bba_0\|_{\bB^{\eps-\beta}_{\bbq;\bba}}
\lesssim\sup_{t\in[0,1]}\sum_{j\leq |\eps-\beta|+1}\|\nabla^j_v\Gamma_{-t}\check\phi^\bba_0\|_\bbq<\infty.
$$
Hence,
$$
\sup_{t\in[0,1]}\|\cR^\bba_0P_t\div_v F\|_{\bbp'}\lesssim\|F\|_{\bB^{\beta}_{\bbp;\bba}},
$$
which together with \eqref{AA3} yields \eqref{AD0316}.
\end{proof}
\br\rm
The differentiability index $\beta$ in \eqref{AD0396} and \eqref{AD0316}  does not make any contribution about the decay rate of large time of $\|\cR^\bba_0 P_t\div_vF\|_{\bbp'}$.
This is natural since we are using the non-homoegenous norm in the right hand side.
\er

Fix $\alpha\in(0,2]$. For given $\bbg=(\gamma_0,\gamma_1)\in\mR^2$, $q\in[1,\infty]$ and a Banach space $\mB$, 
let $\mL^q_\bbg(\mB)$ be the space of all functions $f:(0,\infty)\to\mB$ with finite norm
\begin{align*}
\|f\|_{\mL^{q}_{\bbg}(\mB)}:=\left(\int^\infty_0\(\|f(t)\|_{\mB}(1\wedge t)^{\frac{\gamma_0}\alpha}(1\vee t)^{\frac{\gamma_1}\alpha}\)^q\dif t\right)^{1/q}<\infty.
\end{align*}
If $\gamma_0,\gamma_1$ are positive,
then $\gamma_0$ and $\gamma_1$ stand for the small time singularity and large time decay weights, respectively.

We shall use the following relation for two weights $\bbg,\bbg'\in\mR^2$.
\bd\label{P1}
For  $\bbg=(\gamma_0,\gamma_1)$ and $\bbg'=(\gamma'_0,\gamma'_1)\in\mR^2$, 
one says $\bbg'\wleqq\bbg$ if $\gamma'_0\geq\gamma_0$ and $\gamma'_1\leq\gamma_1$,
 and $\bbg'\wle\bbg$ if $\gamma'_0\geq\gamma_0$ and $\gamma'_1<\gamma_1$.
In particular, if $\bbg'\wleqq\bbg$, then
\begin{align}\label{P2}
\|f\|_{\mL^{q}_{\bbg'}(\mB)}\leq\|f\|_{\mL^{q}_{\bbg}(\mB)}.
\end{align}
\ed
The following proposition is simple, but quite useful for proving the large time decay estimates and 
the smoothness of solutions by shifting the time variable (see Lemma \ref{Le38} below).
\bp
For any $\bbg,\bbg'\in[0,\infty)^2$, it holds that
\begin{align}\label{Sh1}
\|f(t+\cdot)\|_{\mL^q_{\gamma'}(\mB)}\leq (1\wedge t)^{-\frac{\gamma_0}\alpha}(1\vee t)^{-\frac{\gamma_1}\alpha}\|f\|_{\mL^q_{\bbg+\bbg'}(\mB)},\ \ \forall t>0.
\end{align}
\ep
\begin{proof}
Noting that for $\bbg=(\gamma_0,\gamma_1),\bbg'=(\gamma'_0,\gamma'_1)\in[0,\infty)^2$ and $s,t>0$,
$$
(1\wedge s)^{\frac{\gamma'_0}\alpha} \leq(1\wedge (t+s))^{\frac{\gamma_0+\gamma'_0}\alpha} (1\wedge t)^{-\frac{\gamma_0}\alpha}
$$
and
$$
(1\vee s)^{\frac{\gamma'_0}\alpha} \leq(1\vee (t+s))^{\frac{\gamma_0+\gamma'_0}\alpha} (1\vee t)^{-\frac{\gamma_0}\alpha},
$$
by definition, we have
\begin{align*}
\|f(t+\cdot)\|_{\mL^q_{\bbg'}(\mB)}&=\left(\int_0^\infty\Big(\|f(t+s)\|_{\mB}(1\wedge s)^{\frac{\gamma'_0}\alpha}(1\vee s)^{\frac{\gamma'_1}\alpha}\Big)^q\dif s\right)^{1/q}\\
&\leq\left(\int^\infty_t\Big(\|f(s)\|_{\mB}(1\wedge s)^{\frac{\gamma_0+\gamma_0'}\alpha}(1\vee s)^{\frac{\gamma_1+\gamma_1'}\alpha}\Big)^q\dif s \right)^{1/q}(1\wedge t)^{-\frac{\gamma_0}\alpha}(1\vee t)^{-\frac{\gamma_1}\alpha}\\
&\leq(1\wedge t)^{-\frac{\gamma_0}\alpha}(1\vee t)^{-\frac{\gamma_1}\alpha}\|f\|_{\mL^q_{\bbg+\bbg'}(\mB)}.
\end{align*}
The proof is complete.
\end{proof}
The following lemma provides useful estimates about the kinetic semigroup $P_t$.
 \bl\label{Le215}
Let $\bbp_0,\bbp_1,\bbp'\in[1,\infty]^2$ with $\bbp_0,\bbp_1\leq\bbp'$.  Write $\sA_i:=\bba\cdot(\tfrac d{\bbp_i}-\tfrac d{\bbp'})$, $i=0,1$. For any $\beta,\beta'\in\mR$, there are constants $C_i=C_i(\alpha,d,\beta,\beta',\bbp_i,\bbp')>0$, $i=0,1$ such that
for all $t>0$,
\begin{align}\label{AD0446}
\b1_{\beta'\not=\beta-\sA_0}\|P_tf\|_{\bB^{\beta',1}_{\bbp';\bba}}+\|P_tf\|_{\bB^{\beta'}_{\bbp';\bba}}
\lesssim_{C_0}(1\wedge t)^{-\frac{(\beta'-\beta+\sA_0)\vee0}\alpha}\|f\|_{\bB^{\beta}_{\bbp_0;\bba}},
\end{align}
where $\bB^{\beta}_{\bbp;\bba}=\bB^{\beta,\infty}_{\bbp;\bba}$, and for $t\geq 1$,
\begin{align}\label{AD0456}
\|P_tf\|_{\bB^{\beta',1}_{\bbp';\bba}}\lesssim_{C_1}t^{-\frac{\sA_1}\alpha} \|f\|_{\bB^{\beta}_{\bbp_1;\bba}}.
\end{align}
\el
\begin{proof}
(i) Let $\beta_0:=\beta'-\beta+\sA_0\not=0$.
By \eqref{AD0396} and \eqref{AD0306} with $\gamma>\frac{\beta_0\vee 0}\alpha$, we have for all $t>0$,
\begin{align*}
\|P_tf\|_{\bB^{\beta',1}_{\bbp';\bba}}&=\|\cR^\bba_0 P_tf\|_{\bbp'}+\sum_{j\geq 1}2^{j\beta'}\|\cR^\bba_j P_tf\|_{\bbp'}\\
&\lesssim\|f\|_{\bB^{\beta}_{\bbp_0;\bba}}+\sum_{j\geq 1}2^{j\beta_0}((2^{j\alpha}t)^{-\gamma}\wedge 1)\|f\|_{\bB^{\beta}_{\bbp_0;\bba}}.
\end{align*}
Since $\beta_0\not=0$ and $\gamma\alpha>\beta_0\vee0$, by Lemma \ref{LA1} we have
$$
\sum_{j\geq 1}2^{j\beta_0}((2^{j\alpha}t)^{-\gamma}\wedge 1)
\lesssim t^{-\frac{\beta_0\vee 0}\alpha},\ t>0.
$$
Thus we get
$$
\|P_tf\|_{\bB^{\beta',1}_{\bbp';\bba}}\lesssim \|f\|_{\bB^{\beta}_{\bbp_0;\bba}}+t^{-\frac{\beta_0\vee 0}\alpha}\|f\|_{\bB^{\beta}_{\bbp_0;\bba}}\leq 2 (1\wedge t)^{-\frac{(\beta'-\beta+\sA_0)\vee0}\alpha}\|f\|_{\bB^{\beta}_{\bbp_0;\bba}}.
$$
If $\beta_0=0$, we directly have
$$
\|P_tf\|_{\bB^{\beta'}_{\bbp';\bba}}\leq\sup_{j\geq 0}2^{j\beta'}\|\cR^\bba_j P_tf\|_{\bbp'}\lesssim\|f\|_{\bB^{\beta}_{\bbp_0;\bba}}.
$$
Combining the above two estimates, we obtain \eqref{AD0446}.

(ii) Since $\sA_1=\bba\cdot(\frac\bbd{\bbp_1}-\tfrac d{\bbp'})$,
by \eqref{AD0396} and \eqref{AD0306}  with $\gamma>\frac{\sA_1}\alpha$ large enough, we have for $t\geq 1$,
\begin{align*}
\|P_tf\|_{\bB^{\beta',1}_{\bbp';\bba}}
&=\|\cR^\bba_0 P_tf\|_{\bbp'}+\sum_{j\geq 1}2^{j\beta'}\|\cR^\bba_j P_tf\|_{\bbp'}\\
&\lesssim t^{-\frac{\sA_1}\alpha} \|f\|_{\bB^{\beta}_{\bbp_1;\bba}}
+t^{-\gamma}\sum_{j\geq 1}2^{j(\beta'+\sA_1-\alpha\gamma)}\|f\|_{\bB^{\beta}_{\bbp_1;\bba}}
\lesssim t^{-\frac{\sA_1}\alpha} \|f\|_{\bB^{\beta}_{\bbp_1;\bba}}.
\end{align*}
The proof is complete.
\end{proof}
\br\rm
Estimates \eqref{AD0446} and \eqref{AD0456} provide short time and large time estimates of the kinetic semigroup, respectively.
In particular, the  rate of large time decay only depends on the integrability indices. Note that 
$$
\|P_tf\|_{\bbp}\stackrel{\eqref{AB2}}{\lesssim}\|P_tf\|_{\bB^{0,1}_{\bbp;\bba}}
\stackrel{\eqref{AD0456}}{\lesssim}\|f\|_{\bB^{0,\infty}_{\bbp;\bba}}\stackrel{\eqref{AB2}}{\lesssim}_C\|f\|_\bbp,\ \ t\geq 1. 
$$
It is the same reason that \eqref{AD0446} is stronger than the classical estimate:
$$
\|\nabla_v P_tf\|_{\bbp}\lesssim(1\wedge t)^{-\frac{1}\alpha}\|f\|_\bbp.
$$
Moreover, for any $\beta,\beta'\in\mR$ and $\bbp\in[1,\infty]^2$, 
noting that by \eqref{AD0446},
$$
\|P_tf\|_{\bB^{\beta+\beta'}_{\bbp;\bba}}\lesssim(1\wedge t)^{-\frac{\beta'\vee0}{\alpha}}\|f\|_{\bB^{\beta}_{\bbp;\bba}},\ \ t>0,
$$
by the real interpolation theorem (see \cite[Theorem 6.4.5, (1)]{BL76}), we have for any $q\in[1,\infty]$,
\begin{align}\label{AD335}
\|P_tf\|_{\bB^{\beta+\beta',q}_{\bbp;\bba}}\lesssim_C(1\wedge t)^{-\frac{\beta'\vee0}{\alpha}}\|f\|_{\bB^{\beta,q}_{\bbp;\bba}},\ \ t>0.
\end{align}
\er
For $f\in\cS'$ and $t>0$, let
$$
u(t,x,v):=\sI_tf(x,v):=\int^t_0 P_{t-s}f(s,x,v)\dif s.
$$
By the definition \eqref{SCL11} of $P_t$, it is easy to see that in the distributional sense,
$$
\p_t u=\Delta^{\alpha/2}_vu-v\cdot\nabla_x u+f.
$$
We now show the following basic Schauder estimate in time-weighted spaces.

\bt\label{Th26}
(Time-weighted Schauder's estimate) 
Let $\alpha\in(0,2)$, $q\in(1,\infty]$ and $\bbp,\bbp'\in[1,\infty]^2$ with $\bbp\leq\bbp'$.
Write
$$
\sQ:=\alpha-\tfrac\alpha q,\ \ \sA:=\bba\cdot(\tfrac d{\bbp}-\tfrac d{\bbp'}).
$$ 
Recall $\bB^{\beta}_{\bbp;\bba}:=\bB^{\beta,\infty}_{\bbp;\bba}$. Suppose that $\theta\in(0,\sQ)$ and
$\bbg=(\gamma_0,\gamma_1),\bbg'=(\gamma'_0,\gamma'_1)\in\mR^2$ satisfy
\begin{align}\label{CC1}
\bbg\in [0,\sQ)^2,\ \ \bbg'\wle\bbg-\big(\theta, (\sQ-1-\sA)\vee 0\big),
\end{align}
see Definition \ref{P1} for notation $\wle$.
Let $\Theta:=(\alpha,d,q,\bbp,\bbp',\bbg,\bbg',\theta)$ be the parameter set and
$$
\beta':=\sQ-1-\sA-\theta.
$$
 For any  $\beta\in\mR$,  there is a
constant $C=C(\beta,\Theta)>0$ such that for any $F\in (\mL^q_{\bbg}(\bB^{\beta}_{\bbp;\bba}))^d$,
\begin{align}\label{Scha}
\|\sI\div_v F\|_{\mL^\infty_{\bbg'}(\bB^{\beta+\beta',1}_{\bbp';\bba})}
\lesssim_C\|F\|_{\mL^q_{\bbg}(\bB^{\beta}_{\bbp;\bba})}.
\end{align}
Moreover, if $\theta=0$ or $\sQ$, then
\begin{align}\label{Scha0}
\|\sI\div_v F\|_{\mL^\infty_{\bbg'}(\bB^{\beta+\beta'}_{\bbp';\bba})}
\lesssim_C\|F\|_{\mL^q_{\bbg}(\bB^{\beta}_{\bbp;\bba})}.
\end{align}
\et
\begin{proof}
First of all, by \eqref{AD0316}, H\"older's inequality and Lemma \ref{lemA2}, we have
\begin{align}
\|\cR^\bba_0 \sI_t\div_v F\|_{\bbp'}
&\lesssim \int^t_0\|\cR^\bba_0 P_{t-s}\div_v F_s\|_{\bbp'}\dif s
\lesssim\int^t_0(1\vee (t-s))^{-\frac{\sA+1}\alpha} \|F_s\|_{\bB^{\beta}_{\bbp;\bba}}\dif s\no\\
&\lesssim \left(\int^t_0 \left((1\wedge (t-s)^{-\frac{\sA+1}\alpha})  (1\wedge s)^{-\frac{\gamma_0}{\alpha}}(1\vee s)^{-\frac{\gamma_1}{\alpha}}\right)^{r}\dif s\right)^{\frac1{r}}
\|F\|_{\mL^q_\bbg(\bB^{\beta}_{\bbp;\bba})}\no\\
&\lesssim (1\wedge t)^{\frac1{r}-\frac{\gamma_0}{\alpha}}(1\vee t)^{\frac{1}{r}-{\frac{\gamma_1}{\alpha}}}\ell^{1/r}_\vartheta(t)
\|F\|_{\mL^q_\bbg(\bB^{\beta}_{\bbp;\bba})}\label{MX1},
\end{align}
where $\vartheta:=r(\sA+1)/\alpha$ and
$$
\ell_\vartheta(t):=(1\vee t)^{-\vartheta\wedge 1}(1+\b1_{\vartheta=1}\ln(1\vee t)).
$$
%we used $\sA+1\not=\frac\alpha r=\sQ$ in the last $\lesssim$, and in the equality we used that
Noting that by \eqref{CC1} and  $\frac \alpha r=\alpha-\frac\alpha q=\sQ$,
$$
\alpha(\tfrac1{r}-\tfrac{\gamma_1}{\alpha}-\tfrac{\sA+1}{\alpha}\wedge\tfrac1{r})=\sQ-\gamma_1-(\sA+1)\wedge\sQ=(\sQ-1-\sA)\vee0-\gamma_1<-\gamma'_1,
$$
we have
$$
(1\vee t)^{\frac{1}{r}-{\frac{\gamma_1}{\alpha}}}\ell^{1/r}_\vartheta(t)
=(1\vee t)^{\frac{(\sQ-1-\sA)\vee0-\gamma_1}\alpha}(1+\b1_{\vartheta=1}\ln(1\vee t))^{\frac1r}
\lesssim(1\vee t)^{-\frac{\gamma'_1}\alpha}.
$$
Substituting this into \eqref{MX1}, we obtain
\begin{align}\label{Mz1}
\|\cR^\bba_0 \sI_t\div_v F\|_{\bbp'}
\lesssim(1\wedge t)^{\frac{\sQ-\gamma_0}{\alpha}}(1\vee t)^{-\frac{\gamma'_1}{\alpha}}
\|F\|_{\mL^q_\bbg(\bB^{\beta}_{\bbp;\bba})}.
\end{align}
Next, let $\frac1{q}+\frac1{r}=1$. By \eqref{AD0306}, H\"older's inequality and Lemma \ref{lemA2}, we have
for any $j\ge1$,
\begin{align}
\|\cR^\bba_j\sI_t\div_vF\|_{\bbp'}
&\leq\int^t_0\|\cR^\bba_jP_{t-s}\div_vF_s\|_{\bbp'}\dif s
 \lesssim 2^{j(\sA+1-\beta)}\int^t_0 \left[((t-s) 2^{j\alpha})^{-2}\wedge1\right]
\|\div_vF_s\|_{\bB^{\beta-1}_{\bbp;\bba}}\dif s\no\\
&\leq 2^{j(\sA+1-\beta)}
\left[\int^t_0 \Big([((t-s) 2^{j\alpha})^{-2}\wedge1] (1\wedge s)^{-\frac{\gamma_0}{\alpha}}(1\vee s)^{-\frac{\gamma_1}{\alpha}}\Big)^{r}\dif s\right]^{\frac1{r}}
\|F\|_{\mL^q_\bbg(\bB^{\beta}_{\bbp;\bba})} \no\\
&\lesssim 2^{j(\sA+1-\beta)} ((2^{j\alpha} t)^{-\frac1{r}}\wedge 1)(1\wedge t)^{\frac1{r}-\frac{\gamma_0}{\alpha}}(1\vee t)^{\frac1{r}-{\frac{\gamma_1}{\alpha}}}\|F\|_{\mL^q_\bbg(\bB^{\beta}_{\bbp;\bba})}.\label{Bz1}
\end{align}
\iffalse
Since $\frac\alpha{r}=\sQ$, for any $\theta\in(0,\sQ)$, by Lemma \ref{LA1} we get
\begin{align*}
&\sum_{j\geq 1}2^{j(\beta+\sQ-1-\sA-\theta)}\|\cR^\bba_j\sI_t\div_vF\|_{\bbp'}\\
&\qquad \lesssim \int^t_0 \sum_{j\geq 1} 2^{j(\sQ-\theta)}\left[((t-s) 2^{j\alpha})^{-2}\wedge1\right]
\|\div_vF_s\|_{\bB^{\beta-1}_{\bbp;\bba}}\dif s\no\\
&\qquad \lesssim \int^t_0 (t-s)^{-\frac{\sQ-\theta}{\alpha}}
\|\div_vF_s\|_{\bB^{\beta-1}_{\bbp;\bba}}\dif s\no\\
&\qquad \lesssim 
\left[\int^t_0\Big((t-s)^{-\frac{\sQ-\theta}{\alpha}} (1\wedge s)^{-\frac{\gamma_0}{\alpha}}(1\vee s)^{-\frac{\gamma_1}{\alpha}}\Big)^{r}\dif s\right]^{\frac1{r}}
\|F\|_{\mL^q_\bbg(\bB^{\beta}_{\bbp;\bba})}\\
&\qquad\lesssim
(1\vee t)^{-\frac{(\sQ-\theta)\wedge\sQ}{\alpha}}(1\wedge t)^{\frac1{r}-\frac{\gamma_0}{\alpha}}(1\vee t)^{\frac1{r}-{\frac{\gamma_1}{\alpha}}}\|F\|_{\mL^q_\bbg(\bB^{\beta}_{\bbp;\bba})}\\
&\qquad\lesssim
(1\wedge t)^{\frac{\sQ-\gamma_0}{\alpha}}(1\vee t)^{{\frac{\sQ-\gamma_1}{\alpha}}}\|F\|_{\mL^q_\bbg(\bB^{\beta}_{\bbp;\bba})}\\
\end{align*}
\fi
Since $\frac\alpha{r}=\sQ$, for any $\theta\in(0,\sQ)$, by Lemma \ref{LA1} we get
\begin{align*}
&\sum_{j\geq 1}2^{j(\beta+\sQ-1-\sA-\theta)}\|\cR^\bba_j\sI_t\div_vF\|_{\bbp'}\\
&\qquad\lesssim
\sum_{j\geq 1}2^{j(\sQ-\theta)}((2^{j\alpha} t)^{-\frac1{r}}\wedge 1)(1\wedge t)^{\frac1{r}-\frac{\gamma_0}{\alpha}}(1\vee t)^{\frac1{r}-{\frac{\gamma_1}{\alpha}}}\|F\|_{\mL^q_\bbg(\bB^{\beta}_{\bbp;\bba})}\\
&\qquad\lesssim
t^{-\frac{\sQ-\theta}{\alpha}}(1\wedge t)^{\frac{\sQ-\gamma_0}{\alpha}}(1\vee t)^{{\frac{\sQ-\gamma_1}{\alpha}}}\|F\|_{\mL^q_\bbg(\bB^{\beta}_{\bbp;\bba})}\\
&\qquad=
(1\wedge t)^{\frac{\theta-\gamma_0}{\alpha}}(1\vee t)^{{\frac{\theta-\gamma_1}{\alpha}}}\|F\|_{\mL^q_\bbg(\bB^{\beta}_{\bbp;\bba})}.
\end{align*}
Since $\gamma'_1<\gamma_1$, one can choose $\theta'\in(0,\theta]$ so that $\theta'-\gamma_1\leq-\gamma_1'$
and for $t\geq 1$, 
$$
\sum_{j\geq 1}2^{j(\beta+\sQ-1-\sA-\theta')}\|\cR^\bba_j\sI_t\div_vF\|_{\bbp'}
\lesssim t^{{\frac{\theta'-\gamma_1}{\alpha}}}\|F\|_{\mL^q_\bbg(\bB^{\beta}_{\bbp;\bba})}
\leq t^{-{\frac{\gamma'_1}{\alpha}}}\|F\|_{\mL^q_\bbg(\bB^{\beta}_{\bbp;\bba})}.
$$
Thus, we further have
$$
\sum_{j\geq 1}2^{j(\beta+\sQ-1-\sA-\theta)}\|\cR^\bba_j\sI_t\div_vF\|_{\bbp'}
\lesssim (1\wedge t)^{\frac{\theta-\gamma_0}{\alpha}}(1\vee t)^{-{\frac{\gamma'_1}{\alpha}}}\|F\|_{\mL^q_\bbg(\bB^{\beta}_{\bbp;\bba})},
$$
which, together with \eqref{Mz1} and $\theta<\sQ$, yields 
\begin{align*}
\|\sI_t\div_v F\|_{\bB^{\beta+\sQ-1-\sA-\theta,1}_{\bbp';\bba}}
&=\|\cR^\bba_0 \sI_t\div_v F\|_{\bbp'}+\sum_{j\geq 1}2^{j(\beta+\sQ-1-\sA-\theta)}\|\cR^\bba_j\sI_t\div_vF\|_{\bbp'}\\
&\lesssim
(1\wedge t)^{\frac{\theta-\gamma_0}{\alpha}}(1\vee t)^{-\frac{\gamma'_1}{\alpha}}\|F\|_{\mL^q_\bbg(\bB^{\beta}_{\bbp;\bba})}.
\end{align*}
Moreover, by  \eqref{Bz1} we directly have
$$
\|\cR^\bba_j\sI_t\div_vF\|_{\bbp'}\lesssim 
\left\{
\begin{aligned}
&2^{j(\sA+1-\beta-\sQ)}(1\wedge t)^{-\frac{\gamma_0}{\alpha}}(1\vee t)^{-{\frac{\gamma_1}{\alpha}}}\|F\|_{\mL^q_\bbg(\bB^{\beta}_{\bbp;\bba})},\\
&2^{j(\sA+1-\beta)}(1\wedge t)^{\frac{\sQ-\gamma_0}{\alpha}}(1\vee t)^{{\frac{\sQ-\gamma_1}{\alpha}}}\|F\|_{\mL^q_\bbg(\bB^{\beta}_{\bbp;\bba})},
\end{aligned}
\right.
$$
which together with \eqref{Mz1} gives \eqref{Scha0} for $\theta=0$ and $\theta=\sQ$ by definition.
\end{proof}
\br\rm
From the proof, one sees that if $\sA\not=\sQ-1$, then $\vartheta\not=1$ and we can take 
$$
\bbg'=\bbg-\big(\theta, (\sQ-1-\sA)\vee 0\big).
$$
\er
\br\rm
About the regularity index $\beta'=\alpha-\frac{\alpha}{q}-1-\sA-\theta$ in \eqref{Scha}, $\alpha$ is the gain from the kinetic semigroup,
and $(\frac\alpha q,1,\sA,\theta)$ are, respectively, the loss from 
the time integrablility of $f$, the divergence, Sobolev's embedding and  the singularity of small time weight.
\er

\section{Well-posedness of nonlinear kinetic Fokker-Planck equation}
\label{Sec3}

  In this section we aim to establish the global well-posedness to the following Fokker-Planck equation with small initial value:
\begin{align}\label{FPE07}
\p_t u=(\Delta^{\frac\alpha 2}_v-v\cdot\nabla_x)u-\div_v ((b*u)u),\ u(0)=u_0.
\end{align}
Firstly we introduce the definition of weak solutions to the above equation.
\bd
Let $u_0\in\cS'$ and $u\in L^1_{loc}([0,\infty); L^p(\mR^{2d}))$ for some $p\in[1,\infty]$. 
We call $u$ a  weak solution of \eqref{FPE07}  with initial value $u_0$ if 
$$
(b*u)u\in L^1_{loc}([0,\infty);\cS'),
$$ 
and for any $\varphi\in C^\infty_c(\mR^{2d})$ and $t>0$,
\begin{align}\label{EST2}
\<u_t,\varphi\>=\<u_0,\varphi\>+\int_0^t\left[\<u_s,(\Delta_v^{\alpha/2}+v\cdot\nabla_x)\varphi\>+\<(b*u_s)u_s,\nabla_v\varphi\>\right]\dif s,
\end{align}
where $\<u_t,\varphi\>:=\int_{\mR^{2d}}u_t(z)\varphi(z)\dif z=\int u_t\varphi$.
\ed
\br\rm
For $\alpha\in(0,2)$, note that up to a constant $c_{d,\alpha}>0$ (see \cite{Sato}),
\begin{align}\label{Mb3}
\Delta^{\alpha/2}_vu(v)=\lim_{\eps\downarrow0}\int_{|w|\geq\eps}\frac{u(v+w)-u(v)}{|w|^{d+\alpha}}\dif w=\frac12
\int_{\mR^d}\frac{u(v+w)+u(v-w)-2u(v)}{|w|^{d+\alpha}}\dif w.
\end{align}
From this representation, it is easy to see that for $\varphi\in C^\infty_c(\mR^{2d})$, 
$$
\Delta^{\alpha/2}_v\varphi\in (L^1\cap L^\infty)(\mR^{2d}).
$$
In particular, $\int^t_0\<u_s,\Delta^{\alpha/2}_v\varphi\>\dif s$ is well-defined for  $u\in L^1_{loc}([0,\infty); L^p(\mR^{2d}))$.
Moreover, by Duhamel's formula, \eqref{EST2} is equivalent to the following integral equation
\begin{align}\label{AB00}
u_t=P_tu_0-\int^t_0P_{t-s}\div_v ((b_s*u_s)u_s)\dif s.
\end{align}
\er
\subsection{Global well-posedness of small initial value}

In this subsection we prove a basic result about the global well-posedness for \eqref{FPE07} with small initial data
based on the contraction mapping theorem.
First of all, we introduce the main indices and their assumptions.
\begin{enumerate}[{\bf (H$_{0}$)}]
\item Let $\alpha\in{(1,2]}$, $q_b\in(\frac{\alpha}{\alpha-1},\infty]$
and $\bbp_0,\bbrho_0,\bbrho_1\in[1,\infty]^2$ with $\bbrho_0\leq\bbrho_1$ and $\bbb1\leq\tfrac1{\bbp_0}+\tfrac1{\bbrho_0}$. Let $\bar\bbrho_0,\bar\bbrho_1\in[1,\infty]^2$ be defined by
$\tfrac1{\bar\bbrho_0}+\tfrac1{\bbrho_0}=\bbb1$ and $\tfrac1{\bar\bbrho_1}+\tfrac1{\bbrho_1}=\bbb1.$
Let $\bbp_1:=\bbp_0\wedge\bar{\boldsymbol{\varrho}}_1$ and
$$
\sQ:=\alpha-\tfrac\alpha{q_b},\ \ \sA_i:=\bba\cdot(\tfrac d{\bbp_i}+\tfrac d{\bbrho_i}-\bbd)=\bba\cdot(\tfrac d{\bbp_i}-\tfrac d{\bar\bbrho_i}),\ i=0,1.
$$
Let $\beta_0\in(-\sQ,0)$ and $\beta_b\in\mR$. Suppose that $\beta_0\not=\sA_i-\beta_b, i=0,1$ and
\begin{align}\label{AD066} 
0<\sA_0-\beta_0-\beta_b\leq\sQ-1\ \mbox{ and }\ \sA_0-2\beta_0-\beta_b<\sQ.
\end{align}
and  $\bbg_{b}:=(\gamma_{b0},\gamma_{b1})$,
\begin{align}\label{KB}
\gamma_{b0}=0,\ \ \gamma_{b1}\in((\sQ-1-\sA_1)\vee0,(\sQ-\sA_1)\vee1).
\end{align}
\end{enumerate}
\br\rm
Note that $(\beta_0,\bbp_0,\bbp_1)$ is referred to the regularity of initial data $u_0$ and $(q_b,\gamma_b,\beta_b,\bbrho_0,\bbrho_1)$
is referred to the regularity of kernel function $b$. Since $\bbrho_0\leq\bbrho_1$ and $\frac1{\bbp_0}\geq\bbb1-\frac1{\bbrho_0}=\frac1{\bar\bbrho_0}$, we have
$$
\tfrac 1{\bbp_1}-\tfrac 1{\bar\bbrho_1}=\tfrac 1{\bbp_0\wedge\bar\bbrho_1}-\tfrac 1{\bar\bbrho_1}\leq\tfrac d{\bbp_0}-\tfrac d{\bar\bbrho_0}\Rightarrow\sA_1\leq\sA_0.
$$
Moreover,  condition \eqref{AD066} is equivalent to the following conditions:
$$
\left\{
\begin{aligned}
&0<\sA_0-\beta_0-\beta_b\leq\sQ-1,\ &\beta_0\in(-1,0),\\
&0<\sA_0-\beta_0-\beta_b<\sQ+\beta_0,\ &\beta_0\in(-\sQ,-1].
\end{aligned}
\right.
$$
\iffalse
The assumption that $\beta_0\notin\{0,\sA_i-\beta_b,i=0,1\}$ is not essential and purely technical since we can always choose
different $\beta'_0\in(-\sQ,\beta_0)$ such that $\beta'_0\notin\{0,\sA_i-\beta_b,i=0,1\}$ and
$$
\sA_0-\beta'_0-\beta_b<\sQ+\beta'_0\wedge(-1).
$$
\fi
\er

For simplicity of notations, in the following we shall use the following parameter set
\begin{align}\label{Pa1}
\Theta:=\big(\alpha,d,\beta_0,\bbp_0,q_b,\gamma_b,\beta_b,\bbrho_0,\bbrho_1\big).
\end{align}
Recalling $\bB^{\beta}_{\bbp;\bba}:=\bB^{\beta,\infty}_{\bbp;\bba}$, with the parameters in {\bf (H$_0$)},
we introduce three basic spaces for later use: for $\beta\geq 0$,
\begin{align}\label{MY}
\mX_\beta:=\bigcap_{i=0,1}\bB^{\beta-\beta_b,1}_{\bar\bbrho_i;\bba},\ \
\mY_\beta:=\bigcap_{i=0,1}\bB^{\beta,1}_{\bbp_i;\bba},\ \ 
\mZ_\beta:=\bigcap_{i=0,1}\bB^{\beta}_{\bbp_i;\bba}.
\end{align}
Note that the dual space of $\mX_\beta$ is given by 
$$
\mX_\beta'=\bB^{\beta_b-\beta}_{\bbrho_0;\bba}+\bB^{\beta_b-\beta}_{\bbrho_1;\bba}.
$$
More precisely, for any $b\in\mX_\beta'$, it can be written as
$$
b=b_0+b_1,\ \ b_i\in\bB^{\beta_b-\beta}_{\bbrho_i;\bba}, i=0,1.
$$

By Lemma \ref{Le215}, we have the following estimates about the kinetic semigroup $P_t$ in $\mX_\beta$, $\mY_\beta$ and $\mZ_\beta$.
\bl
Under {\bf (H$_0$)}, for any $\beta\geq0$, there is a constant $C=C(\beta,\Theta)>0$ such that
for all $t>0$,
\begin{align}\label{AD646}
\|P_tf\|_{\mX_\beta}&\lesssim_C (1\wedge t)^{-\frac{\beta+\sA_0-\beta_0-\beta_b}\alpha}\|f\|_{\mZ_{\beta_0}},\ \ 
\|P_tf\|_{\mY_\beta}\lesssim_C (1\wedge t)^{-\frac{\beta-\beta_0}\alpha}\|f\|_{\mZ_{\beta_0}},
\end{align}
and for all $t\geq 1$,
\begin{align}\label{AD756}
\|P_tf\|_{\mX_\beta}\lesssim_Ct^{-\frac{\sA_1}\alpha}\|f\|_{\mZ_{\beta_0}},\ \
\|P_tf\|_{\mY_\beta}\lesssim_C\|f\|_{\mZ_{\beta_0}}.
\end{align}
Moreover, we also have
\begin{align}\label{AD856}
\|P_tf\|_{\mX_\beta}\lesssim_{C}(1\wedge t)^{-\frac{\beta}{\alpha}}\|f\|_{\mX_0},\ \
\|P_tf\|_{\mY_\beta}\lesssim_{C}(1\wedge t)^{-\frac{\beta}{\alpha}}\|f\|_{\mY_0}.
\end{align}
\el
\begin{proof}
For $i=0,1$, since $\beta_0\not=0, \sA_i-\beta_b$, by  \eqref{AD0446}, it  is easy to see that
$$
\|P_t f\|_{\bB^{\beta,1}_{\bbp_i;\bba}}
\lesssim(1\wedge t)^{-\frac{\beta-\beta_0}\alpha}\|f\|_{\bB^{\beta_0}_{\bbp_i;\bba}},
$$
and for any $0\leq\beta\not=\beta_0+\beta_b-\sA_i$, 
$$
\|P_t f\|_{\bB^{\beta-\beta_b,1}_{\bar\varrho_i;\bba}}
\lesssim(1\wedge t)^{-\frac{\sA_i+\beta-\beta_b-\beta_0}\alpha}\|f\|_{\bB^{\beta_0}_{\bbp_i;\bba}}
\leq(1\wedge t)^{-\frac{\beta+\sA_0-\beta_0-\beta_b}\alpha}\|f\|_{\bB^{\beta_0}_{\bbp_i;\bba}},
$$
where the second inequality is due to $\sA_1\leq\sA_0$.
Moreover, by interpolation inequality \eqref{Sob}, the above inequality holds for all $\beta\geq 0$.
Thus we get \eqref{AD646}. Similarly, estimate \eqref{AD756} follows by \eqref{AD0456} and $\sA_1\leq\sA_0$.
Estimate \eqref{AD856} follows by \eqref{AD335}.
\end{proof}
Now we introduce the following bilinear form for later use:
 \begin{align}\label{Bi0}
 \cH^b_t(u,w)(x,v):=\int^t_0P_{t-s}\div_v ((b_s*u_s)w_s)(x,v)\dif s.
 \end{align}
Note that by the change of variables and the semigroup property of $(P_t)_{t\geq 0}$,
\begin{align}\label{SH2}
\cH^b_{t+s}(u,w)=P_s\cH^b_{t}(u,w)+\cH^{b^t}_{s}(u^t,w^t),\ \ s,t>0,
\end{align}
where for a function $f$,
$$
f^t(s)=f(t+s).
$$
The following lemma provides important estimates about the above bilinear form,
which are consequences of Theorem \ref{Th26} and Lemma \ref{lemB1}.
\bl\label{lem240}
Suppose {\bf (H$_{0}$)}.
Let $\theta\in(0,\sQ)$, $\bbp,\bbp'\in[1,\infty]^2$ with $\bbp\leq\bbp'$ and write $\sA:=\bba\cdot(\tfrac d{\bbp}-\tfrac d{\bbp'})$.
For any $\beta\geq0$ and  $\bbg,\bar\bbg,\bbg',\bbg_b\in\mR^2$ satisfying $\bbg+\bar\bbg+\bbg_b\in[0,\sQ)^2$ and
$$
\bbg'\wle\bbg+\bar\bbg+\bbg_b-(\theta, (\sQ-1-\sA)\vee0),
$$
there is a constant $C=C(\beta,\theta,\bbg,\bar\bbg,\bbg_b,\bbg',\Theta)>0$ such that for all $b\in\mL^{q_b}_{\bbg_b}(\mX_0')$,
$u\in{\mL^\infty_{\bar\bbg}(\mX_\beta)}$ and $w\in{\mL^\infty_{\bbg}(\mL^{\bbp}\cap \bB^{\beta}_{\bbp;\bba})}$,
\begin{align}\label{AS1}
\begin{split}
&\|\cH^b(u,w)\|_{\mL^\infty_{\bbg'}(\bB^{\beta+\beta',1}_{\bbp';\bba})}\lesssim_C
\|b\|_{\mL^{q_b}_{\bbg_b}(\mX_0')}\Big[\|u\|_{\mL^\infty_{\bar\bbg}(\mX_0)}\|w\|_{\mL^\infty_{\bbg}(\mL^{\bbp})}\b1_{\beta=0}\\
&\qquad+\big(\|u\|_{\mL^\infty_{\bar\bbg}(\mX_\beta)}
\|w\|_{\mL^\infty_{\bbg}(\mL^{\bbp})}+\|u\|_{\mL^\infty_{\bar\bbg}(\mX_0)}
\|w\|_{\mL^\infty_{\bbg}(\bB^{\beta}_{\bbp;\bba})}\big)\b1_{\beta>0}\Big],
\end{split}
\end{align}
where $\beta':=\sQ-1-\sA-\theta$.
\el
\begin{proof}
For any $\beta\geq 0$, by  Theorem \ref{Th26}, we have
\begin{align}\label{AS980}
\|\cH^b(u,w)\|_{\mL^\infty_{\bbg'}(\bB^{\beta+\beta',1}_{\bbp';\bba})}
\lesssim\|(b*u)w\|_{\mL^{q_b}_{\bbg+\bar\bbg+\bbg_b}(\bB^{\beta}_{\bbp;\bba})}.
\end{align}
Since $b=b_0+b_1\in \mX_0'=\bB^{\beta_b}_{\bbrho_0;\bba}+\bB^{\beta_b}_{\bbrho_1;\bba}$,
by \eqref{AB2} and \eqref{Con}, we have for $i=0,1$,
\begin{align}\label{SE2}
\begin{split}
\|(b_i*u)w\|_{\bB^{0}_{\bbp;\bba}}
&\lesssim\|(b_i*u)w\|_{\mL^{\bbp}}
\leq\|b_i*u\|_{\mL^\infty}\|w\|_{\mL^{\bbp}}\\
&\lesssim\|b_i*u\|_{\bB^{0,1}_{\infty;\bba}}\|w\|_{\mL^{\bbp}}
\lesssim\|b_i\|_{\bB^{\beta_b}_{\bbrho_i;\bba}}\|u\|_{\bB^{-\beta_b,1}_{\bar\bbrho_i;\bba}}
\|w\|_{\mL^{\bbp}}.
\end{split}
\end{align}
Hence,
$$
\|(b*u)w\|_{\bB^{0}_{\bbp;\bba}}
\lesssim\|b\|_{\mX_0'}\|u\|_{\mX_0}\|w\|_{\mL^{\bbp}}
$$
and
\begin{align}\label{Mn1}
\|(b*u)w\|_{\mL^{q_b}_{\bbg+\bar\bbg+\bbg_b}(\bB^{0}_{\bbp;\bba})}\lesssim\|b\|_{\mL^{q_b}_{\bbg_b}(\mX_0')}\|u\|_{\mL^\infty_{\bar\bbg}(\mX_0)}
\|w\|_{\mL^\infty_{\bbg}(\mL^{\bbp})}.
\end{align}
For $\beta>0$, by (iii) and (ii) of Lemma \ref{lemB1}, we have for $i=0,1$,
\begin{align*}
\|(b_i*u)w\|_{\bB^{\beta}_{\bbp;\bba}}
&\lesssim\|b_i*u\|_{\bB^{\beta}_{\infty;\bba}}\|w\|_{\mL^{\bbp}}
+\|b_i*u\|_{\mL^\infty}\|w\|_{\bB^{\beta}_{\bbp;\bba}}\\
&\lesssim\|b_i\|_{\bB^{\beta_b}_{\bbrho_i;\bba}}
\Big(\|u\|_{\bB^{\beta-\beta_b,1}_{\bar\bbrho_i;\bba}}\|w\|_{\mL^{\bbp}}+\|u\|_{\bB^{-\beta_b,1}_{\bar\bbrho_i;\bba}}\|w\|_{\bB^{\beta}_{\bbp;\bba}}\Big).
\end{align*}
Hence,
$$
\|(b*u)w\|_{\bB^{\beta}_{\bbp;\bba}}
\lesssim \|b\|_{\mX_0'}\Big(\|u\|_{\mX_\beta}\|w\|_{\mL^{\bbp}}+\|u\|_{\mX_0}\|w\|_{\bB^{\beta}_{\bbp;\bba}}\Big)
$$
and
\begin{align}\label{Mn2}
\|(b*u)w\|_{\mL^{q_b}_{\bbg+\bar\bbg+\bbg_b}(\bB^{\beta}_{\bbp;\bba})}
\lesssim \|b\|_{\mL^{q_b}_{\bbg_b}(\mX_0')}
\Big(\|u\|_{\mL^\infty_{\bar\bbg}(\mX_\beta)}
\|w\|_{\mL^\infty_{\bbg}(\mL^{\bbp})}+\|u\|_{\mL^\infty_{\bar\bbg}(\mX_0)}\|w\|_{\mL^\infty_{\bbg}(\bB^{\beta}_{\bbp;\bba})}\Big).
\end{align}
Substituting \eqref{Mn1} and \eqref{Mn2} into \eqref{AS980}, we obtain \eqref{AS1}. 
\end{proof}

\br\rm\label{Re219}
If we assume that $-1<\beta_b<0$ and $|b_t(x,y)|\leq K_t(x-y)$ for some $K\in \mL^{q_b}_{\bbg_b}(\mX_0')$, and $b*u$ is replaced by $b\circledast u=\int_{\mR^d} b(x,y)u(y)\dif y$, then
\eqref{AS1} still holds for $\beta=0$.
Indeed, for the corresponding term $b\circledast u$ in \eqref{SE2}, we can estimate it as follows:
$$
\|b\circledast u\|_{\mL^\infty}\leq\|K*|u|\|_{\mL^\infty}\lesssim \sum_{i=0,1}\|K_i*|u|\|_{\bB^{0,1}_{\infty;\bba}}
\lesssim\sum_{i=0,1}\|K_i\|_{\bB^{\beta_b}_{\bbrho_i;\bba}}\||u|\|_{\bB^{-\beta_b,1}_{\bar\bbrho_i;\bba}}.
$$
Since $-1<\beta_b<0$, by characterization \eqref{CH1} and $||a|-|b||\leq|a-b|$, it is easy to see that
$$
\||u|\|_{\bB^{-\beta_b,1}_{\bar\bbrho_i;\bba}}\lesssim \|u\|_{\bB^{-\beta_b,1}_{\bar\bbrho_i;\bba}}.
$$
Hence, 
$$
\|b\circledast u\|_{\mL^\infty}\lesssim \|K\|_{\mX'_0}\|u\|_{\mX_0}.
$$
\er
\br\rm
It shall be seen below that estimate \eqref{AS1} for $\beta=0$ is enough to show 
the existence and uniqueness of a weak solution for equation \eqref{AB00}. While,
estimate \eqref{AS1} for $\beta>0$ shall be used to show the smoothness of solutions.
\er

The following lemma is a direct application of Lemma \ref{lem240}.
\bl\label{Cor22}
Suppose {\bf (H$_{0}$)}.
Let $\beta\geq 0$ and $\theta\in(0,\sQ)$. Let $b\in\mL^{q_b}_{\bbg_b}(\mX_0')$ and $\bbg,\bar\bbg\in\mR^2$ satisfy 
$$
\bbg+\bar\bbg+\bbg_b\in[0,\sQ)^2.
$$
\begin{enumerate}[(i)]
\item Let $\beta':=\sQ-1-\theta$. For any $\bbg'\in\mR^2$ with
$$
\bbg'\wle\bbg+\bar\bbg+\bbg_b-(\theta, \sQ-1),
$$
there is a $C=C(\Theta,\beta,\theta,\bbg,\bar\bbg,\bbg')>0$ such that for all 
$u\in{\mL^\infty_{\bar\bbg}(\mX_\beta)}$ and $w\in\mL^\infty_{\bbg}(\mY_\beta)$,
\begin{align}\label{AS12}
\|\cH^b(u,w)\|_{\mL^\infty_{\bbg'}(\mY_{\beta+\beta'})}\lesssim_C
\|b\|_{\mL^{q_b}_{\bbg_b}(\mX_0')}\big(\|u\|_{\mL^\infty_{\bar\bbg}(\mX_\beta)}
\|w\|_{\mL^\infty_{\bbg}(\mY_0)}+\|u\|_{\mL^\infty_{\bar\bbg}(\mX_0)}
\|w\|_{\mL^\infty_{\bbg}(\mY_\beta)}\big).
\end{align}
\item Let $\beta':=\sQ-1-\sA_0+\beta_b-\theta$. For any $\bbg'\in\mR^2$ with
$$
\bbg'\wle\bbg+\bar\bbg+\bbg_b-(\theta, (\sQ-1-\sA_1)\vee0),
$$
there is a  $C=C(\Theta,\beta,\theta,\bbg,\bar\bbg,\bbg')>0$ such that for all 
 $u\in{\mL^\infty_{\bar\bbg}(\mX_\beta)}$ and $w\in\mL^\infty_{\bbg}(\mY_\beta)$,
\begin{align}\label{AS14}
\|\cH^b(u,w)\|_{\mL^\infty_{\bbg'}(\mX_{\beta+\beta'})}\lesssim_C
\|b\|_{\mL^{q_b}_{\bbg_b}(\mX_0')}\big(\|u\|_{\mL^\infty_{\bar\bbg}(\mX_\beta)}
\|w\|_{\mL^\infty_{\bbg}(\mY_0)}+\|u\|_{\mL^\infty_{\bar\bbg}(\mX_0)}
\|w\|_{\mL^\infty_{\bbg}(\mY_\beta)}\big).
\end{align}
\end{enumerate}
\el
\begin{proof}
(i) For $i=0,1$, by \eqref{AS1} with $(\bbp',\bbp)=(\bbp_i,\bbp_i)$, we have
\begin{align*}
\|\cH^b(u,w)\|_{\mL^\infty_{\bbg'}(\bB^{\beta+\beta',1}_{\bbp_i;\bba})}
\lesssim\|b\|_{\mL^{q_b}_{\bbg_b}(\mX_0')}
\Big(\|u\|_{\mL^\infty_{\bar\bbg}(\mX_\beta)}\|w\|_{\mL^\infty_{\bbg}(\bB^{0,1}_{\bbp_i})}
+\|u\|_{\mL^\infty_{\bar\bbg}(\mX_0)}\|w\|_{\mL^\infty_{\bbg}(\bB^{\beta,1}_{\bbp_i})}\Big).
\end{align*}
Thus we obtain \eqref{AS12} by the definition of $\mY_\beta$ in \eqref{MY}.

(ii) For $i=0,1$, noting that by $\sA_1\leq\sA_0$,
$$
\sQ-1-\sA_1-\theta\geq\sQ-1-\sA_0-\theta=\beta'-\beta_b,
$$ 
and 
$$
\bbg'\wle\bbg+\bar\bbg+\bbg_b-(\theta, (\sQ-1-\sA_i)\vee0),
$$
by \eqref{AS1} with $(\bbp',\bbp)=(\bar\bbrho_i,\bbp_i)$, we have
\begin{align*}
&\|\cH^b(u,w)\|_{\mL^\infty_{\bbg'}(\bB^{\beta+\beta'-\beta_b,1}_{\bar\bbrho_i;\bba})}
\lesssim\|\cH^b(u,w)\|_{\mL^\infty_{\bbg'}(\bB^{\beta+\sQ-1-\sA_i-\theta,1}_{\bar\bbrho_i;\bba})}\\
&\quad\lesssim\|b\|_{\mL^{q_b}_{\bbg_b}(\mX_0')}
\Big(\|u\|_{\mL^\infty_{\bar\bbg}(\mX_\beta)}\|w\|_{\mL^\infty_{\bbg}(\bB^{0,1}_{\bbp_i})}
+\|u\|_{\mL^\infty_{\bar\bbg}(\mX_0)}\|w\|_{\mL^\infty_{\bbg}(\bB^{\beta,1}_{\bbp_i})}\Big),
\end{align*}
which implies \eqref{AS14}  by  the definition of $\mX_\beta$ in \eqref{MY}.
\end{proof}

Before stating and proving our main result, under {\bf (H$_0$)}, we introduce the following time-weighted spaces where the solutions 
of equation \eqref{FPE07} shall stay: for $\beta\geq 0$,
\begin{align}\label{AS59}
\mS_\beta:=\mL^\infty_{\bbl_\beta}(\mX_{\beta})\cap \mL^\infty_{\bbg_\beta}(\mY_\beta),
\end{align}
where 
\begin{align}\label{AS58}
\bbl_\beta:=\big(\beta+\sA_0-\beta_0-\beta_b, (\sQ-1)\wedge\sA_1\big),\ \ \bbg_\beta:=\(\beta-\beta_0,0\).
\end{align}
It should be noticed that the norm $\|\cdot\|_{\mS_\beta}$ does not have the monotonicity with respect to $\beta$, i.e.,
$$
\|f\|_{\mS_0}\not\leq C\|f\|_{\mS_\beta}, \ \ \beta>0.
$$
But the following interpolation inequality holds: for any $\theta\in[0,1]$ and $\beta>0$,
\begin{align}\label{MX4}
\|f\|_{\mS_{\theta\beta}}\lesssim \|f\|^{1-\theta}_{\mS_{0}}\|f\|^\theta_{\mS_\beta}.
\end{align}
Indeed, by the interpolation inequality \eqref{Sob}, we have $\|f\|_{\mX_{\theta\beta}}\lesssim \|f\|_{\mX_\beta}^\theta\|f\|_{\mX_0}^{1-\theta}$  and
\begin{align*}
\|f\|_{\mL^\infty_{\bbl_{\theta\beta}}(\mX_{\theta\beta})}
&=\sup_{t>0}\left(\|f\|_{\mX_{\theta\beta}}(1\wedge t)^{\theta\beta+\sA_0-\beta_0-\beta_b}(1\vee t)^{(\sQ-1)\wedge\sA_1}\right)\\
&\lesssim \sup_{t>0}\left(\|f\|_{\mX_{\beta}}(1\wedge t)^{\beta+\sA_0-\beta_0-\beta_b}(1\vee t)^{(\sQ-1)\wedge\sA_1}\right)^\theta\\
&\quad\times\sup_{t>0}\left(\|f\|_{\mX_0}(1\wedge t)^{\sA_0-\beta_0-\beta_b}(1\vee t)^{(\sQ-1)\wedge\sA_1}\right)^{1-\theta}\\
&=\|f\|^\theta_{\mL^\infty_{\bbl_{\beta}}(\mX_{\beta})}\|f\|^{1-\theta}_{\mL^\infty_{\bbl_0}(\mX_0)}.
\end{align*}
Similarly, 
$$
\|f\|_{\mL^\infty_{\bbg_{\theta\beta}}(\mY_{\theta\beta})}\lesssim \|f\|^\theta_{\mL^\infty_{\bbg_{\beta}}(\mY_{\beta})}\|f\|^{1-\theta}_{\mL^\infty_{\bbg_0}(\mY_0)}.
$$
Combining the above two estimates, we get \eqref{MX4}.

Let 
$$
\Lambda:=(\sQ-1)\wedge(\sQ-1-\sA_0+\beta_b)>0.
$$
Now we use Lemma \ref{Cor22} to show  the following crucial lemma.
\bl\label{Le38}
Under {\bf (H$_{0}$)}, for any $\beta\geq 0$ and $\beta'\in[0,\Lambda)$, there is a constant $C=C(\beta',\beta,\Theta)>0$ such that for all $b\in\mL^{q_b}_{\bbg_b}(\mX_0')$ and $u,w\in\mS_0\cap\mS_\beta$,
\begin{align}\label{AD122}
    \|\cH^{b}(u,w)\|_{\mS_{\beta+\beta'}}\lesssim_C\|b\|_{\mL^{q_b}_{\bbg_b}(\mX_0')}\big(\|u\|_{\mS_\beta}+\|u\|_{\mS_0}\big)\big(\|w\|_{\mS_\beta}+\|w\|_{\mS_0}\big)=:\sU,
\end{align}
where $\mS_\beta$ is defined in \eqref{AS59}.
\el
\begin{proof}
We divide the proof into three steps.

({\it Step 1}) In this step we first show \eqref{AD122} for $\beta=\beta'=0$.
Observe that
$$
\gamma_{b1}\in((\sQ-1-\sA_1)\vee0,(\sQ-\sA_1)\vee1)
=(\sQ-1-(\sQ-1)\wedge\sA_1,\sQ-(\sQ-1)\wedge\sA_1).
$$ 
By the definition \eqref{AS58} and the assumption \eqref{AD066}, it is easy to see that
$$
\bbg_0+\bbl_0+\bbg_b=(\sA_0-2\beta_0-\beta_b,(\sQ-1)\wedge\sA_1+\gamma_{b1})
\in[0,\sQ)^2.
$$
Since $\bbg_0\wle\bbg_0+\bbl_0+\bbg_b-(\sQ-1, \sQ-1)$, 
by (i) of Lemma \ref{Cor22} with the choice of $\theta=\sQ-1$, $\beta=0$ and $\bbg=\bbg_0, \bar\bbg=\bbl_0$, 
$\bbg'=\bbg_0$, we have
$$
\|\cH^{b}(u,w)\|_{\mL^\infty_{\bbg_0}(\mY_0)}\lesssim\|b\|_{\mL^{q_b}_{\bbg_b}(\mX_0')}\|u\|_{\mL^\infty_{\bbl_0}(\mX_0)}
\|w\|_{\mL^\infty_{\bbg_0}(\mY_0)},
$$
Moreover, for $\theta:=\sQ-1-\sA_0+\beta_b$, since $\bbl_0\wle\bbg_0+\bbl_0+\bbg_b-(\theta, (\sQ-1-\sA_1)\vee0),$
by (ii) of Lemma \ref{Cor22} with the choice of the above $\theta$, $\beta=0$ and 
$\bbg=\bbg_0$, $\bar\bbg=\bbg'=\bbl_0$, we have
$$
\|\cH^{b}(u,w)\|_{\mL^\infty_{\bbl_0}(\mX_0)}\lesssim\|b\|_{\mL^{q_b}_{\bbg_b}(\mX_0')}\|u\|_{\mL^\infty_{\bbl_0}(\mX_0)}
\|w\|_{\mL^\infty_{\bbg_0}(\mY_0)}.
$$
Hence,
\begin{align}\label{SH3}
\|\cH^{b}(u,w)\|_{\mS_0}=\|\cH^{b}(u,w)\|_{\mL^\infty_{\bbg_0}(\mY_0)}+\|\cH^{b}(u,w)\|_{\mL^\infty_{\bbl_0}(\mX_0)}\lesssim\|b\|_{\mL^{q_b}_{\bbg_b}(\mX_0')}\|u\|_{\mS_0}\|w\|_{\mS_0}.
\end{align}

({\it Step 2}) Let $\beta\geq 0$. For fixed $\beta'\in[0,\sQ-1)$, let
$$
\theta=\sQ-1-\beta'\in(0,\sQ).
$$
By (i) of Lemma \ref{Cor22} with the choice of the above $\theta$ , $\bbg=\bar\bbg=0$ and $\bbg'=(-\theta,-\ell_{01})$, we have
$$
\|\cH^{b^t}_\cdot(u^t,w^t)\|_{\mL^\infty_{\bbg'}(\mY_{\beta+\beta'})}
\lesssim \|b\|_{\mL^{q_b}_{\bbg_b}(\mX_0')}\Big(\|u^t\|_{\mL^\infty_{0}(\mX_{\beta})}
\|w^t\|_{\mL^\infty_{0}(\mY_0)}+\|u^t\|_{\mL^\infty_{0}(\mX_0)}
\|w^t\|_{\mL^\infty_{0}(\mY_\beta)}\Big).
$$
Note that by \eqref{Sh1} and \eqref{AS58},
\begin{align}\label{AS8}
\|u^t\|_{\mL^\infty_{0}(\mX_{\beta})}\|w^t\|_{\mL^\infty_{0}(\mY_0)}
&\lesssim(1\wedge t)^{-\frac{\ell_{\beta0}+\gamma_{00}}\alpha}(1\vee t)^{-\frac{\ell_{01}}\alpha}\|u\|_{\mL^\infty_{\bbl_\beta}(\mX_{\beta})}
\|w\|_{\mL^\infty_{\bbg_0}(\mY_0)},
\end{align}
and
\begin{align}\label{AS9}
\|u^t\|_{\mL^\infty_{0}(\mX_0)}\|w^t\|_{\mL^\infty_{0}(\mY_\beta)}
&\lesssim(1\wedge t)^{-\frac{\ell_{00}+\gamma_{\beta0}}\alpha}(1\vee t)^{-\frac{\ell_{01}}\alpha}\|u\|_{\mL^\infty_{\bbl_0}(\mX_0)}\|w\|_{\mL^\infty_{\bbg_\beta}(\mY_\beta)}.
\end{align}
Thanks to $\ell_{\beta0}+\gamma_{00}=\ell_{00}+\gamma_{\beta0}=\ell_{00}+\beta-\beta_0$, we further have
\begin{align*}
\|\cH^{b^t}_t(u^t,w^t)\|_{\mY_{\beta+\beta'}}
&\lesssim (1\wedge t)^{\frac{\sQ-1-\beta'}\alpha}(1\vee t)^{\frac{\ell_{01}}\alpha}
\|\cH^{b^t}_\cdot(u^t,w^t)\|_{\mL^\infty_{\bbg'}(\mY_{\beta+\beta'})}\\
&\lesssim (1\wedge t)^{\frac{\sQ-1-\beta'-\ell_{00}-\beta+\beta_0}\alpha}\sU
\leq (1\wedge t)^{-\frac{\beta'+\beta-\beta_0}\alpha}\sU,
\end{align*}
where $\sU$ is defined by \eqref{AD122}.
On the other hand, by \eqref{AD856} and \eqref{SH3}, we have
\begin{align*}
\|P_t\cH^b_t(u,w)\|_{\mY_{\beta+\beta'}}&\lesssim (1\wedge t)^{-\frac{\beta+\beta'}{\alpha}}\|\cH^b_t(u,w)\|_{\mY_0}\lesssim  (1\wedge t)^{-\frac{\beta+\beta'+\gamma_{00}}\alpha}\|\cH^b(u,w)\|_{\mL^\infty_{\bbg_0}(\mY_0)}\\
&\lesssim  (1\wedge t)^{-\frac{\beta+\beta'+\gamma_{00}}\alpha}
\|b\|_{\mL^{q_b}_{\bbg_b}(\mX_0')}\|u\|_{\mS_0}\|w\|_{\mS_0}\leq(1\wedge t)^{-\frac{\beta+\beta'+\gamma_{00}}\alpha}\sU.
\end{align*}
Combining the above two estimates and by \eqref{SH2} with $s=t$, we obtain
\begin{align}\label{AS11}
\|\cH^b_{2t}(u,w)\|_{\mY_{\beta+\beta'}} 
\lesssim  (1\wedge t)^{-\frac{\beta+\beta'+\gamma_{00}}\alpha}\sU.
\end{align}

({\it Step 3}) 
Let $\beta\geq 0$. For fixed $\beta'\in[0,\sQ-1-\sA_0+\beta_b)$, let
$$
\theta:=\sQ-1-\sA_0+\beta_b-\beta'\in(0,\sQ).
$$
By (ii) of Lemma \ref{Cor22} with the choice of  the above $\theta$, $\bbg=\bar\bbg=0$ and $\bbg'=(-\theta,0)$, we have
$$
\|\cH^{b^t}_\cdot(u^t,w^t)\|_{\mL^\infty_{\bbg'}(\mX_{\beta+\beta'})}
\lesssim \|b\|_{\mL^{q_b}_{\bbg_b}(\mX_0')}\Big(\|u^t\|_{\mL^\infty_{0}(\mX_{\beta})}
\|w^t\|_{\mL^\infty_{0}(\mY_0)}+\|u^t\|_{\mL^\infty_{0}(\mX_0)}
\|w^t\|_{\mL^\infty_{0}(\mY_\beta)}\Big).
$$
By \eqref{AS8} and \eqref{AS9}, we further have
\begin{align*}
\|\cH^{b^t}_t(u^t,w^t)\|_{\mX_{\beta+\beta'}}
&\lesssim(1\wedge t)^{\frac{\sQ-1-\sA_0+\beta_b-\beta'}\alpha}
\|\cH^{b^t}_\cdot(u^t,w^t)\|_{\mL^\infty_{\bbg'}(\mX_{\beta+\beta'})}\\
&\lesssim (1\wedge t)^{\frac{\sQ-1-\sA_0+\beta_b-\beta'-\ell_{00}-\beta+\beta_0}\alpha}(1\vee t)^{-\frac{\ell_{01}}\alpha}\sU\\
&\lesssim (1\wedge t)^{-\frac{\ell_{00}+\beta+\beta'}\alpha}(1\vee t)^{-\frac{\ell_{01}}\alpha}\sU.
\end{align*}
Moreover, by \eqref{AD856} we have
\begin{align*}
&\|P_t\cH^b_t(u,w)\|_{\mX_{\beta+\beta'}}
\lesssim (1\wedge t)^{-\frac{\beta+\beta'}{\alpha}}
\|\cH^b_t(u,w)\|_{\mX_0}\\
&\qquad\lesssim (1\wedge t)^{-\frac{\beta+\beta'+\ell_{00}}{\alpha}}(1\vee t)^{-\frac{\ell_{01}}\alpha}
\|\cH^b(u,w)\|_{\mL^\infty_{\bbl_0}(\mX_0)}\\
&\qquad\!\!\!\stackrel{\eqref{SH3}}{\lesssim} (1\wedge t)^{-\frac{\beta+\beta'+\ell_{00}}{\alpha}}(1\vee t)^{-\frac{\ell_{01}}\alpha}
\|b\|_{\mL^{q_b}_{\bbg_b}(\mX_0')}\|u\|_{\mS_0}\|w\|_{\mS_0}.
\end{align*}
Hence, by \eqref{SH2} with $s=t$ again,
$$
\|\cH^b_{2t}(u,w)\|_{\mX_{\beta+\beta'}}
\lesssim  (1\wedge t)^{-\frac{\beta+\beta'+\ell_{00}}{\alpha}}(1\vee t)^{-\frac{\ell_{01}}\alpha}\sU,
$$
which together with \eqref{AS11} yields \eqref{AD122}.
\end{proof}

%For $\beta\geq 0$ and $\bbp\in[1,\infty]^2$, we introduce the following Banach space:
%$$\mZ_{\beta;\bba}_{\bbp,\infty}:=\bB^{\beta;\bba}_{\bbp,\infty}\cap\bB^{\beta;\bba}_{1,\infty},$$

\br\rm
In the situation of Remark \ref{Re219}, \eqref{AD122} still holds for $\beta=0$.
\er

Now we can state and prove our main result of this section.
\bt\label{Main}
Under {\bf (H$_{0}$)}, there are two constants $C_0, C_1>0$ only depending on $\Theta$ such that 
for any $b\in \mL^{q_b}_{\bbg_b}(\mX'_0)$ and $u_0\in\mZ_{\beta_0}$ satisfying 
\begin{align}\label{Small}
\|u_0\|_{\mZ_{\beta_0}}\|b\|_{\mL^{q_b}_{\bbg_b}(\mX'_0)}\leq C_0,
\end{align}
there exists a unique weak solution to Fokker-Planck equation \eqref{FPE07} such that
\begin{align}\label{AS6}
\|u\|_{\mS_0}\leq C_1\|u_0\|_{\mZ_{\beta_0}}.
\end{align}
Moreover, we have the following conclusions:
\begin{enumerate}[(i)]
\item{\bf (Smoothness)} For any $\beta\geq 0$, there is a constant $C_\beta=C_\beta(\Theta,\|b\|_{\mL^{q_b}_{\bbg_b}(\mX'_0)})>0$ such that
\begin{align}\label{SM9}
\|u\|_{\mS_\beta}\leq C_\beta.
\end{align}
\item {\bf (Stability)} For another $\wt b\in \mL^{q_b}_{\bbg_b}(\mX'_0)$
and $\wt u_0\in\mZ_{\beta_0}$ satisfying \eqref{Small},
let $\wt u$ be the unique solution corresponding to $\wt b$ and $\wt u_0$.  
Then for any $\beta\geq 0$,
there is a constant $C_\beta>0$ only depending on $\beta,\Theta$ and 
$\|b\|_{\mL^{q_b}_{\bbg_b}(\mX'_0)}, \|\wt b\|_{\mL^{q_b}_{\bbg_b}(\mX'_0)}$ such that
\begin{align}\label{AS609}
\|u-\wt u\|_{\mS_\beta}\leq C_\beta \big(\|u_0-\wt u_0\|_{\mZ_{\beta_0}}+\|b-\wt b\|_{\mL^{q_b}_{\bbg_b}(\mX_0')}\big).
\end{align}
\end{enumerate}
\et
\begin{proof}
%We divide the proof into four steps.

We use the contraction mapping theorem to prove the existence and uniqueness of a weak solution.
First of all, by \eqref{AD646} and \eqref{AD756}, we have for any $\beta\geq 0$,
\begin{align}
\|P_\cdot u_0\|_{\mS_\beta}\lesssim\|u_0\|_{\mZ_{\beta_0}}.
\label{AD02}
\end{align}
Thus, by \eqref{AD02} and \eqref{AD122}, we can define a map $\cU: {\mS_0}\to{\mS_0}$ by
\begin{align}\label{AD301}
\cU(u)(t):=P_tu_0{\red-}\cH^{b}_t(u,u)
\end{align}
so that for some constants $C_2, C_3>0$ only depending on $\Theta$ and for all $u\in\mS_0$,
\begin{align}\label{AD31}
\|\cU(u)\|_{{\mS_0}}
\leq C_2\|u_0\|_{\mZ_{\beta_0}}+C_3\|b\|_{\mL^{q_b}_{\bbg_b}(\mX_0')}\|u\|^2_{{\mS_0}},
\end{align}
and  for all $u,\wt u\in\mS_0$,
\begin{align}\label{AD01}
\|\cU(u)-\cU(\wt u)\|_{{\mS_0}}&=\|\cH^{b}(u,u)-\cH^{b}(\wt u,\wt u)\|_{{\mS_0}}
\leq C_3\|b\|_{\mL^{q_b}_{\bbg_b}(\mX_0')}\Big(\|u\|_{{\mS_0}}+
\|\wt u\|_{{\mS_0}}\Big)\|u-\wt u\|_{{\mS_0}}.
\end{align}
Now let us choose 
$$
\lambda=4C_2\|u_0\|_{\mZ_{\beta_0}}/3=C_1\|u_0\|_{\mZ_{\beta_0}},
$$
and suppose 
$$
\|u_0\|_{\mZ_{\beta_0}}\|b\|_{\mL^{q_b}_{\bbg_b}(\mX_0')}\leq3/(16 C_3C_2)=:C_0.
$$
Then it is easy to see that
\begin{align}\label{AD41}
C_2\|u_0\|_{\mZ_{\beta_0}}+C_3\|b\|_{\mL^{q_b}_{\bbg_b}(\mX_0')}\lambda^2\leq\lambda,\ \ 
2C_3\|b\|_{\mL^{q_b}_{\bbg_b}(\mX_0')}\lambda\leq 1/2,
\end{align}
and by \eqref{AD31} and \eqref{AD01}, $\cU$ is a contraction mapping on the closed ball 
$$
\cB^\lambda:=\Big\{u\in {{\mS_0}}:\|u\|_{{\mS_0}}\leq\lambda\Big\}.
$$
Therefore, there is a unique fixed point $u\in\cB^\lambda$ so that for all $t>0$,
\begin{align}\label{FF9}
u(t)=\cU(u)(t)=P_tu_0{\red-}\cH^{b}_t(u,u).
\end{align}

({\bf Smoothness}) In this step we show the smoothness of the solution by induction and \eqref{AD122}. 
Suppose that we have shown \eqref{SM9} for some $\beta\geq 0$.
Then by \eqref{AD02} and \eqref{AD122}, we have for any $\beta'\in[0,\Lambda)$,
\begin{align*}
\|u\|_{\mS_{\beta+\beta'}}\leq \|P_\cdot u_0\|_{\mS_{\beta+\beta'}}+\|\cH^b(u,u)\|_{\mS_{\beta+\beta'}}
\lesssim \|u_0\|_{\mZ_{\beta_0}}+\|b\|_{\mL^{q_b}_{\bbg_b}(\mX'_0)}\big(\|u\|_{\mS_\beta}+\|u\|_{\mS_0}\big)^2<\infty.
\end{align*}
 
({\bf Stability}) Let  $\wt u$ be the unique solution corresponding to $(\wt b,\wt u_0)$.
By \eqref{FF9} we have
\begin{align*}
u(t)-\wt u(t)=P_tu_0-P_t\wt u_0{\red-}\cH^{b}_t(u,u){\red+}\cH^{\wt b}_t(\wt u,\wt u).
\end{align*}
By symmetry, we may assume $\|\wt u_0\|_{\mZ_{\beta_0}}\leq \|u_0\|_{\mZ_{\beta_0}}$.
By \eqref{AD01} and \eqref{AD41}, we have
$$
\|\cH^{b}(u,u)-\cH^{b}(\wt u,\wt u)\|_{\mS_0}\leq C_3\|b\|_{\mL^{q_b}_{\bbg_b}(\mX_0')}\Big(\|u\|_{{\mS_0}}+
\|\wt u\|_{{\mS_0}}\Big)\|u-\wt u\|_{{\mS_0}}
\leq\tfrac12\|u-\wt u\|_{{\mS_0}},
$$
and
$$
\|\cH^{b}(\wt u,\wt u)-\cH^{\wt b}(\wt u,\wt u)\|_{\mS_0}
=\|\cH^{b-\wt b}(\wt u,\wt u)\|_{\mS_0}\leq C_3\|b-\wt b\|_{\mL^{q_b}_{\bbg_b}(\mX_0')}\|\wt u\|_{\mS_0}^2.
$$
Thus by \eqref{AD02},
$$
\|u-\wt u\|_{\mS_0}\leq C_2\|u_0-\wt u_0\|_{\mZ_{\beta_0}}+\tfrac12\|u-\wt u\|_{\mS_0}+C_3\|b-\wt b\|_{\mL^{q_b}_{\bbg_b}(\mX_0')},
$$
which implies that
$$
\|u-\wt u\|_{\mS_0}\leq  2C_2\|u_0-\wt u_0\|_{\mZ_{\beta_0}}+2C_3\|b-\wt b\|_{\mL^{q_b}_{\bbg_b}(\mX_0')}.
$$
Next we use induction and \eqref{AD122} to show \eqref{AS609} for any $\beta>0$. Suppose that we have shown \eqref{AS609}
for some $\beta\geq 0$. Then by \eqref{AD02} and \eqref{AD122}, we have for any $\beta'\in[0,\Lambda)$,
\begin{align*}
\|u-\wt u\|_{\mS_{\beta+\beta'}}
&\leq\|P_\cdot(u_0-\wt u_0)\|_{\mS_{\beta+\beta'}}
+\|\cH^b(u,u)-\cH^{\wt b}(\wt u,\wt u)\|_{\mS_{\beta+\beta'}}\\
&\lesssim\|u_0-\wt u_0\|_{\mZ_{\beta_0}}
+\|\cH^b(u, u)-\cH^b(\wt u,\wt u)\|_{\mS_{\beta+\beta'}}+\|\cH^{b-\wt b}(\wt u,\wt u)\|_{\mS_{\beta+\beta'}}\\
&\lesssim\|u_0-\wt u_0\|_{\mZ_{\beta_0}}+\|b\|_{\mL^{q_b}_{\bbg_b}(\mX_0')}\Big(\|u-\wt u\|_{\mS_\beta}\|u\|_{\mS_0}+\|u\|_{\mS_\beta}\|u-\wt u\|_{\mS_0}\\
&+ \|\wt u\|_{\mS_\beta}\|u-\wt u\|_{\mS_0}+\|u-\wt u\|_{\mS_\beta}\|\wt u\|_{\mS_0}\Big)
+\|b-\wt b\|_{\mL^{q_b}_{\bbg_b}(\mX_0')}\|\wt u\|_{\mS_\beta}\|\wt u\|_{\mS_0},
\end{align*}
which,  by the induction hypothesis and \eqref{SM9}, in turn implies that
$$
\|u-\wt u\|_{\mS_{\beta+\beta'}}\lesssim\|u_0-\wt u_0\|_{\mZ_{\beta_0}}+\|b-\wt b\|_{\mL^{q_b}_{\bbg_b}(\mX_0')}.
$$
Thus we complete the proof.
 \end{proof}
\br\rm\label{Re35}
Estimate \eqref{AS609} provides the continuous dependence of the solutions with respect to the kernel function $b$ and 
initial value $u_0$. 
We assume $b=b_0+b_1\in\mL^{q_b}_0(\mX'_0)$, 
where $b_0$ has low integrability and $b_1$ has high integrability, which is motivated by example \eqref{HH2}.
Moreover, we require $\gamma_b>(\sQ-1-\sA_1)\vee0$, which means that $b$ has some time decay property as 
$t\to\infty$. This of course rules out the {\it time-independent} $b$.
However, when considering the finite time interval, we can drop this assumption. Indeed,
for fixed $T\geq 1$, by definition it is easy to see that
\begin{align}\label{MX8}
\|b\b1_{[0,T]}\|_{\mL^{q_b}_{\bbg_b}(\mX'_0)}
\leq T^{\frac{\gamma_b}\alpha}\|b\|_{\mL^{q_b}_0(\mX'_0)}.
\end{align}
Thus, we can use Theorem \ref{Main} to $b\b1_{[0,T]}$. In the next subsection,
we shall consider the large time decay estimate under extra assumptions about the initial value $u_0$ and integrability index $\varrho_1$.
\er
\br\rm\label{Re313}
Remark \ref{Re219} is applicable to Theorem \ref{Main} except that \eqref{SM9} and \eqref{AS609} only hold for $\beta\in[0,\Lambda)$
since \eqref{AD122} still holds for $\beta=0$ under the assumptions of Remark \ref{Re219}. 
\er

\subsection{Large time decay estimates}\label{Sub33}
In this subsection we drop the condition $\gamma_b>(\sQ-1-\sA_1)\vee0$ by assuming extra assumptions about the initial value 
and  integrability index $\varrho_1$. 
More precisely, we have

\bt\label{thm230}
In addition {\bf (H$_{0}$)}, we assume $\bba\cdot\frac d{\bbrho_1}>\sQ-1$, $\bbp_0\not=\bbb1$ and $b\in\mL^{q_b}_0(\mX_0')$. 
For any $u_0\in\mZ_{\beta_0}\cap \bB^{\beta_0}_{\bbb1;\bba}$,
there are two constants $\delta\geq 1$ and $C_0>0$ only depending on $\Theta$ such that if
\begin{align}\label{Mb10}
\|u_0\|_{\mZ_{\beta_0}}\|b\|_{\mL^{q_b}_0(\mX_0')}
\Big(\big(\|b\|_{\mL^{q_b}_0(\mX_0')}\|u_0\|_{\bB^{\beta_0}_{\bbb1;\bba}}\big)\vee 1\Big)^\delta
\leq C_0,
\end{align}
then there is a unique global weak solution  $u$ to FPE \eqref{FPE07} such that
\begin{align}\label{SM219}
\sup_{t>0}\left((1\wedge t)^{-\frac{\beta_0}\alpha}\big(\|u(t)\|_{\mY_0}+\|u(t)\|_{\bB^{0,1}_{\bbb1;\bba}}\big)
+(1\wedge t)^{\frac{\sA_0-\beta_0-\beta_b}\alpha}\|u(t)\|_{\mX_0}\right)<\infty.
\end{align}
Moreover, we have the following conclusions:
\begin{enumerate}[(i)]
\item{\bf (Smoothness)} For any $\beta\geq 0$, it holds that
\begin{align}\label{SM29}
\sup_{t>0}\left((1\wedge t)^{\frac{\beta-\beta_0}\alpha}\big(\|u(t)\|_{\mY_\beta}+\|u(t)\|_{\bB^{\beta,1}_{\bbb1;\bba}}\big)+(1\wedge t)^{\frac{\beta+\sA_0-\beta_0-\beta_b}\alpha}\|u(t)\|_{\mX_\beta}\right)
<\infty.
\end{align}
\item {\bf (Stability)} For another $\wt b\in\mL^{q_b}_0(\mX_0')$ and $\wt u_0\in\mZ_{\beta_0}$ satisfying the same condition \eqref{Mb10},
let $\wt u$ be the unique solution corresponding to $\wt b$ and $\wt u_0$.  For any $\beta\geq 0$,
there is a constant $C_\beta>0$ only depending on $\Theta,\|b\|_{\mL^{q_b}_0(\mX'_0)},\|\wt b\|_{\mL^{q_b}_0(\mX'_0)}$ 
such that
\begin{align}\label{SM39}
\sup_{t>0}\left((1\wedge t)^{\frac{\beta-\beta_0}\alpha}\|u(t)-\wt u(t)\|_{\mY_\beta}\right)
\leq C_\beta \big(\|u_0-\wt u_0\|_{\mZ_{\beta_0}}+\|b-\wt b\|_{\mL^{q_b}_{0}(\mX_0')}\big).
\end{align}
\item {\bf (Large time decay estimate)} For any $\beta\geq 0$ and $i=0,1$, it holds that
\begin{align}\label{aAA51}
\sup_{t\ge1}\left(t^{\frac{\bba\cdot(\bbd-d/\bbp_i)}\alpha}\|u(t)\|_{\bB^{\beta,1}_{\bbp_i;\bba}}\right)<\infty.
\end{align}
\end{enumerate}
\et
\begin{proof}%[Proof of Theorem \ref{thm230}]
The proof is similar to Theorem \ref{Main}. But we need to carefully treat the large time decay estimates. 
Since $\bba\cdot\frac d{\bbrho_1}>\sQ-1$ and $\bbp_0\not=1$,
$$
\sQ-1-(\sQ-1)\wedge\sA_1=(\sQ-1-\sA_1)\vee 0<\bba\cdot\tfrac{d}{\bbrho_1}-\sA_1=\bba\cdot(\bbd-\tfrac{d}{\bbp_1}),
$$ 
one can fix two large time decay weights $\gamma_1,\ell_1\in[0,\infty)$ with $\gamma_1+\ell_1\in[0,\sQ)$ and
\begin{align}\label{Decay}
\gamma_1\in((\sQ-1-\sA_1)\vee 0,\bba\cdot(\bbd-\tfrac{d}{\bbp_1})),\ \ \ell_1\in(\sQ-1,\bba\cdot\tfrac{d}{\bbrho_1}).
\end{align}
Let the weights $\bbg_\beta$ and $\bbl_\beta$ in the definition \eqref{AS59} of solution space $\mS_\beta$  be replaced by
\begin{align}\label{AS058}
\bbg_\beta:=\(\beta-\beta_0,\gamma_1\),\ \ \bbl_\beta:=\big(\beta+\sA_0-\beta_0-\beta_b, \ell_1\big).
\end{align}
In the following we divide the proof into three steps.

({\it Step 1})  In this step we show that for any $\beta\geq 0$, 
there are constants $C_1,C_2,\delta_0,\delta_1>0$  only depending on $\Theta$ such that for any $T\geq 1$,
\begin{align}\label{AD433}
\|P_\cdot u_0\|_{\mS_\beta}\leq C_1T^{\delta_0}\|u_0\|_{\mZ_{\beta_0}}+C_2 T^{-\delta_1}\|u_0\|_{\bB^{\beta_0}_{\bbb1;\bba}}.
\end{align}
As in proving \eqref{AD756}, by \eqref{AD0456} and $\bbp_1\leq\bbp_0$, we have
\begin{align*}
\sup_{t\geq 1}\left(t^{\frac{\gamma_1}\alpha}\|P_t u_0\|_{\mY_\beta}\right)
&\leq\sup_{t\in[1,T]}\left(t^{\frac{\gamma_1}\alpha}\|P_t u_0\|_{\mY_\beta}\right)
+\sup_{t\geq T}\left(t^{\frac{\gamma_1}\alpha}\|P_t u_0\|_{\mY_\beta}\right)\\
&\lesssim T^{\frac{\gamma_1}\alpha}\|u_0\|_{\mZ_{\beta_0}}
+\sup_{t\geq T}\left(t^{\frac{\gamma_1-\bba\cdot(\bbd-d/\bbp_1)}\alpha}\|u_0\|_{\bB^{\beta_0}_{\bbb1;\bba}}\right)\\
&\leq T^{\frac{\gamma_1}\alpha}\|u_0\|_{\mZ_{\beta_0}}
+T^{\frac{\gamma_1-\bba\cdot(\bbd-d/\bbp_1)}\alpha}\|u_0\|_{\bB^{\beta_0}_{\bbb1;\bba}},
\end{align*}
and by $\bar\bbrho_1\leq\bar\bbrho_0$,
\begin{align*}
\sup_{t\geq 1}\left(t^{\frac{\ell_1}\alpha}\|P_t u_0\|_{\mX_\beta}\right)
&\leq\sup_{t\in[1,T]}\left(t^{\frac{\ell_1}\alpha}\|P_t u_0\|_{\mX_\beta}\right)
+\sup_{t\geq T}\left(t^{\frac{\ell_1}\alpha}\|P_t u_0\|_{\mX_\beta}\right)\\
&\lesssim T^{\frac{\ell_1}\alpha}\|u_0\|_{\mZ_{\beta_0}}
+\sup_{t\geq T}\left(t^{\frac{\ell_1-\bba\cdot d/\bbrho_1}\alpha}\|u_0\|_{\bB^{\beta_0}_{\bbb1;\bba}}\right)\\
&\leq T^{\frac{\ell_1}\alpha}\|u_0\|_{\mZ_{\beta_0}}
+T^{\frac{\ell_1-\bba\cdot d/\bbrho_1}\alpha}\|u_0\|_{\bB^{\beta_0}_{\bbb1;\bba}}.
\end{align*}
Moreover, by \eqref{AD646}, we also have
\begin{align*}
\sup_{t\in(0,1]}\left(t^{\frac{\gamma_{\beta0}}\alpha}\|P_t u_0\|_{\mY_\beta}
+t^{\frac{\ell_{\beta 0}}\alpha}\|P_t u_0\|_{\mX_\beta}\right)
\lesssim\|u_0\|_{\mZ_{\beta_0}}.
\end{align*}
Combining the above estimates and by definition \eqref{AS59}, we obtain \eqref{AD433}.

({\it Step 2}) In this step we show the existence and uniqueness of a global solution under some  smallness condition \eqref{Mb10}.
With the weights \eqref{AS058}, by Lemma \ref{Cor22}, it is completely the same as Lemma \ref{Le38} to derive 
that for any $\beta\geq 0$ and $\beta'\in[0,\Lambda)$,
\begin{align}\label{AD422}
    \|\cH^{b}(u,w)\|_{\mS_{\beta+\beta'}}\lesssim\|b\|_{\mL^{q_b}_0(\mX_0')}\big(\|u\|_{\mS_\beta}+\|u\|_{\mS_0}\big)\big(\|w\|_{\mS_\beta}+\|w\|_{\mS_0}\big).
\end{align}
Thus by \eqref{AD433} and \eqref{AD422}, 
we can define a map $\cU: {\mS_0}\to{\mS_0}$ by
$$
\cU(u)(t):=P_tu_0{\red-}\cH^{b}_t(u,u)
$$
so that for some constants $C_1, C_2, C_3>0$ only depending on $\Theta$ and for all $u\in\mS_0$,
$$
\|\cU(u)\|_{{\mS_0}}
\leq C_1T^{\delta_0}\|u_0\|_{\mZ_{\beta_0}}
+C_2 T^{-\delta_1}\|u_0\|_{\bB^{\beta_0}_{\bbb1;\bba}}
+C_3\|b\|_{\mL^{q_b}_0(\mX_0')}\|u\|^2_{{\mS_0}},
$$
and  for all $u,\wt u\in\mS_0$,
$$
\|\cU(u)-\cU(\wt u)\|_{{\mS_0}}=\|\cH^{b}(u,u)-\cH^{b}(\wt u,\wt u)\|_{{\mS_0}}
\leq C_3\|b\|_{\mL^{q_b}_0(\mX_0')}\big(\|u\|_{{\mS_0}}+
\|\wt u\|_{{\mS_0}}\big)\|u-\wt u\|_{{\mS_0}}.
$$
Now let us choose $T\geq 1$ large enough and be fixed so that
$$
C_2 T^{-\delta_1}\|u_0\|_{\bB^{\beta_0}_{\bbb1;\bba}}\leq1/(12C_3\|b\|_{\mL^{q_b}_0(\mX_0')}).
$$
Suppose that for $\delta=\delta_0/\delta_1$ and some $C_0=C_0(\Theta)>0$,
$$
\|u_0\|_{\mZ_{\beta_0}}\leq\tfrac{T^{-\delta_0}}{12 C_1C_3\|b\|_{\mL^{q_b}_0(\mX_0')}}\leq
\tfrac{C_0}{\|b\|_{\mL^{q_b}_0(\mX_0')}((\|b\|_{\mL^{q_b}_0(\mX_0')}\|u_0\|_{\bB^{\beta_0}_{1;\bba}}\vee 1)^\delta},
$$
where we have used $T\geq 1$ in  the second inequality.
It is easy to see that for $\lambda=\tfrac{1}{4C_3\|b\|_{\mL^{q_b}_0(\mX_0')}}$,
$$
C_1T^{\delta_0}\|u_0\|_{\mZ_{\beta_0}}
+C_2 T^{-\delta_1}\|u_0\|_{\bB^{\beta_0}_{\bbb1;\bba}}
+C_3\|b\|_{\mL^{q_b}_0(\mX_0')}\lambda^2
\leq\lambda.
$$
In particular, 
\begin{align}\label{CS2}
\mbox{$\cU$ is a contraction mapping on }\cB^\lambda:=\big\{u\in {{\mS_0}}:\|u\|_{{\mS_0}}\leq\lambda\big\}.
\end{align}
On the other hand,  letting
$$
\bbp=\bbp'=\bbb1,\ \ \bbg'=\bbg=(-\beta_0,0),\ \ \theta=\sQ-1,\ \ \bar\bbg=\bbl_0,
$$ 
by \eqref{AS1} with the above choice of parameters, we also have
\begin{align}\label{AD402}
\|\cH^b(u,w)\|_{\mL^\infty_\bbg(\bB^{0,1}_{\bbb1;\bba})}\lesssim
\|b\|_{\mL^{q_b}_0(\mX_0')}\|u\|_{\mL^\infty_{\bbl_0}(\mX_0)}\|w \|_{\mL^\infty_\bbg(\bB^{0,1}_{\bbb1;\bba})}.
\end{align}
Moreover, by \eqref{AD0446} and \eqref{AD0456},
$$
\|P_\cdot u_0\|_{\mL^\infty_\bbg(\bB^{0,1}_{\bbb1;\bba})}\lesssim\|u_0\|_{\bB^{\beta_0}_{\bbb1;\bba}}.
$$
Hence,  there are two constants $C_4, C_5>0$
such that for all $u,\wt u\in\mS_0\cap \mL^\infty_\bbg(\bB^{0,1}_{\bbb1;\bba})$,
\begin{align}\label{CS5}
\|\cU(u)\|_{\mL^\infty_\bbg(\bB^{0,1}_{\bbb1;\bba})}
\leq C_4\|u_0\|_{\bB^{\beta_0}_{\bbb1;\bba}}
+C_5\|b\|_{\mL^{q_b}_0(\mX_0')}\|u\|_{{\mS_0}}\|u\|_{\mL^\infty_\bbg(\bB^{0,1}_{\bbb1;\bba})},
\end{align}
and 
\begin{align}\label{CS6}
\|\cU(u)-\cU(\wt u)\|_{\mL^\infty_\bbg(\bB^{0,1}_{\bbb1;\bba})}
\leq C_5\|b\|_{\mL^{q_b}_0(\mX_0')}
\big(\|u-\wt u\|_{{\mS_0}}\|u\|_{\mL^\infty_\bbg(\bB^{0,1}_{\bbb1;\bba})}+\|\wt u\|_{{\mS_0}}\|u-\wt u\|_{\mL^\infty_\bbg(\bB^{0,1}_{\bbb1;\bba})}\big).
\end{align}
Without loss of generality, we may assume $C_5\leq C_3$, where $C_3$ is the same as above.
Otherwise, we may replace above $C_3$ by bigger $C_5$.
\iffalse
Let $u^{(0)}:=u_0$ and for $n\in\mN$, we define recursively
$$
u^{(n)}:=\cU(u^{(n-1)}).
$$
By \eqref{CS2}, we have
$$
\|u^{(n)}\|_{\mS_0}\leq\lambda,\ \ \|u^{(n)}-u^{(m)}\|_{\mS_0}\leq \tfrac12\|u^{(n-1)}-u^{(m-1)}\|_{\mS_0}\leq\cdots\leq\tfrac\lambda {2^{n\wedge m}}.
$$
Moreover, by \eqref{CS5} and \eqref{CS6}, we also have
$$
\|u^{(n)}\|_{\mL^\infty_0(\mL^1)}\leq \|u_0\|_{\bB^{\beta_0}_{\bbb1;\bba}}
+\tfrac12\|u^{(n-1)}\|_{\mL^\infty_0(\mL^1)}
$$
and
$$
\|u^{(n)}-u^{(m)}\|_{\mL^\infty_0(\mL^1)}\leq \tfrac\lambda {2^{n\wedge m}}
C_3\|b\|_{\mL^{q_b}_0(\mX_0')}\|u^{(n-1)}\|_{\mL^\infty_0(\mL^1)}
+\tfrac12\|u^{(n-1)}-u^{(m-1)}\|_{\mL^\infty_0(\mL^1)}
$$
\fi
Thus by \eqref{CS2}, \eqref{CS5} and \eqref{CS6}, it is standard to derive that there is a unique $u\in\mS_0\cap\mL^\infty_\bbg(\bB^{0,1}_{\bbb1;\bba})$ such that
$$
u(t)=\cU(u)(t)=P_tu_0{\red-}\cH^{b}_t(u,u)
$$
and
\begin{align}\label{Cx1}
\|u\|_{\mL^\infty_\bbg(\bB^{0,1}_{\bbb1;\bba})}\lesssim \|u_0\|_{\bB^{\beta_0}_{\bbb1;\bba}}.
\end{align}
Moreover, the proofs of smoothness and stability are completely the same as Theorem \ref{Main}.

({\it Step 3})  In this step we prove the optimal large time decay estimate \eqref{aAA51}
by shifting time variable as used in Lemma \ref{Le38}. Let
$$
\eta_i:=\bba\cdot\big(\bbd-\tfrac{d}{\bbp_i}\big),\ i=0,1.
$$
Note that by \eqref{AB00} and \eqref{SH2},
$$
u(t+s)=P_su(t){\red-}\cH_s^{b^t}(u^t,u^t),\ \ s,t>0,
$$
where $u^t(s):=u(t+s)$. Thus for fixed $\beta\geq 0$ and any $t\geq 1$ and $i=0,1$,
\begin{align*}
\|u(2t)\|_{\bB^{\beta,1}_{\bbp_i;\bba}}\leq\|P_tu(t)\|_{\bB^{\beta,1}_{\bbp_i;\bba}}
+\|\cH_t^{b^t}(u^t,u^t)\|_{\bB^{\beta,1}_{\bbp_i;\bba}}.
\end{align*}
For the first term, by   \eqref{AD0456} and \eqref{Cx1} we have
$$
\|P_tu(t)\|_{\bB^{\beta,1}_{\bbp_i;\bba}}\lesssim t^{-\frac{\eta_i}\alpha}\|u(t)\|_{\bB^0_{\bbb1;\bba}}
\lesssim t^{-\frac{\eta_i}\alpha}\|u_0\|_{\bB^{\beta_0}_{\bbb1;\bba}},\ \ t\geq 1.
$$
For the second term,  by \eqref{AS1} with $\theta=\sQ-1$ and
$$
\bar\bbg=\bbg=\bbg_b=0,\ \ \bbg'=(-\theta,\gamma'_1),\ \gamma'_1\in(-\ell_1,1-\sQ),
$$
we derive that for any $t\geq 1$ and $i=0,1$,
\begin{align*}
\|\cH_\cdot^{b^t}(u^t,u^t)\|_{\mL^\infty_{\bbg'}(\bB^{\beta,1}_{\bbp_i;\bba})}
&\lesssim \|b^t\|_{\mL^{q_b}_0(\mX_0')}\|u^t\|_{\mL^\infty_0(\mX_\beta)}
\|u^t\|_{\mL^\infty_0(\bB^{\beta,1}_{\bbp_i;\bba})}\\
&\!\!\!\stackrel{\eqref{Sh1}}{\lesssim}\|b\|_{\mL^{q_b}_0(\mX'_0)}
t^{-\frac{\ell_1}\alpha}\|u\|_{\mL^\infty_{\bbl_\beta}(\mX_\beta)}\|u^t\|_{\mL^\infty_0(\bB^{\beta,1}_{\bbp_i;\bba})}.
\end{align*}
Combining the above two estimates, we obtain that for any $t\geq 1$,
$$
\|u(2t)\|_{\bB^{\beta,1}_{\bbp_i;\bba}}
\lesssim t^{-\frac{\eta_i}\alpha}\|u_0\|_{\bB^{\beta_0}_{\bbb1;\bba}}
+t^{-\frac{\gamma'_1+\ell_1}\alpha}
\|u^t\|_{\mL^\infty_0(\bB^{\beta,1}_{\bbp_i;\bba})},\ i=0,1.
$$
In particular, there are two constants $C_6=C_6(\beta,\Theta)>0$ and $C_7=C_7(\beta,\Theta,\|b\|_{\mL^{q_b}_0(\mX'_0)})>0$ 
such that for any $t\geq 2$,
$$
\|u(t)\|_{\bB^{\beta,1}_{\bbp_i;\bba}}
\leq C_6t^{-\frac{\eta_i}\alpha}\|u_0\|_{\bB^{\beta_0}_{\bbb1;\bba}}
+C_7 t^{-\frac{\gamma'_1+\ell_1}\alpha}
\|u^{t/2}\|_{\mL^\infty_0(\bB^{\beta,1}_{\bbp_i;\bba})},\ i=0,1.
$$
Since $\gamma'_1+\ell_1>0$, one can choose  $T_0=T_0(\beta,\Theta,\|b\|_{\mL^{q_b}_0(\mX'_0)},\gamma_1',\ell_1)\geq 2$ large enough so that 
$$
C_7 T_0^{-\frac{\gamma'_1+\ell_1}\alpha}\leq 2^{-1-\frac{\eta_i}\alpha},\ \ i=0,1,
$$
and so,
\begin{align*}
\sup_{t\geq T_0}\left(t^{\frac{\eta_i}\alpha}\|u(t)\|_{\bB^{\beta,1}_{\bbp_i;\bba}}\right)
&\leq C_6\|u_0\|_{\bB^{\beta_0}_{\bbb1;\bba}}
+2^{-1-\frac{\eta_i}\alpha}\sup_{t\geq T_0}\left(t^{\frac{\eta_i}\alpha}
\|u^{t/2}\|_{\mL^\infty_0(\bB^{\beta,1}_{\bbp_i;\bba})}\right)\\
& = C_6\|u_0\|_{\bB^{\beta_0}_{\bbb1;\bba}}
+2^{-1-\frac{\eta_i}\alpha}\sup_{t\geq T_0}\left(t^{\frac{\eta_i}\alpha}
\sup_{s\geq t/2}\|u(s)\|_{\bB^{\beta,1}_{\bbp_i;\bba}}\right)\\
& \leq C_6\|u_0\|_{\bB^{\beta_0}_{\bbb1;\bba}}
+\frac12\sup_{t\geq T_0}\sup_{s\geq t/2}\left( s^{\frac{\eta_i}\alpha}\|u(s)\|_{\bB^{\beta,1}_{\bbp_i;\bba}}\right)\\
& =C_6\|u_0\|_{\bB^{\beta_0}_{\bbb1;\bba}}
+\frac{T_0^{\frac{\eta_i}\alpha}}2\sup_{s\in[\frac{T_0}2, T_0]}\|u(s)\|_{\bB^{\beta,1}_{\bbp_i;\bba}} 
+\frac12\sup_{s\geq T_0}\left( s^{\frac{\eta_i}\alpha}\|u(s)\|_{\bB^{\beta,1}_{\bbp_i;\bba}}\right).
\end{align*}
Hence, 
\begin{align*}
\sup_{t\geq T_0}\left(t^{\frac{\eta_i}\alpha}\|u(t)\|_{\bB^{\beta,1}_{\bbp_i;\bba}}\right)
&\leq 2C_6\|u_0\|_{\bB^{\beta_0}_{\bbb1;\bba}}+T_0^{\frac{\eta_i}\alpha}\sup_{s\in[\frac{T_0}2, T_0]}\|u(s)\|_{\bB^{\beta,1}_{\bbp_i;\bba}}\leq C_8.
\end{align*}
The large time decay estimate now follows.
\end{proof}

\subsection{Local well-posedness of large initial value}
In this section we study the local well-posedness for FPE \eqref{FPE07} with large initial values. 
First of all we show the following asymptotic estimate of short time for $P_tf$ in Besov spaces.
\bl\label{Lem315H}
Let $\beta,\beta'\in\mR$ and $\bbp,\bbp'\in[1,\infty]^2$ with $\bbp\leq\bbp'$. 
Let $\sA:=\bba\cdot(\frac d{\bbp}-\frac{d}{\bbp'})$. Suppose that
$$
\beta<0,\ \ \sA+\beta'-\beta>0.
$$
For any $f\in\cup_{q\in[1,\infty)}\bB^{\beta,q}_{\bbp;\bba}$, it holds that
$$
\lim_{t\downarrow 0}\left(t^{\frac{\sA+\beta'-\beta}\alpha}\|P_tf\|_{\bB^{\beta',1}_{\bbp';\bba}}\right)=0
$$
\el
\begin{proof}
By the Bernstein inequality \eqref{Ber} and \eqref{Mm1}, we have 
\begin{align*}
\|\cR^\bba_jP_tf\|_{\bbp'}\lesssim 2^{j\sA}\|\cR^\bba_jP_tf\|_{\bbp}
\leq2^{j\sA}\sum_{\ell\in\Theta^t_j}\|\cR^\bba_j\Gamma_tp_t*\Gamma_t\cR^\bba_\ell f\|_\bbp.
\end{align*}
Here and below the implicit constant does not depend on $t\in(0,1)$.
Thus by Young's inequality,
\begin{align*}
t^{\frac{\sA+\beta'-\beta}\alpha}\sum_{j\geq 1}2^{j\beta'}\|\cR^\bba_jP_tf\|_{\bbp'}
\lesssim t^{\frac{\sA+\beta'-\beta}\alpha}\sum_{j\geq 1}\sum_{\ell\in\Theta^t_j}2^{j(\sA+\beta')}\|\cR^\bba_j\Gamma_tp_t\|_1\|\cR^\bba_\ell f\|_\bbp=\sum_{\ell\geq 0}c_{\ell,t}2^{\ell\beta}\|\cR^\bba_\ell f\|_\bbp,
\end{align*}
where
$$
c_{\ell,t}:=t^{\frac{\sA+\beta'-\beta}\alpha}\sum_{j\geq 1}2^{j(\sA+\beta')}2^{-\ell\beta}\b1_{\ell\in\Theta^t_j}\|\cR^\bba_j\Gamma_tp_t\|_1.
$$
By the definition of $\Theta^t_j$, Lemmas \ref{lem001} and \ref{LA1}, for $l\geq \frac{2(\sA+\beta'-\beta)}\alpha$, we have
\begin{align*}
c_{\ell,t}&=\sum_{j\geq 1}(t2^{j\alpha})^{\frac{\sA+\beta'}\alpha}
(t2^{\ell\alpha})^{-\frac{\beta}\alpha}\b1_{\ell\in\Theta^t_j}\|\cR^\bba_j\Gamma_tp_t\|_1\\
&\lesssim \sum_{j\geq 1}(t2^{j\alpha})^{\frac{\sA+\beta'}\alpha}
(t2^{j\alpha}(1+t2^{j\alpha}))^{-\frac{\beta}\alpha}((t 2^{j\alpha})^{-l}\wedge 1)\lesssim 1,
\end{align*}
and
\begin{align*}
\sum_{\ell\geq 1}c_{\ell,t}&=t^{\frac{\sA+\beta'-\beta}\alpha}\sum_{j\geq 1}2^{j(\sA+\beta')}\sum_{\ell\in\Theta^t_j}2^{-\ell\beta}\|\cR^\bba_j\Gamma_tp_t\|_1\\
&\!\!\!\stackrel{\eqref{AS88}}{\lesssim} t^{\frac{\sA+\beta'-\beta}\alpha}\sum_{j\geq 1}2^{j(\sA+\beta'-\beta)}(1+(t2^{j\alpha}))^{-\beta}((t 2^{j\alpha})^{-l}\wedge 1)\lesssim 1.
\end{align*}
Hence, for any $q\in[1,\infty)$, by H\"older's inequality,
$$
t^{\frac{\sA+\beta'-\beta}\alpha}\sum_{j\geq 1}2^{j\beta'}\|\cR^\bba_jP_tf\|_{\bbp'}
\lesssim\left[\sum_{\ell\geq 0}c_{\ell,t}(2^{\ell\beta}\|\cR^\bba_\ell f\|_\bbp)^q\right]^{\frac1q}
\left[\sum_{\ell\geq 1}c_{\ell,t}\right]^{1-\frac1q}
\lesssim\left[\sum_{\ell\geq 0}c_{\ell,t}(2^{\ell\beta}\|\cR^\bba_\ell f\|_\bbp)^q\right]^{\frac1q}.
$$
In particular, if $f\in\bB^{\beta_0,q}_{\bbp;\bba}$, then by the dominated convergence theorem,
\begin{align}\label{Mm2}
\lim_{t\downarrow 0}\left(t^{\frac{\sA+\beta'-\beta}\alpha}\sum_{j\geq 1}2^{j\beta'}\|\cR^\bba_jP_tf\|_{\bbp'}\right)
\lesssim\left[\sum_{\ell\geq 0}\lim_{t\downarrow 0}c_{\ell,t}(2^{\ell\beta}\|\cR^\bba_\ell f\|_\bbp)^q\right]^{\frac1q}=0.
\end{align}
On the other hand, by \eqref{AD0396}, for $t\in(0,1)$,
$$
\|\cR^\bba_0P_tf\|_{\bbp'}\lesssim \|f\|_{\bB^\beta_{\bbp;\bba}},
$$
which together with \eqref{Mm2} yields the desired limit.
\end{proof}
The following lemma is a direct application of the above lemma.
\bl
Under {\bf (H$_0$)}, for any $f\in \cup_{q\in[1,\infty)}(\bB^{\beta_0,q}_{\bbp_0;\bba}
\cap\bB^{\beta_0,q}_{\bbp_1;\bba})$, it holds that
\begin{align}\label{AD616}
\lim_{t\downarrow 0}\left(t^{\frac{\sA_0-\beta_0-\beta_b}\alpha}\|P_tf\|_{\mX_0}+t^{-\frac{\beta_0}\alpha}\|P_tf\|_{\mY_0}\right)=0.
\end{align}
\el

Now we can show the following local well-posedness for FPE \eqref{FPE07} with large initial value.
\bt\label{Cor35}
Suppose that {\bf (H$_{0}$)} and one of the following two conditions hold:

\begin{enumerate}[(i)]
\item $b\in\mL^{q_b}_0(\mX'_0)$ for some $q_b\in[1,\infty)$ and $u_0\in\mZ_{\beta_0}=\bB^{\beta_0,\infty}_{\bbp_0;\bba}\cap\bB^{\beta_0,\infty}_{\bbp_1;\bba}$;

\item $b\in\mL^\infty_0(\mX'_0)$ and $u_0\in\wt\mZ_{\beta_0}$, where 
$\wt\mZ_{\beta_0}:=\cup_{q\in[1,\infty)}(\bB^{\beta_0,q}_{\bbp_0;\bba}\cap\bB^{\beta_0,q}_{\bbp_1;\bba})\subset \mZ_{\beta_0}.$
\end{enumerate}
Then there is a small time $T_0>0$ depending only on $\Theta$ and the corresponding norms of $b$ and $u_0$ 
so that \eqref{FPE07} admits a unique smooth solution $u$ on $(0,T_0]$ and for any $\beta\geq 0$,
$$
\|u\b1_{[0,T_0]}\|_{\mS_\beta}<\infty,
$$
where $\mS_\beta$ is defined in \eqref{AS59}.
Moreover, for another $(\wt b,\wt u_0)$ satisfying the above (i) or (ii),
let $\wt u$ be the unique local solution corresponding to $(\wt b,\wt u_0)$ on the time interval $[0,\wt T_0]$.
Then for any $\beta\geq 0$,
\begin{align}\label{AS639}
\|(u-\wt u)\b1_{[0,T_0\wedge\wt T_0]}\|_{\mS_\beta}
\leq C_\beta \big(\|u_0-\wt u_0\|_{\mZ_{\beta_0}}+\|b-\wt b\|_{\mL^{q_b}_0(\mX_0')}\big),
\end{align}
where $C_\beta=C_\beta\big(T_0,\wt T_0,\Theta,\|b\|_{\mL^{q_b}_0(\mX'_0)},\|\wt b\|_{\mL^{q_b}_0(\mX'_0)}\big)>0$.
\et
\begin{proof}
(i) Since for $q_b\in[1,\infty)$,
$$
\lim_{T\downarrow 0}\|b\b1_{[0,T]}\|_{\mL^{q_b}_{\bbg_b}(\mX'_0)}=0,
$$
for any $u_0\in\mZ_{\beta_0}$, one can choose $T_0$ small enough so that
for the constant $C_0$ in Theorem \ref{Main},
$$
\|u_0\|_{\mZ_{\beta_0}}\leq C_0/\|b\b1_{[0,T_0]}\|_{\mL^{q_b}_{\bbg_b}(\mX'_0)}.
$$
Thus, we  conclude the proof by Theorem \ref{Main}.

(ii) Suppose $u_0\in \wt\mZ_{\beta_0}$. By the definition of $\mS_0$ and \eqref{AD616}, we have
$$
\lim_{T\downarrow 0}\|\b1_{(0,T]}(P_\cdot u_0)\|_{\mS_0}
=\lim_{T\downarrow 0}\sup_{t\in(0,T]}\left(t^{\frac{\sA_0-\beta_0-\beta_b}\alpha}\|P_tu_0\|_{\mX_0}+t^{-\frac{\beta_0}\alpha}\|P_tu_0\|_{\mY_0}\right)=0.
$$
From the proof of Theorem \ref{Main}, one can choose $T_0$ small enough so that
$\cU$ defined by \eqref{AD301} is a contraction  mapping when we limit on the time interval $[0,T_0]$.
Thus we can repeat the proof of Theorem \ref{Main} to conclude the proof.
\end{proof}

\br\rm\label{Re318}
If the condition \eqref{AD066} is replaced with
$$
0<\sA_0-\beta_0-\beta_b<\sQ-1\ \mbox{ and }\ \sA_0-2\beta_0-\beta_b<\sQ,
$$
which is equivalent to 
\begin{align}\label{AD0696} 
0<\sA_0-\beta_0-\beta_b<\sQ+\beta_0\wedge(-1),
\end{align}
then the cases (i) and (ii) in Theorem \ref{Cor35} can be combined with only one case that
 $b\in\mL^{q_b}_0(\mX'_0)$ for some $q_b\in[1,\infty]$ and $u_0\in\mZ_{\beta_0}$ since
 $\mZ_{\beta_0}\subset \cap_{\eps>0}\wt\mZ_{\beta_0+\eps}$. We would like to point out that the above (i) and (ii) 
 allow us to consider some critical cases.
\er

\subsection{Global well-posedness of large initial value}
In this subsection we shall consider two cases that the smallness condition \eqref{Small} in Theorem \ref{Main} can be dropped.

Our first result of global well-posedness  with large initial value is
\bt\label{TH1}
Let $\alpha\in{(1,2]}$, $q_b\in(\frac{\alpha}{\alpha-1},\infty]$
and $\bbrho_0,\bbrho_1\in[1,\infty]^2$ with $\bbrho_0\leq\bbrho_1$. 
Let $\beta_b\in\mR$ and $\beta_0\in(\tfrac\alpha{q_b}-\alpha,0)$. Suppose that $\beta_b\not=\bba\cdot\tfrac d{\bbrho_i}-\beta_0$, $i=0,1$ and
\begin{align}\label{AD0696} 
0<\bba\cdot\tfrac d{\bbrho_0}-\beta_0-\beta_b
\leq\alpha-\tfrac\alpha{q_b}-1,\ \ \bba\cdot\tfrac d{\bbrho_0}-2\beta_0-\beta_b
<\alpha-\tfrac\alpha{q_b},
\end{align}
and one of the following two conditions holds:

\begin{enumerate}[(i)]
\item $b\in\mL^{q_b}_0(\mX'_0)$ for some $q_b\in[1,\infty)$ and $u_0\in\bB^{\beta_0,\infty}_{1;\bba}$;

\item $b\in\mL^\infty_0(\mX'_0)$ and $u_0\in\cup_{q\in[1,\infty)}\bB^{\beta_0,q}_{1;\bba}\subset \bB^{\beta_0,\infty}_{1;\bba}.$
\end{enumerate}
Then there is a unique smooth solution $u$ to \eqref{FPE07} on $\mR_+$ so that for any $T,\beta>0$,
$$
\sup_{t\in(0,T]}\left(t^{(\beta-\beta_0)/\alpha}\|u(t)\|_{\bB^{\beta,1}_{\bbb1;\bba}}\right)<\infty.
$$
Moreover, for another $(\wt b,\wt u_0)$ satisfying the above (i) or (ii),
let $\wt u$ be the unique solution corresponding to $(\wt b,\wt u_0)$. For any $T,\beta>0$, there is a constant 
$C_\beta=C_\beta(T,\Theta,\|b\|_{\mL^{q_b}_0(\mX_0')},\|\wt b\|_{\mL^{q_b}_0(\mX_0')})>0$ such that
$$
\sup_{t\in(0,T]}\left(t^{(\beta-\beta_0)/\alpha}\|u(t)-\wt u(t)\|_{\bB^{\beta,1}_{\bbb1;\bba}}\right)
\leq C_\beta \big(\|u_0-\wt u_0\|_{\bB^{\beta_0}_{\bbb1;\bba}}+\|b-\wt b\|_{\mL^{q_b}_0(\mX_0')}\big).
$$
\et
\begin{proof}
If we take $\bbp_0=1$ in {\bf (H$_0$)}, then \eqref{AD066} reduces to \eqref{AD0696}. By Theorem \ref{Cor35},
let $T_0>0$ be the short time of the existence of  a smooth solution.  For any fixed $t_0\in(0,T_0]$, we have
$$
u\in \cap_{\beta>0}C([t_0,T_0];\bB^{\beta}_{\bbb1;\bba})
$$ 
and
$$
\p_t u=\Delta^{\alpha/2}_vu-v\cdot\nabla_x u-\div_v((b*u)u),\ t\in[t_0,T_0].
$$
Now, starting from the time $t_0$, since $\beta_0<0$ and $\|u(t_0)\|_{\bB^{\beta_0,q}_{\bbb1;\bba}}\leq C\|u(t_0)\|_{\mL^1}$ for any $q\in[1,\infty]$,
by Theorem \ref{Cor35} again, there is a time length $\delta>0$ only depending on $\Theta$, $\|b\|_{\mL^{q_b}_0(\mX_0')}$
and the bound of $\|u(t_0)\|_{\mL^1}$ so that
there is a unique smooth solution on the time interval $[t_0,t_0+\delta]$. 
On the other hand, by Lemma \ref{LEC1} in appendix we have
$$
\|u(t_0+\delta/2)\|_{\mL^1}\leq\|u(t_0)\|_{\mL^1}.
$$
Next starting from the time $t_0+\delta/2$, we can proceed to find a unique smooth 
solution to \eqref{FPE07}  on the time interval $[t_0+\delta/2,t_0+3\delta/2]$ and so on. Thus we obtain a global solution.
The regularity and stability estimates are direct consequences of Theorem \ref{Cor35}.
\end{proof}
%\br The above result covers \cite{RZ21} and \cite{CJM22}.\er
Our second global well-posedness result of large initial value is
\bt\label{TH2}
Suppose that {\bf (H$_{0}$)} holds with $\bbp_0=(p_0,p_0)$ and $\bbrho_1=(\varrho_1,\varrho_1)$,  
 and one of the following two conditions holds:

\begin{enumerate}[(i)]
\item $b\in\mL^{q_b}_0(\mX'_0)$ for some $q_b\in[1,\infty)$, $\div_v b=0$ and $u_0\in\mZ_{\beta_0}=\bB^{\beta_0,\infty}_{\bbp_0;\bba}\cap\bB^{\beta_0,\infty}_{\bbp_1;\bba}$;

\item $b\in\mL^\infty_0(\mX'_0)$, $\div_v b=0$ and $u_0\in\wt\mZ_{\beta_0}$, where 
$\wt\mZ_{\beta_0}:=\cup_{q\in[1,\infty)}(\bB^{\beta_0,q}_{\bbp_0;\bba}\cap\bB^{\beta_0,q}_{\bbp_1;\bba})\subset \mZ_{\beta_0}.$
\end{enumerate}
Then there is a unique smooth solution $u$ to \eqref{FPE07} on $\mR_+$ so that for any $T,\beta>0$,
$$
\|u\b1_{[0,T]}\|_{\mS_\beta}<\infty,
$$
where $\mS_\beta$ is defined by \eqref{AS59}.
Moreover, for another $(\wt b,\wt u_0)$ satisfying the above (i) or (ii),
let $\wt u$ be the unique solution  on $\mR_+$ corresponding to $(\wt b,\wt u_0)$. 
Then for any $T,\beta>0$, there is a constant $C_\beta=C_\beta(T,\Theta,\|b\|_{\mL^{q_b}_0(\mX_0')},\|\wt b\|_{\mL^{q_b}_0(\mX_0')})>0$ such that
$$
\|(u-\wt u)\b1_{[0,T]}\|_{\mS_\beta}
\leq C_\beta \big(\|u_0-\wt u_0\|_{\mZ_{\beta_0}}+\|b-\wt b\|_{\mL^{q_b}_0(\mX_0')}\big).
$$
\et
\begin{proof}
Let $T_0>0$ be the short time of the existence of a smooth solution in Theorem \ref{Cor35}.  
For any fixed $t_0\in(0,T_0)$, we have
$$
u\in C([t_0,T_0];C^\infty_b(\mR^{2d})),\ \ \DD:=b*u\in L^1([t_0,T_0]; C^1_b(\mR^{2d})),
$$ 
and
$$
\p_t u=\Delta^{\alpha/2}_vu-v\cdot\nabla_x u-\div_v(\DD u),\ t\in[t_0,T_0].
$$
Since $\div_vb =0$, we have $\div_v \DD=0$. Thus by Lemma \ref{LEC1} in appendix, for $t\in(t_0,T_0)$,
\begin{align}\label{AB0}
\|u(t)\|_{\mL^1}\leq\|u(t_0)\|_{\mL^1},\ \ \|u(t)\|_{\mL^\infty}\leq\|u(t_0)\|_{\mL^\infty}.
\end{align}
Noting that for $p\in[1,\infty]$,
$$
\|f\|_{\mL^p}\leq\|f\|_{\mL^\infty}^{1-\frac1{p}}\|f\|_{\mL^1}^{\frac1{p}},
$$
by $\bbp_0=(p_0,p_0)$ and $\bbrho_1=(\varrho_1,\varrho_1)$, recalling $p_1=p_0\wedge\frac{\varrho_1}{\varrho_1-1}$, we have
for any $q\in[1,\infty]$,
\begin{align*}
\sum_{i=0,1}\|u(t_0)\|_{\bB^{\beta_0,q}_{\bbp_i;\bba}}
\stackrel{\eqref{AB2}}{\lesssim}\sum_{i=0,1}\|u(t_0)\|_{\mL^{p_i}}
\leq\sum_{i=0,1}\|u(t_0)\|_{\mL^\infty}^{1-\frac1{p_i}}\|u(t_0)\|_{\mL^1}^{\frac1{p_i}}
=:\flat(t_0).
\end{align*}
Now, starting from $t_0$, by Theorem \ref{Cor35} again, there is a time length $\delta>0$ depending only on $\Theta$, 
$\|b\|_{\mL^{q_b}_0(\mX_0')}$ and the bound of $\flat(t_0)$  so that
there is a unique solution on the time interval $[t_0,t_0+\delta]$. 
On the other hand, by \eqref{AB0} we have
$$
\flat(t_0+\delta/2)\leq\flat(t_0).
$$
Next starting from the time $t_0+\delta/2$, we can proceed to find a unique smooth solution to \eqref{FPE07}  on the time interval $[t_0+\delta/2,t_0+3\delta/2]$ and so on. Thus we obtain a global solution. The regularity and stability estimates are direct consequences of Theorem \ref{Cor35}.
\end{proof}
\br\rm
If $u$ is probability measure, then the condition $\div_v b=0$ can be replaced by $(\div_v b)^-=0$
since $(\div_v\DD)^-=(\div_v b)^-*u=0$, and by Lemma \ref{LEC1} we still have \eqref{AB0}.
\er

As a special case, we also consider the non-degenerate fractional Fokker-Planck equation:
\begin{align}\label{FPE01}
\p_t u=\Delta^{\alpha/2}u-\div(( b*u)u),\ \ u(0)=u_0.
\end{align}
By Duhamel's formula, we have
$$
u(t)=P_t u_0-\int^t_0 P_{t-s}\div((b*u)u)\dif s=P_tu_0-\cH^b_t(u,u),
$$
where $P_t$ is the semigroup associated with $\Delta^{\alpha/2}$.
We assume
\begin{enumerate}[{\bf ($\wt {\bf H}_{0}$)}]
\item Let $\alpha\in{(1,2]}$, $q_b\in(\frac{\alpha}{\alpha-1},\infty]$
and $p_0,\varrho_0,\varrho_1\in[1,\infty]$ with $\varrho_0\leq \varrho_1$ and $1\leq\tfrac1{p_0}+\tfrac1{\varrho_0}$. 
Let $\bar\varrho_0,\bar\varrho_1\in[1,\infty]$ be defined by $\tfrac1{\bar\varrho_0}+\tfrac1{\varrho_0}=1$ and
$\tfrac1{\bar\varrho_1}+\tfrac1{\varrho_1}=1.$
Let $p_1:=p_0\wedge\bar{\varrho}_1$ and
$$
\sQ:=\alpha-\tfrac\alpha{q_b},\ \ \sA_i:=\tfrac d{p_i}+\tfrac d{\varrho_i}-d=\tfrac d{p_i}-\tfrac d{\bar\varrho_i},\ i=0,1.
$$
Let $\beta_0\in(-\sQ,0)$ and $\mR\ni\beta_b\not=\sA_i-\beta_0, i=0,1$. Suppose that 
\begin{align}\label{AD606} 
0<\sA_0-\beta_0-\beta_b\leq\sQ-1\ \mbox{ and }\ \sA_0-2\beta_0-\beta_b<\sQ.\
\end{align}
\end{enumerate}
As in \eqref{MY}, we introduce three isotropic Besov spaces: for $\beta\geq 0$, 
\begin{align}\label{MY0}
\bar\mX_\beta:=\bigcap_{i=0,1}\bB^{\beta-\beta_b,1}_{\bar\varrho_i},\ \
\bar\mY_\beta:=\bigcap_{i=0,1}\bB^{\beta,1}_{p_i},\ \ 
\bar\mZ_\beta:=\bigcap_{i=0,1}\bB^{\beta}_{p_i}.
\end{align}
The following result is a combination of Theorems  \ref{Main}, \ref{thm230}, \ref{TH1} and \ref{TH2}.
Since it can be considered as a special case and the proofs are completely the same, we omit the details.
\bt\label{thm02}
Suppose {\bf ($\wt {\bf H}_{0}$)} and $b\in\mL^{q_b}_0(\bar\mX_0')$, where $\bar\mX_0'=\bB^{\beta_0}_{\varrho_0}+\bB^{\beta_0}_{\varrho_1}$ is the dual space of $\bar\mX_0$. 
For any $\gamma_b>(\sQ-1-\sA_1)\vee 0$, there is a constant $C_0=C_0(\Theta)>0$ 
such that for any fixed $T\geq 1$ and $u_0\in\bar\mZ_{\beta_0}$ satisyfying
\begin{align}\label{Mb11}
\|u_0\|_{\bar\mZ_{\beta_0}}\|b\|_{\mL^{q_b}_0(\bar\mX_0')}\leq C_0T^{-\gamma_b/\alpha},
\end{align}
there is a unique weak solution  $u$ to FPE \eqref{FPE01} on $[0,T]$ such that
$$
\sup_{t\in(0,T]}\left(t^{-\frac{\beta_0}\alpha}\|u(t)\|_{\bar\mY_0}+t^{\frac{\sA_0-\beta_0-\beta_b}\alpha}\|u(t)\|_{\bar\mX_0}\right)<\infty.
$$
Moreover, we have the following conclusions:
\begin{enumerate}[(i)]
\item{\bf (Smoothness)} For any $\beta\geq 0$, it holds that
$$
\sup_{t\in(0,T]}\left(t^{\frac{\beta-\beta_0}\alpha}\|u(t)\|_{\bar\mY_\beta}+t^{\frac{\beta+\sA_0-\beta_0-\beta_b}\alpha}\|u(t)\|_{\bar\mX_\beta}\right)
<\infty.
$$
\item {\bf (Stability)} For another $\wt b\in\mL^{q_b}_0(\bar\mX_0')$ and $\wt u_0\in\bar\mZ_{\beta_0}$ satisfying the same condition \eqref{Mb11},
let $\wt u$ be the unique solution on $[0,T]$ corresponding to $\wt b$ and $\wt u_0$.  For any $\beta\geq 0$,
there is a constant $C_\beta=C_\beta(T,\Theta,\|b\|_{\bar\mL^{q_b}_0(\bar\mX'_0)},\|\wt b\|_{\mL^{q_b}_0(\bar\mX'_0)})>0$ such that
$$
\sup_{t\in(0,T]}\left(t^{\frac{\beta-\beta_0}\alpha}\|u(t)-\wt u(t)\|_{\bar\mY_\beta}\right)
\leq C_\beta \big(\|u_0-\wt u_0\|_{\bar\mZ_{\beta_0}}+\|b-\wt b\|_{\mL^{q_b}_{0}(\bar\mX_0')}\big).
$$
\item {\bf (Large time decay estimate)} If $\frac d{\varrho_1}>\sQ-1$, $p_0\not=1$ and $u_0\in\bar\mZ_{\beta_0}\cap \bB^{\beta_0}_1$, 
then we can take $\gamma_b=0$ in \eqref{Mb11}, and the constant $C_0$ may depend on $\|u_0\|_{\bB^{\beta_0}_1}$, and there is a global solution on $\mR_+$ so that
for any $\beta\geq 0$ and $i=0,1$, 
$$
\sup_{t\ge1}\left(t^{\frac{d-d/p_i}\alpha}\|u(t)\|_{\bB^{\beta,1}_{p_i}}\right)<\infty.
$$

\item {\bf (Global solution)} In addition, we 
suppose that one of the following four conditions holds:
\begin{enumerate}[{\rm (a)}]
\item $q_b\in[1,\infty)$ and $u_0\in \bB^{\beta_0}_1$.
\item $q_b\in[1,\infty)$, $\div b=0$ and $u_0\in \bB^{\beta_0}_{p_0}\cap\bB^{\beta_0}_{p_1}$.
 \item $q_b=\infty$ and $u_0\in \cup_{q\in[1,\infty)}\bB^{\beta_0,q}_1$.
 \item $q_b=\infty$, $\div b=0$ and $u_0\in \cup_{q\in[1,\infty)}(\bB^{\beta_0,q}_{p_0}\cap\bB^{\beta_0,q}_{p_1})$.
\end{enumerate}
Then we can drop the smallness assumption \eqref{Mb11} and have a global solution.
\end{enumerate}
\et

\section{Well-posedness of mean-field kinetic SDEs}

In this section we use previous results to show the weak and strong well-posedness for mean-field kinetic SDEs.
Although the associated nonlinear Fokker-Planck equation has a smooth solution for any positive time, we have to carefully deal with 
the time singularity at initial time.
\subsection{Well-posedness of linear kinetic SDEs}
In this subsection, we first study the weak/strong well-posedness to the following linear kinetic SDE:
for $\alpha\in(1,2]$,
\begin{align}\label{KSDE}
\left\{
\begin{aligned}
&X_t=X_0+\int_0^t V_s\dif s,\\
&V_t=V_0+\int^t_0\DD_s(X_s,V_s)\dif s+L^{(\alpha)}_t,
\end{aligned}
\right.
\end{align}
where $\DD: \mR_+\times\mR^{2d}\to\mR^d$ is a measurable vector field, and  $(L^{(\alpha)}_t)_{t\geq 0}$ is the  $\alpha$-stable process as in the introduction.
Moreover, if we let $Z_t:=(X_t,V_t)$, $\bD_t(z)=\bD_t(x,v)=(v,\DD_t(x,v))^*\in\mR^{2d}$ and $\sigma:=(0,\mI)^*\in\mR^{2d}\otimes\mR^d$, then SDE \eqref{KSDE} can be written as 
\begin{align}\label{MV13}
Z_t=Z_0+\int^t_0\bD_s(Z_s)\dif s+\sigma L^{(\alpha)}_t,\ Z_0\sim\mu_0\in\cP(\mR^{2d}).
\end{align}
For  $q\in[1,\infty]$, $T>0$ and a Banach space $\mB$, we write
$$
\mL^q_T(\mB):=L^q([0,T];\mB),\ \mL^\infty_T:=L^\infty([0,T]\times\mR^{2d}).
$$

We have the following well-posedness result for SDE \eqref{MV13}.
\bt\label{thmW}
Let $\alpha\in(1,2]$, $q\in(\frac{\alpha}{\alpha-1},\infty]$ and $T,\beta>0$.  
\begin{enumerate}[(i)]
\item {\bf (Weak well-posedness)} Suppose $\DD\in \mL^q_T(\bC^\beta_\bba)$. For any initial distribution $\mu_0$,
there is a unique weak solution $Z$ to kinetic SDE \eqref{MV13}.
\item {\bf (Strong well-posedness)} Suppose $\DD\in \mL^q_T(\bC^{\frac{\alpha+\beta}{1+\alpha}}_x\cap\bC^\beta_\bba)$, where $\beta>1-\tfrac\alpha 2$. For any initial value $Z_0$,
there is a unique strong solution $Z$ to kinetic SDE \eqref{MV13}.
\item {\bf (Strong stability)} Let $\wt\DD\in \mL^q_T(\bC^{\frac{\alpha+\beta}{1+\alpha}}_x\cap\bC^\beta_\bba)$ be another vector field with the same $\beta>1-\tfrac\alpha 2$ and $\wt Z_0$ another initial random variable
on the same probability space. Let $\wt Z$ be the strong solution of \eqref{MV13} corresponding to $\wt H$ and $\wt Z_0$. Then we have
\begin{align}\label{SW3}
\mE\left(\sup_{t\in[0,T]}|Z_t-\wt Z_t|^2\right)\lesssim_C \mE|Z_0-\wt Z_0|^2
+\|H-\wt H\|^2_{\mL^q_T(\bC^{\beta}_{\bba})}.
\end{align}
\end{enumerate}
\et

To show this result, we need to study the following linear kinetic equation with transport drift:
\begin{align}\label{LKE}
\p_t u=\Delta_v^{\alpha/2} u+v\cdot\nabla_x u-\lambda u+\DD\cdot\nabla_v u+f,\quad u_0=0,
\end{align}
where $\lambda>0$, $\DD: \mR_+\times\mR^{2d}\to\mR^d$ and $f: \mR_+\times\mR^{2d}\to\mR$ are two Borel measurable functions.

We first show the following simple result.
\bt\label{thm41}
Let $T>0$, $\alpha\in(1,2]$, $\beta,\gamma>0$ and  $q\in(\frac{\alpha}{\alpha-1},\infty]$. 
\begin{enumerate}[(i)]
\item For any $\DD,f\in\mL^q_T(\bC^{\beta}_\bba)$,
there is a unique weak solution $u$ to equation \eqref{LKE} with the regularity that for any $q'\in[q,\infty]$, $\beta'\in[0,\alpha+\frac\alpha{q'}-\frac\alpha q)$
and all $\lambda\geq 1$,
\begin{align}\label{FA00}
\|u\|_{\mL^{q'}_{T}(\bC^{\beta'+\beta}_\bba)}\lesssim_C \lambda^{\frac1q+\frac{\beta'}\alpha-1-\frac1{q'}}\|f\|_{\mL^q_T(\bC^\beta_\bba)},
\end{align}
where $C>0$ only depends on $T,d,\beta,\alpha,q,q',\beta'$ and $\|\DD\|_{\mL^q_T(\bC^\beta_\bba)}$.
\item Let $\DD,f,\wt\DD,\wt f\in\mL^q_T(\bC^{\beta}_{\bba})$ and $u,\wt u$ be the solutions of 
equation \eqref{LKE} corresponding to $(\DD, f)$ and $(\wt\DD,\wt f)$, respectively.
For any $q'\in[q,\infty]$ and $\beta'\in[0,\alpha+\frac\alpha{q'}-\frac\alpha q)$, 
there is a constant $C>0$ depending on $T,d,\beta,\alpha,q,q',\beta'$ and $ \|\DD\|_{\mL^q_T(\bC^{\beta}_{\bba})}$,
$ \|\wt\DD\|_{\mL^q_T(\bC^{\beta}_{\bba})}$  such that for all $\lambda\geq 1$,
\begin{align}\label{FA101}
\|u-\wt u\|_{\mL^{q'}_{T}(\bC^{\beta'+\beta}_{\bba})}\lesssim_C 
\lambda^{\frac1q+\frac{\beta'}\alpha-1-\frac1{q'}}\Big(\|H-\wt H\|_{\mL^q_T(\bC^{\beta}_{\bba})}
+\|f-\wt f\|_{\mL^q_T(\bC^{\beta}_{\bba})}\Big).
\end{align}  
\item Let $\DD,f\in\mL^q_T(\bC^{\gamma,\beta}_{x,\bba})$ and $u$ be the solution of \eqref{LKE}.
For any $q'\in[q,\infty]$ and $\beta'\in[0,\alpha+\frac\alpha{q'}-\frac\alpha q)$,  
there is a constant $C>0$ depending on $T,d,\gamma,\beta,\alpha,q,q',\beta'$ and $ \|\DD\|_{\mL^q_T(\bC^{\gamma,\beta}_{x,\bba})}$ such that for all $\lambda\geq 1$,
\begin{align}\label{FA001}
\|u\|_{\mL^{q'}_{T}(\bC^{\gamma,\beta'+\beta}_{x,\bba})}\lesssim_C \lambda^{\frac1q+\frac{\beta'}\alpha-1-\frac1{q'}}\|f\|_{\mL^q_T(\bC^{\gamma,\beta}_{x,\bba})},
\end{align} 
\end{enumerate}
\et
\begin{proof}
(i) We only prove the a priori estimate \eqref{FA00}. By Duhamel's formula, we have
\begin{align*}
u_t&=\int_0^tP^\lambda_{t-s}(\DD_s\cdot\nabla_v u_s+f_s)\dif s,
\end{align*}
where
$$
P^\lambda_tf(z)=\e^{-\lambda t} P_tf(z),\ \ z\in\mR^{2d},\ t>0.
$$
By \eqref{AD0446} with $\bbp_0=\bbp'=\infty$, for any $\beta'\in[0,\alpha)$, we have
\begin{align}
\|u_t\|_{\bC^{\beta'+\beta}_\bba}&\lesssim\int_0^t\e^{-\lambda(t-s)}(t-s)^{-\frac{\beta'}\alpha}
\Big(\|\DD_s\cdot\nabla_v u_s\|_{\bC^{\beta}_\bba}+\|f_s\|_{\bC^{\beta}_\bba}\Big)\dif s\no\\
&\lesssim\int_0^t\e^{-\lambda(t-s)}(t-s)^{-\frac{\beta'}\alpha}\Big(\|\DD_s\|_{\bC^{\beta}_\bba}\|u_s\|_{\bC^{1+\beta}_\bba}+\|f_s\|_{\bC^{\beta}_\bba}\Big)\dif s.\label{Aw1}
\end{align}
In particular, taking $\beta'=1$, by H\"older's inequality, we have for $r=\frac{q}{q-1}$,
\begin{align*}
\|u_t\|^r_{\bC^{1+\beta}_\bba}&\lesssim\|\DD\|^r_{\mL^q_T(\bC^{\beta}_\bba)}\int_0^t\e^{-r\lambda(t-s)}(t-s)^{-\frac{r}\alpha}\|u_s\|^r_{\bC^{1+\beta}_\bba}\dif s
+\|f\|^r_{\mL^q_T(\bC^{\beta}_\bba)}\int_0^t\e^{-r\lambda s}s^{-\frac{r}\alpha}\dif s.
\end{align*}
Since $\frac 1\alpha<\frac 1r=1-\frac1q$, by Gronwall's inequality, we get
$$
\|u\|_{\mL^\infty_T(\bC^{1+\beta}_\bba)}\lesssim \|f\|_{\mL^q_T(\bC^{\beta}_\bba)} \lambda^{\frac1\alpha-\frac1r},\ \lambda\geq 1.
$$
Let $1+\frac1{q'}=\frac1r+\frac1q$. By \eqref{Aw1} and Young's inequality we have for $\lambda\geq 1$,
\begin{align*}
\|u\|_{\mL^{q'}_T(\bC^{\beta'+\beta}_\bba)}
&\lesssim\left(\|u\|_{\mL^\infty_T(\bC^{1+\beta}_\bba)}\|H\|_{\mL^q_T(\bC^\beta_\bba)}
+\|f\|_{\mL^q_T(\bC^\beta_\bba)}\right)
\left(\int_0^T\e^{-r\lambda s}s^{-\frac{\beta' r}\alpha}\dif s\right)^{1/r}\\
&\lesssim\|f\|_{\mL^q_T(\bC^\beta_\bba)} \lambda^{\frac{\beta'}\alpha-\frac1r}
\left(\int_0^\infty\e^{-r s}s^{-\frac{\beta' r}\alpha}\dif s\right)^{1/r}
\lesssim\|f\|_{\mL^q_T(\bC^\beta_\bba)} \lambda^{\frac1q+\frac{\beta'}\alpha-1-\frac1{q'}}.
\end{align*}
Thus we get \eqref{FA00}.

(ii) Let $U=u-\wt u$. Then $U$ solves the following PDE:
$$
\p_t U=\Delta_v^{\alpha/2} U+v\cdot\nabla_x U-\lambda U+\DD\cdot\nabla_v U+G,\quad U_0=0,
$$
where
$$
G:=f-\wt f+(\DD-\wt \DD)\cdot\nabla_v \wt u.
$$
By \eqref{FA00}, we have
\begin{align}\label{FA003}
\|U\|_{\mL^{q'}_{T}(\bC^{\beta'+\beta}_{\bba})}
\lesssim \lambda^{\frac1q+\frac{\beta'}\alpha-1-\frac1{q'}}\|G\|_{\mL^q_T(\bC^{\beta}_{\bba})}.
\end{align}
Note that
\begin{align*}
\|G\|_{\mL^q_T(\bC^{\beta}_{\bba})}
&\leq \|f-\wt f\|_{\mL^q_T(\bC^{\beta}_{\bba})}+\|(H-\wt H)\cdot\nabla_v\wt u\|_{\mL^q_T(\bC^{\beta}_{\bba})}\\
&\leq \|f-\wt f\|_{\mL^q_T(\bC^{\beta}_\bba)}+\|H-\wt H\|_{\mL^q_T(\bC^{\beta}_{\bba})}
\|\wt u\|_{\mL^\infty_T(\bC^{1+\beta}_{\bba})}.
\end{align*}
Substituting it into \eqref{FA003} and by \eqref{FA00} with $q'=\infty$ and $\beta'=1$, we obtain \eqref{FA101}.

(iii) Noting that for any $j\geq 0$,
$$
\|\cR^x_jP_tf\|_{\bC^{\beta'+\beta}_{\bba}}\lesssim (1\wedge t)^{-\frac{\beta'}\alpha}\|\cR^x_jf\|_{\bC^{\beta}_{\bba}},
$$
we have for any $\gamma\in\mR$,
$$
\|P_tf\|_{\bC^{\gamma,\beta'+\beta}_{x,\bba}}\lesssim (1\wedge t)^{-\frac{\beta'}\alpha}\|f\|_{\bC^{\gamma,\beta}_{x,\bba}}.
$$
By this estimate, the estimate \eqref{FA001} is completely the same as \eqref{FA00}.
\end{proof}

For $\alpha\in(1,2)$, let $N(\dif w,\dif t)$ be the Poisson random measure associated with $\alpha$-stable process $L^{(\alpha)}_t$.
More precisely, for any Borel set $E\in\sB(\mR^d)$ and $t>0$,
$$
N(E\times [0,t]):=\sum_{s\in(0,t]}\b1_{E}(L^{(\alpha)}_s-L^{(\alpha)}_{s-}).
$$
Let $\wt N(\dif w, [0,t]):=N(\dif w, [0,t])-\nu^{(\alpha)}(\dif w) t$ be the compensated Poisson random measure, 
where $\nu^{(\alpha)}(\dif w)=c_{d,\alpha}\dif w/|w|^{d+\alpha}$ 
is the L\'evy measure of $\alpha$-stable process $L^{(\alpha)}_t$, and $c_{d,\alpha}$ is a normalized constant only depending on $d,\alpha$. 
In particular,
$$
L^{(\alpha)}_t=\int_{\mR^d} w \wt N(\dif w,[0,t]).
$$
The following lemma is crucial for weak uniqueness and pathwise uniqueness.
 \bl\label{lem44}
 Let $T>0$, $\alpha\in(1,2]$, $\beta>0$ and  $q\in(\frac{\alpha}{\alpha-1},\infty]$. Suppose that $\DD,f\in\mL^q_T(\bC^{\beta}_\bba)$.
Let $u$ be the unique weak solution of the following backward kinetic PDE (see Theorem \ref{thm41}):
\begin{align}\label{Aw2}
\p_t u+\Delta_v^{\alpha/2} u+v\cdot\nabla_x u+\DD\cdot\nabla_v u=f,\quad u(T)=0,
\end{align}
and $Z$ be a solution of kinetic SDE \eqref{MV13}.   If $\alpha\in(1,2)$, then for any $t\in[0,T]$,
\begin{align}\label{SDE34}
u(t,Z_t)=u(0,Z_0)+\int^t_0 f(s, Z_s)\dif s+\int^t_0\!\!\int_{\mR^d}\left(u(s,Z_{s-}+\sigma w)-u(s,Z_{s-})\right)\wt{N}(\dif w, \dif s);
\end{align}
if $\alpha=2$, then for any $t\in[0,T]$,
\begin{align}\label{SDE304}
u(t,Z_t)=u(0,Z_0)+\int^t_0 f(s, Z_s)\dif s
+\int^t_0\nabla_v u(s,Z_s)\dif W_s.
\end{align}
\el
\begin{proof}
We only prove \eqref{SDE34} since \eqref{SDE304} is similar and easier.
Let $\varGamma\in C^\infty_c(\mR^{2d})$ be a smooth density function with support in $B^\bba_1$. 
For each $n\in\mN$, define 
\begin{align}\label{Mo5}
\varGamma_{n}(z)=\varGamma_{n}(x,v):=n^{(2+\alpha)d}\varGamma(nx,n^{1+\alpha}v)
\end{align}
and
$$
u^n(t,z):=(u(t,\cdot)*\varGamma_{n})(z).
$$
Taking convolutions for both sides of \eqref{Aw2} with $\varGamma_n$, we get
\begin{align*}
\p_t u_n+(\Delta^{\alpha/2}_v+v\cdot\nabla_x)u_n+\DD\cdot\nabla_v u_n=f_n+g_n+h_n,
\end{align*}
where
\begin{align*}
g_n=v\cdot\nabla_x u_n-(v\cdot\nabla_x u)*\varGamma_n,\ \ h_n:=\DD\cdot\nabla_v u_n-(\DD\cdot\nabla_v u)*\varGamma_n.
\end{align*}
By It\^o's formula, we have
\begin{align*}
u_n(t,Z_t)&=u_n(0,Z_0)+\int^t_0 (f_n+h_n+g_n)(s, Z_s)\dif s\\
&+\int^t_0\!\!\int_{\mR^d}\left(u_n(s,Z_{s-}+\sigma w)-u_n(s,Z_{s-})\right)\wt{N}(\dif w, \dif s).
\end{align*}
Noting that by \eqref{FA00} with $q'=\infty$ and $\beta'\in(0,\alpha-\frac\alpha q)$,
\begin{align}
\|u_n-u\|_{\mL^\infty_T}&\leq\sup_{t\in[0,T]}\sup_{x,v}\int_{\mR^{2d}}|u(t,x-y,v-w)-u(t,x,v)|\varGamma_n(y,w)\dif y\dif w\no\\
&\leq\|u\|_{\mL^\infty_T(\bC^{\beta'}_\bba)}\int_{\mR^{2d}}(|y|^{\frac1{1+\alpha}}+|w|)\varGamma_n(y,w)\dif y\dif w
\lesssim n^{-\frac{\beta'}{1+\alpha}},\label{Aw3}
\end{align}
we have
$$
\left(u_n(t,Z_t),u_n(0,Z_0)\right)\to\left(u(t,Z_t),u(0,Z_0)\right)\ \ a.s.,\ n\to\infty.
$$
Let
$$
A_n:=\int^t_0\int_{\mR^d}\left(u_n(s,Z_{s-}+\sigma w)-u_n(s,Z_{s-})\right)\wt N(\dif w, \dif s)
$$
and
$$
A:=\int^t_0\int_{\mR^d}\left(u(s,Z_{s-}+\sigma w)-u(s,Z_{s-})\right)\wt N(\dif w, \dif s).
$$
Noting that
\begin{align*}
&|(u_n(s,z+\sigma w)-u_n(s,z))-(u(s,z+\sigma w)-u(s,z))|\\
&\qquad\leq \|u_n(s)-u(s)\|_\infty\wedge(\|\nabla_vu(s)\|_\infty|\sigma w|),
\end{align*}
by the isometry formula of stochastic integrals, we have
\begin{align*}
\mE|A_n-A|^2\lesssim  \int^t_0\!\!\int_{\mR^d}\left(\|u_n(s)-u(s)\|_{\infty}^2\wedge\big(\|\nabla_v u(s)\|_\infty^2|\sigma w|^{2}\big)\right)\nu^{(\alpha)}(\dif w)\dif s.
\end{align*}
Recalling $\nu^{(\alpha)}(\dif w)=c_{d,\alpha}|w|^{-d-\alpha}\dif w$,
by \eqref{Aw3} and the dominated convergence theorem, we get
\begin{align*}
\lim_{n\to\infty}\mE|A_n-A|^2=0.
\end{align*}
Finally, noting that
\begin{align*}
g_n(t,z)=\int_{\mR^{2d}}u(t,x-y,v-w) w\cdot\nabla_y \varGamma_n(y,w)\dif y\dif w,
\end{align*}
since $\varGamma_{n}(y,w)=n^{(2+\alpha)d}\varGamma(ny,n^{1+\alpha}w)$, we have
\begin{align*}
\|g_n(t)\|_{\infty}\lesssim n^{-\alpha}\|u(t)\|_{\infty}
\end{align*}
and
\begin{align*}
\|h_n(t)\|_{\infty}&\lesssim n^{-\beta/(1+\alpha)}\|\DD_{t}\|_{\bC^{\beta}_\bba}\|\nabla_v u(t)\|_{\infty}.
\end{align*}
Hence, by \eqref{FA00} with $(q',\beta')=(\infty,1)$ and the dominated convergence theorem,
$$
\int^t_0(|f_n-f|+|h_n|+|g_n|)(s,Z_s)\dif s\to 0,\ \ a.s.,\ n\to\infty.
$$
Thus we complete the proof of \eqref{SDE34}.
\end{proof}
\br\rm
The stochastic integrals in \eqref{SDE34} and \eqref{SDE304} are martingales by the regularity estimates of $u$ in \eqref{FA00}.
 \er
Now, we can give the
\begin{proof}[Proof of Theorem \ref{thmW}:]
({\bf Existence of weak solutions}) Since  $\DD\in\mL^q_T(\bC^\beta_\bba)$ is continuous in the spatial variable, 
the existence of a weak solution is standard by weak convergence method. We referred to \cite{EK86} and \cite{CZZ21}.

({\bf Weak uniqueness}) Let $Z^1$ and $Z^2$ be two weak solutions of SDE \eqref{MV13} with the same initial distribution. 
For any $f\in  \bC^{\beta}_\bba$, by Lemma \ref{lem44} with $T=t$, we have
$$
\mE \int_0^t f(Z^1_s)\dif s=-\mE u(0,Z^1_0)=-\mE u(0,Z^2_0)=\mE \int_0^t f(Z^2_s)\dif s.
$$
Hence,
$$
\mE f(Z^1_t)=\mE f(Z^2_t).
$$
In particular, $Z^1$ and $Z^2$ have the same one dimensional marginal distribution. Thus,
 as in \cite[Theorem 4.4.3]{EK86}, by a standard method, 
we obtain the uniqueness.

({\bf Strong uniqueness and stability}) 
We use the well-known Zvonkin transformation to prove the strong stability, which automatically imply the pathwise uniqueness.
The strong uniqueness then follows by Yamada-Watanabe's theorem. Below we only prove it for $\alpha\in(1,2)$. 
Recall that
$$
\DD,\wt\DD\in \mL^q_T(\bC^{\frac\gamma{1+\alpha}}_x\cap\bC^\beta_\bba),\ \ 
\beta>1-\tfrac\alpha 2,\ \ \gamma=\alpha+\beta.
$$
By the interpolation inequality \eqref{IA1}, we have for any $\theta\in[0,1]$,
$$
\DD,\wt\DD\in \mL^q_T(\bC^{\theta\gamma,(1-\theta)\beta}_{x,\bba}).
$$
Fix $T>0$ and $\lambda\geq 1$. Let $\bu:=(u_1,\cdots, u_d)$ solve the following backward PDE:
$$
\p_t \bu+(\Delta^{\alpha/2}_v+v\cdot\nabla_x-\lambda)\bu+\DD\cdot\nabla_v \bu+\DD=0,\quad \bu(T)=0.
$$
%Similarly, we have $\wt \bu$.
By \eqref{FA001}, we have for any $\beta'\in[0,\alpha+\frac\alpha{q'}-\frac\alpha q)$ and $\theta\in[0,1]$,
\begin{align}\label{AM3}
\|\bu\|_{\mL^{q'}_{T}(\bC^{\theta\gamma,\beta'+(1-\theta)\beta}_{x,\bba})}\lesssim_C \lambda^{\frac1q+\frac{\beta'}\alpha-1-\frac1{q'}}\|\DD\|_{\mL^q_T(\bC^{\theta\gamma,(1-\theta)\beta}_{x,\bba})}.
\end{align}
Similarly, we define $\wt\bu$ by $\wt\DD$. Then \eqref{AM3} also holds for $\wt H$ and $\wt u$. 
Moreover, by \eqref{FA101} we have
\begin{align}\label{FA102}
\|\bu-\wt \bu\|_{\mL^{q'}_{T}(\bC^{\beta'+\beta}_{\bba})}\lesssim_C 
\lambda^{\frac1q+\frac{\beta'}\alpha-1-\frac1{q'}}\|H-\wt H\|_{\mL^q_T(\bC^{\beta}_{\bba})}.
\end{align} 
In particular, letting $q'=\infty$, $\theta=1$ in \eqref{AM3} and $\beta'\in(1,\alpha-\frac\alpha q)$, 
by \eqref{IA2}, we obtain
$$
\|\nabla \bu\|_{\mL^\infty_T}\vee \|\nabla \wt\bu\|_{\mL^\infty_T}\lesssim \|\bu\|_{\mL^\infty_{T}(\bC^{\gamma,\beta'}_{x,\bba})}\vee \|\wt\bu\|_{\mL^\infty_{T}(\bC^{\gamma,\beta'}_{x,\bba})}\lesssim
\lambda^{\frac1q+\frac{\beta'}\alpha-1},
$$
where the implicit constant is independent of $\lambda\geq 1$.
Since $\frac1q+\frac{\beta'}\alpha<1$, one can choose $\lambda$ large enough so that
$$
\|\nabla \bu\|_{\mL^\infty_T}\vee \|\nabla \wt\bu\|_{\mL^\infty_T}\leq\tfrac12.
$$
Now if we define
$$
\Phi(t,z)=\Phi(t,x,v):=(x,v+\bu(t,x,v)),
$$
then for $z=(x,v)$ and $z'=(x',v')$,
$$
|z-z'|\leq |\Phi(t,z)-\Phi(t,z')|+|\bu(t,z)-\bu(t,z')|\leq
|\Phi(t,z)-\Phi(t,z')|+\tfrac12|z-z'|,
$$
which implies
\begin{align}\label{Am1}
|z-z'|\leq 2|\Phi(t,z)-\Phi(t,z')|.
\end{align}
Similarly, we define $\wt\Phi$ in terms of $\wt\bu$.
On the other hand,  since $\beta>1-\tfrac\alpha 2$, by \eqref{IA2} and \eqref{AM3} with $q'=q$ and $\theta=0,1$,
it is easy to see that for some $\eps_0>0$,
\begin{align}\label{FB07}
\|\nabla \bu\|_{\mL^q_{T}(\bC^{\alpha/2+\eps_0}_v)}\vee \|\nabla \wt\bu\|_{\mL^q_{T}(\bC^{\alpha/2+\eps_0}_v)}<\infty.
\end{align}
By \eqref{SDE34} with $f=-\DD$, one finds that
\begin{align}\label{SDE324}
\Phi(t, Z_t)=\Phi(0, Z_0)+\int^t_0 (V_s,0)\dif s
+\int^t_0\int_{\mR^d}g(s,Z_{s-},w)\wt{N}(\dif w,\dif s),
\end{align}
where 
$$
g(s,z,w):=\Phi(s,z+\sigma w)-\Phi(s,z)=(0,w+\bu(s,z+\sigma w)-\bu(s,z)).
$$
Similarly, we define
$$
\wt g(s,z,w):=\wt \Phi(s,z+\sigma w)-\wt\Phi(s,z)=(0,w+\wt\bu(s,z+\sigma w)-\wt\bu(s,z)).
$$
Let $\wt Z$ be any solution of SDE \eqref{SDE304} corresponding to $\wt\DD$ and initial value $\wt Z_0$. Define
$$
U_t:=Z_t-\wt Z_t,\ \ \Lambda_t:=\Phi(t,Z_t)-\wt\Phi(t,\wt Z_t),\ G_s(w):=g(s,Z_{s-},w)-\wt g(s,\wt Z_{s-},w).
$$
By \eqref{SDE324}, we have
$$
\Lambda_t=\Lambda_0+\int^t_0(V_s-\wt V_s,0)\dif s+\int^t_0\!\!\int_{\mR^d}G_s(w)\wt{N}(\dif w,\dif s).
$$
By BDG's inequality, we have
\begin{align}
\mE\left(\sup_{s\in[0,t]}|\Lambda_s|^2\right)
&\lesssim \mE|\Lambda_0|^2+\int^t_0\mE|V_s-\wt V_s|^2\dif s
+\mE\left(\sup_{s'\in[0,t]}\left|\int^{s'}_0\!\!\int_{\mR^d}G_s(w)\wt{N}(\dif w,\dif s)\right|^2\right)\no\\
&\lesssim \mE|\Lambda_0|^2+\int^T_0\mE|U_s|^2\dif s
+\mE\left(\int^t_0\!\!\int_{\mR^d}|G_s(w)|^2\nu^{(\alpha)}(\dif w)\dif s\right).\label{SW5}
\end{align}
Noting that
\begin{align*}
g(s,z,w)-g(s,\wt z,w)=(z-\wt z)\cdot\int_0^1\nabla_z g(s,\theta z+(1-\theta)\wt z,w)\dif \theta,
\end{align*}
and
$$
\nabla_z g(s,z,w)=
\left(
\begin{array}{cc}
0,&\nabla_x\bu(s,z+\sigma w)-\nabla_x\bu(s,z)\\
0,&\nabla_v\bu(s,z+\sigma w)-\nabla_v\bu(s,z)
\end{array}
\right),
$$
by \eqref{FB07}, we have for $\delta=\frac\alpha 2+\eps_0$,
$$
|g(s,z,w)- g(s,\wt z,w)|\lesssim|z-\wt z|\left(\left(\|\nabla \bu(s)\|_{\bC^\delta_v}|w|^\delta\right)\wedge\|\nabla \bu(s)\|_{\infty}\right).
$$
Similarly, we have
$$
|g(s,z,w)-\wt g(s, z,w)|\lesssim\left(\|\bu(s)-\wt\bu(s)\|_{\bC^\delta_v}|w|^\delta\right)\wedge\|\bu(s)-\wt \bu(s)\|_{\infty}.
$$
Hence,
\begin{align*}
|G_s(w)|&\lesssim|U_s|\left(\left(\|\nabla \bu(s)\|_{\bC^\delta_v}|w|^\delta\right)\wedge\|\nabla \bu(s)\|_{\infty}\right)\\
&\quad+\left(\|\bu(s)-\wt\bu(s)\|_{\bC^\delta_v}|w|^\delta\right)\wedge\|\bu(s)-\wt \bu(s)\|_{\infty}.
\end{align*}
Recalling $\nu^{(\alpha)}(\dif w)=c_{d,\alpha}|w|^{-d-\alpha}\dif w$,  for $\delta>\frac\alpha 2$ and any $a,b>0$, we have
$$
\int_{\mR^d} (a|w|^\delta\wedge b)^2\nu^{(\alpha)}(\dif w)\leq C(a^2+b^2).
$$
Moreover, by \eqref{Am1} we also have
$$
|U_s|\leq 2|\Lambda_s|+2\|\bu(s)-\wt\bu(s)\|_\infty.
$$
Combining the above calculations and by \eqref{SW5}, we get for all $t\in[0,T]$,
\begin{align*}
\mE\left(\sup_{s\in[0,t]}|U_s|^2\right)
&\lesssim \mE|U_0|^2+\|\bu-\wt\bu\|^2_{\mL^\infty_T}+\int^T_0\|\bu(s)-\wt\bu(s)\|^2_{\bC^\delta_v}\dif s\\
&\quad+\int^t_0\Big(1+\|\nabla \bu(s)\|^2_{\bC^\delta_v}\Big)\mE|U_s|^2\dif s.
\end{align*}
By Gronwall's inequality and \eqref{FB07}, we obtain 
$$
 \mE\left(\sup_{s\in[0,T]}|U_s|^2\right)\lesssim\mE|U_0|^2+\|\bu-\wt\bu\|^2_{\mL^\infty_T}+\|\bu-\wt\bu\|^2_{\mL^2_T(\bC^\delta_v)},
$$
which in turn implies \eqref{SW3} by \eqref{FA102}.
\end{proof}

\subsection{Proofs of Theorems \ref{Main0} and \ref{Main1}}\label{Sec4.2H}
Suppose that \eqref{AD067} and \eqref{KB9} hold. Let
\begin{align}\label{Mb0}
\gamma_b>(\sQ-1-\sA_1)\vee0=\sQ-1-(\sQ-1)\wedge\sA_1.
\end{align}
Recall the parameter set $\Theta$ defined in \eqref{Pa1}. Under \eqref{AD067} and \eqref{KB9},
by Theorem \ref{Main}, there is a constant $C_0=C_0(\Theta)>0$ such that
for  any fixed $T\geq 1$ and $\mu_0\in \mZ_{\beta_0}\cap\cP(\mR^{2d})$ with
\begin{align}\label{Pa2}
\|\mu_0\|_{\mZ_{\beta_0}}\leq C_0 T^{-\gamma_b/\alpha}/\kappa_b,
\end{align}
where $\kappa_b$ is from \eqref{KB9},
there is a unique smooth solution $u_T$ to FPE \eqref{FPE00} with $b_T=b\b1_{[0,T]}$ and initial value $u_T(0)=\mu_0$
so that for any $\beta\geq 0$,
$$
\|u_T\|_{\mS_\beta}\leq C_\beta,
$$
where $\mS_\beta$ is defined by \eqref{AS59}.
In the following we fix time $T\geq 1$ and drop the subscript $T$ in $u_T$.

Recalling the notations in Section \ref{Sec3}, we first show the following $L^1$-estimate.
\bl
For any $t\in[0,T]$, it holds that
\begin{align}\label{Aa3}
\|u(t)\|_{\mL^1}\leq 1.
\end{align}
\el
\begin{proof}
We use the stability in Theorem \ref{Main} to show \eqref{Aa3}.
Let $b^n:=b*\varGamma_n$ and $u^n_0:=\mu_0*\varGamma_n$ be the modifying approximation of $b$ and $\mu_0$,
where $\varGamma_n$ is defined in \eqref{Mo5}. By definitions, it is easy to see that
$$
b^n\in\mL^{q_b}_{\bbg_b}(C^\infty_b(\mR^{2d})),\ \ u^n_0\in C^\infty_b(\mR^{2d})\cap L^1(\mR^{2d}),
$$
and
$$
\sup_{n\in\mN}\|b^n\|_{\mL^{q_b}_{\bbg_b}(\mX'_0)}\leq\|b\|_{\mL^{q_b}_{\bbg_b}(\mX'_0)},\ \ 
\sup_{n\in\mN}\|u^n_0\|_{\mZ_{\beta_0}}\leq\|u_0\|_{\mZ_{\beta_0}},
$$
and for any $\eps>0$,
\begin{align}\label{Cx6}
\|b^{n}-b\|_{\mL^{q_b}_{\bbg_b}(\mX_{\eps}')}+\|u^{n}_0-u_0\|_{\mZ_{\beta_0-\eps}}\to0\quad \text{as } n\to\infty.
\end{align}
Let $u^n\in C([0,T];C^2_b(\mR^{2d}))$ be the smooth solution to the following kinetic Fokker-Planck  equation
$$
\p_t u^n=(\Delta^{\frac\alpha 2}_v-v\cdot\nabla_x)u^n-\div_v ((b^n*u^n)u^n),\ u^n(0)=u^n_0.
$$
By Lemma \ref{LEC1} in appendix, we have
\begin{align}\label{Cx5}
\|u^n(t)\|_{\mL^1}\leq\|u^n(0)\|_{\mL^1}=\|\mu_0*\varGamma_n\|_{\mL^1}=1.
\end{align}
By \eqref{AD067}, one can fix $\eps>0$ being small enough so that
\begin{align*}
\sA_0-(\beta_0-\eps)-(\beta_b-\eps)<\sQ+(\beta_0-\eps)\wedge(-1).
\end{align*}
Using $(\beta_0-\eps,\beta_b-\eps)$ to replace $(\beta_0,\beta_b)$ in Theorem \ref{Main} 
and slightly adjusting the smallness constant in \eqref{Pa2}, by \eqref{Cx6} and the stability estimate \eqref{AS609}, 
for each $t\in(0,T]$, one has 
$$
\lim_{n\to\infty}\|u^n(t)-u(t)\|_{\bbp_0}=0.
$$
Thus by \eqref{Cx5} and Fatou's lemma, we conclude \eqref{Aa3}.
\end{proof}

The following lemma improves the regularity estimate in \eqref{SM9}.
\bl
For any $\beta>\beta_0$, there is a constant $C=C(\Theta,\beta, T, \kappa_b,\|\mu_0\|_{\mZ_{\beta_0}})>0$ such that 
\begin{align}\label{AZ3}
\|u(t)\|_{\bB^{\beta,1}_{\bbp_0;\bba}}+
\|u(t)\|_{\bB^{\beta,1}_{\bbp_1;\bba}}\leq Ct^{-\frac{\beta-\beta_0}\alpha}, \ t\in(0,T].
\end{align}
\el
\begin{proof}
For $\theta=\sA_0-2\beta_0-\beta_b,$
by \eqref{AD067} and \eqref{Mb0}, we have
$$
\theta\in(0,\sQ),\ \ \beta_0\leq\sQ-1-\theta,\ \  0\wle\bbg_0+\bbl_0+\bbg_b-(\theta,\sQ-1).
$$
Thus for $i=0,1$, by \eqref{AS1} with $(\bbp',\bbp)=(\bbp_i,\bbp_i)$,
\begin{align*}
\|\cH^b(u,u)\|_{\mL^\infty_0(\bB^{\beta_0,1}_{\bbp_i;\bba})}
&\stackrel{\eqref{P2}}{\lesssim}\|\cH^b(u,u)\|_{\mL^\infty_{\bbg'}(\bB^{\sQ-1-\theta,1}_{\bbp_i;\bba})}
\lesssim \|b_T\|_{\mL^{q_b}_{\bbg_b}(\mX_0')}\|u\|_{\mL^\infty_{\bbl_0}(\mX_0)}\|u\|_{\mL^\infty_{\bbg_0}(\mL^{\bbp_i})}\\
&\stackrel{\eqref{AB2}}{\lesssim}\|b\|_{\mL^{q_b}_0(\mX_0')} T^{\gamma_b/\alpha}\|u\|_{\mS_0}^2\stackrel{\eqref{AS6}}{\lesssim}
\|b\|_{\mL^{q_b}_0(\mX_0')} T^{\gamma_b/\alpha}\|\mu_0\|^2_{\mZ_{\beta_0}}.
\end{align*}
On the other hand, by Lemma \ref{Le215} we have
$$
\|P_\cdot u_0\|_{\mL^\infty_0(\bB^{\beta_0}_{\bbp_i;\bba})}\lesssim \|\mu_0\|_{\bB^{\beta_0}_{\bbp_i;\bba}},\ i=0,1.
$$
Combining the above two estimates, we obtain 
$$
\|u(t)\|_{\bB^{\beta_0}_{\bbp_0;\bba}}+
\|u(t)\|_{\bB^{\beta_0}_{\bbp_1;\bba}}\leq C, \ t\in(0,T],
$$
For general case $\beta>\beta_0$, it follows by \eqref{AS6}
and the real interpolation theorem  (see \cite[Theorem 6.4.5, (1)]{BL76}).
\end{proof}
Define
\begin{align}\label{Aa1}
\DD(s,z):=(b(s)*u(s))(z).
\end{align}
To use Theorem \ref{thmW}, the following lemma provides necessary estimates of $\DD$ for linear SDE.
\bl\label{Le45}
(i) Under \eqref{AD067} and \eqref{KB9}, there are $q>\tfrac{\alpha}{\alpha-1}$ and $\beta>0$ such that
$$
\DD\in\mL^q_T(\bC^\beta_\bba).
$$
(ii) Under \eqref{JC08} and \eqref{Conb}, there are $q>\tfrac{\alpha}{\alpha-1}$ and $\beta>1-\tfrac\alpha2$ such that
$$
\DD\in\mL^q_T(\bC^{\gamma}_x\cap\bC^\beta_\bba),\ \ \gamma:=\tfrac{\beta+\alpha}{1+\alpha}.
$$
(iii) Let $\wt b$ be another kernel satisfying \eqref{KB9} and \eqref{Conb}, and $\wt \mu_0\in\mZ_{\beta_0}$ satisfy \eqref{Pa2}. 
Let $\wt u$ be the 
unique smooth solution of  FPE \eqref{FPE00} corresponding to $\wt b$ and $\wt u_0$. Define $\wt\DD_s(z):=(\wt b(s)*\wt u(s))(z)$. Then 
there are $q>\tfrac{\alpha}{\alpha-1}$, $\beta>0$ and constant $C=C(T,\Theta,q,\beta)>0$ such that
$$
\|\DD-\wt \DD\|_{\mL^q_T(\bC^{\beta}_{\bba})}\lesssim_C
\|\mu_0-\wt \mu_0\|_{\mZ_{\beta_0}}+\|b-\wt b\|_{\mL^{q_b}_0(\mX_0')}.
$$
\el
\begin{proof}
For any $\beta\geq 0$, since $\sA_1\leq\sA_0$,  we have for any $t\in(0,T]$,
\begin{align}\label{SW6}
\|u(t)\|_{\mX_\beta}=\sum_{i=0,1}\|u(t)\|_{\bB^{\beta-\beta_b,1}_{\bar\bbrho_i;\bba}}
\stackrel{\eqref{Sob} }{\lesssim}
\sum_{i=0,1}\|u(t)\|_{\bB^{\beta-\beta_b+\sA_i,1}_{\bbp_i;\bba}}
\stackrel{\eqref{AZ3} }{\lesssim}t^{-\frac{\beta+\sA_0-\beta_0-\beta_b}{\alpha}}.
\end{align}
(i)  Recall $b=b_0+b_1$. For $\beta\geq 0$, by \eqref{Con} we have
$$
\|b(t)*u(t)\|_{\bC^\beta_\bba}=\|b(t)*u(t)\|_{\bB^{\beta}_{\infty;\bba}}\lesssim 
\Big(\|b_0(t)\|_{\bB^{\beta_b}_{\bbrho_0;\bba}}+\|b_1(t)\|_{\bB^{\beta_b}_{\bbrho_1;\bba}}\Big)\|u(t)\|_{\mX_\beta}
=\|b(t)\|_{\mX'_0}\|u(t)\|_{\mX_\beta}.
$$
Since $0<\sA_0-\beta_0-\beta_b<\alpha-\frac\alpha{q_b}-1$ and $q_b>\frac\alpha{\alpha-1}$, 
one can choose $q\in(\frac\alpha{\alpha-1},q_b)$ and $\beta>0$ so that
$$
0<\beta+\sA_0-\beta_0-\beta_b<\tfrac\alpha {q}-\tfrac\alpha{q_b}=:\tfrac\alpha r.
$$
Thus by H\"older's inequality and \eqref{SW6}, we have 
$$
\|H\|_{\mL^{q}_{T}(\bC^\beta_\bba)}=\|b*u\|_{\mL^{q}_{T}(\bC^\beta_\bba)}\lesssim  
\|b\|_{\mL^{q_b}_0(\mX'_0)}\left(\int^T_0t^{-r(\beta+\beta_0-\beta_b+\sA_0)/\alpha}\dif t\right)^{1/r}<\infty.
$$
(ii) Since $0<\sA_0-\beta_0-\beta_b<\tfrac32\alpha-\tfrac\alpha{q_b}-2$, 
one can choose $q>\frac\alpha{\alpha-1}$ and $\beta>1-\frac\alpha2$  so that
$$
0<\beta+\sA_0-\beta_0-\beta_b<\tfrac\alpha q-\tfrac\alpha{q_b}=:\tfrac\alpha r.
$$
Let $\delta=\frac\alpha 2-1$. Fix $i=0,1$. By \eqref{Con}, we have
$$
\|b_i(t)*u(t)\|_{\bB^{\beta}_{\infty;\bba}}
\lesssim \|b_i(t)\|_{\bB^{\beta_b}_{\bbrho_i;\bba}}\|u(t)\|_{\bB^{\beta-\beta_b,1}_{\bar\bbrho_i;\bba}},
$$
and for all $j\geq 0$,
$$
\|\cR^x_j b_i(t)*u(t)\|_{\bB^{\delta+\beta}_{\infty;\bba}}
\lesssim \|\cR^x_j b_i(t)\|_{\bB^{\delta+\beta_b}_{\bbrho_i;\bba}}\|u(t)\|_{\bB^{\beta-\beta_b,1}_{\bar\bbrho_i;\bba}}.
$$
In particular,  by Definition \ref{bs} and \eqref{SW6}, for $i=0,1$,
\begin{align*}
\|b_i(t)*u(t)\|_{\bB^{\alpha-\delta,\delta+\beta}_{\infty;x,\bba}\cap \bB^{\beta}_{\infty;\bba}}
\lesssim \| b_i(t)\|_{\bB^{\alpha-\delta,\delta+\beta_b}_{\bbrho_i;x,\bba}\cap \bB^{\beta_b}_{\bbrho_i;\bba}}t^{-\frac{\beta+\sA_0-\beta_0-\beta_b}{\alpha}}.
\end{align*}
By H\"older's inequality, we get
$$
\|b_i*u\|_{\mL^q_T(\bB^{\alpha-\delta,\delta+\beta}_{\infty;x,\bba}\cap \bB^{\beta}_{\infty;\bba})}
\lesssim \| b_i\|_{\mL^{q_b}_T(\bB^{\alpha-\delta,\delta+\beta_b}_{\bbrho_i;x,\bba}\cap \bB^{\beta_b}_{\bbrho_i;\bba})}
\left(\int^T_0t^{-\frac{r(\beta+\sA_0-\beta_0-\beta_b)\vee0}{\alpha}}\dif t\right)^{1/r}<\infty.
$$
Since $\delta+\beta>0$, by \eqref{IA3}, we have for $\gamma=\frac{\alpha+\beta}{1+\alpha}$,
$$
\|b_i*u\|_{\mL^q_T(\bC^{\gamma}_{x}\cap \bC^{\beta}_{\bba})}=
\|b_i*u\|_{\mL^q_T(\bB^{\gamma}_{\infty;x}\cap \bB^{\beta}_{\infty;\bba})}\lesssim 
\|b_i*u\|_{\mL^q_T(\bB^{\alpha-\delta,\delta+\beta}_{\infty;x,\bba}\cap \bB^{\beta}_{\infty;\bba})}<\infty.
$$
(iii) It follows by the calculations in (i) and \eqref{AS609}.
\end{proof}
\br\rm\label{Re47}
In the situation of Remark \ref{Re219}, by Remark \ref{Re313}, the above (i) still holds.
\er

Now we are in a position to give

\begin{proof}[Proof of Theorem \ref{Main0}]
Let $\alpha\in(1,2]$. Consider the following linear SDE:
$$
\dif Z_t=\bD_t(Z_t)\dif t+\sigma\dif L^{(\alpha)}_t,\ \ Z_0\sim\mu_0,
$$
where $\bH_t(x,v)=(v,H_t(x,v))$ is defined in \eqref{Aa1}.
By Lemma \ref{Le45} and Theorem \ref{thmW}, there is a unique weak/strong solution to the above SDE
under \eqref{AD067}-\eqref{KB9} and \eqref{JC08}-\eqref{Conb}, respectively.
It remains to show that
\begin{align}\label{UN1}
\mu_t(\dif z):=\mbox{Law of }Z_t=u(t,z)\dif z.
\end{align}
Let 
$$
w(t,\dif z):=\mu_t(\dif z)-u(t,z)\dif z\in\cM(\mR^{2d}).
$$
Clearly, by \eqref{Aa3} one has
$$
\|w(t)\|_\cM=\sup_{\|f\|_\infty\leq1}\left|\int_{\mR^{2d}}f(z)w(t,\dif z)\right|\leq 2,
$$ 
and $w$ solves the following FPE with initial value $w(0)=0$ in the distributional sense:
$$
\p_tw=(\Delta^{\frac\alpha 2}_v-v\cdot\nabla_x)w-\div_v (\DD w).
$$
By Duhamel's formula, we have
$$
w(t,z)+\int^t_0P_{t-s}\div_v (\DD_s w_s)\dif s=0,\ t\in(0,T].
$$
Let $P_t^*$ be the adjoint semigroup of $P_t$ (see \eqref{SCL11}). 
By the definition of total variation norm, 
\begin{align*}
\|w(t)\|_\cM&=\sup_{\|f\|_\infty\leq1}\left|\int^t_0\int_{\mR^{2d}}P_{t-s}\div_v (\DD_s w_s)(z) f(z)\dif z\dif s\right|\\
&=\sup_{\|f\|_\infty\leq1}\left|\int^t_0\int_{\mR^{2d}}\DD_s(z)\cdot \nabla_vP^*_{t-s}f(z)w(s,\dif z)\dif s\right|\\
&\leq\sup_{\|f\|_\infty\leq1}\int^t_0\|\DD_s\|_\infty\|\nabla_vP^*_{t-s}f\|_\infty\|w(s)\|_\cM\dif s.
\end{align*}
Let $q>\frac\alpha{\alpha-1}$ be as in (i) of Lemma \ref{Le45} and $\frac1r+\frac1q=1$. Since
$\|\nabla_vP^*_{t-s}f\|_\infty\lesssim (t-s)^{-\frac1\alpha}\|f\|_\infty$ (see Lemma \ref{Le10}), by H\"older's inequality we 
further have
\begin{align*}
\|w(t)\|_\cM
&\lesssim\int^t_0(t-s)^{-\frac1\alpha}\|\DD_s\|_\infty\|w(s)\|_\cM\dif s\
\leq \|\DD\|_{\mL^q_T(\mL^\infty)}
\left(\int^t_0(t-s)^{-\frac{r}\alpha}\|w(s)\|^r_\cM\dif s\right)^{1/r},
\end{align*}
which, by $\frac r\alpha<1$ and Gronwall's inequality, yields $w(t)=0$.  Thus we obtain \eqref{UN1}.

(i) The smoothness of density and short time singularity estimates follow by \eqref{AZ3}.

(ii) The large time decay estimate follows by Theorem \ref{thm230}.

(iii) The stability of the density follows by \eqref{AS609}.

(iv) The strong stability follows by Theorem \ref{thmW} and (iii) of Lemma \ref{Le45}.

(v) The global existence follows by Theorems \ref{TH1}, \ref{TH2} and Remark \ref{Re318}.
\end{proof}

\br\rm\label{Re48}
Under the assumptions of Remark \ref{Re219}, by Remark \ref{Re47}, it is easy to see that the weak well-posedness still holds.
\er

Next we give the
\begin{proof}[Proof of Theorem \ref{Main0}]
We only show the two-sides estimates of the density since the others are completely the same as Theorem \ref{Main0}.
Let $u$ be the unique solution of \eqref{FPE01} and define $H_t(x)=b(t)*u(t)(x)$.
Consider the following linear SDE:
$$
\dif X_t=H_t(X_t)\dif t+\dif L^{(\alpha)}_t.
$$
Since $\DD\in\mL^q_T(\bC^\beta)$ for some $q>\frac\alpha{\alpha-1}$ and $\beta>0$,
it is well known that for each starting point $X_0=x$, $X_t$ has a density $p(t,x,y)$ enjoying 
two-sides estimates: for $\alpha\in(1,2)$ (see \cite{CHXZ18}),
$$
p(t,x,y)\asymp t(t^{1/\alpha}+|x-y|)^{-d-\alpha},
$$
and for $\alpha=2$ (see \cite[Theorem 1.1]{CHXZ18}),
$$
C^{-1}_1t^{-d/2} \e^{-\frac{|x-y|^2}{C_2t}}\leq p(t,x,y)\leq 
C_1t^{-d/2}\e^{-\frac{C_2|x-y|^2}{t}}.
$$
For general starting point $X_0\sim\mu_0$, we have
$$
\rho(t,y)=\int_{\mR^d}p(t,x,y)\mu_0(\dif x).
$$
The desired estimates now follows.
\end{proof}

\section{Application to fractional Vlasov-Poisson-Fokker-Planck equation}
Let $d\geq 3$.
In this section, we apply Theorem \ref{thm230} to the following VPFP equation in $\mR^{2d}$:
\begin{align}\label{VPFP}
\p_t f=\Delta^{\alpha/2}_v f-v\cdot\nabla_x f+\gamma \nabla U\cdot\nabla_vf,
\end{align}
where $\alpha\in(1,2]$ and $\gamma=\pm1$ stands for the attractive or repulsive force in physics, respectively,
\begin{align*}
U(t,x)=(|\cdot|^{2-d}*\<f\>(t))(x),\quad \<f\>(t,x):=\int_{\mR^d}f(t,x,v)\dif v.
\end{align*}
Let 
$$
K(x):=\nabla |\cdot|^{2-d}(x),\ \ \wt K(x,v):=K(x).
$$
With these notations, we have
$$
E(x):=\nabla U*\<f\>(x)=K*\<f\>(x)=(\wt K*f)(x,v).
$$
By Lemma \ref{lemA02} and Proposition \ref{Pro1}, we have
\bl\label{Le51}
Let $\varrho_0, \varrho_1\in[1,\infty]$ and $\bbrho_i:=(\varrho_i,\infty)$, $i=0,1$.
If either $\varrho_0>\frac{d}{d-1}$ and $\varrho_1\in[1,\infty]$ or  
$\varrho_0<\frac{d}{d-1}$ and $\varrho_1\in(\tfrac d{d-1},\infty]$,
then for any $\beta_b\leq (1+\alpha)(\tfrac{d}{\varrho_0}+1-d)$,
$$
\wt K\in\bB^{\beta_b}_{\bbrho_0;\bba}+\bB^{\beta_b}_{\bbrho_1;\bba}.
$$ 
Moreover, for $s_0=1+\frac\alpha 2$ and $s_1=\frac\alpha2-1+\beta_b$,
$$
\wt K\in\bB^{s_0,s_1}_{\bbrho_0;x,\bba}.
$$
\el

As an application of  Theorem \ref{thm230}, we have
\bt\label{thm:VP01}
Let $\alpha\in(1,2]$, $p_0\in[2,\infty]$, $\bbp_0:=(p_0,1)$ and $\beta_0\in(-\alpha,0)$.
Suppose that
\begin{align}\label{CD7}
\tfrac{(1+\alpha)d}{p_0}=2\alpha+\beta_0,\ \beta_0\in(-1,0),\ \ \tfrac{(1+\alpha)d}{p_0}<2\alpha+1+2\beta_0,\ \beta_0\in(-\alpha,-1].
\end{align}
For any $f_0\in \bB^{\beta_0}_{\bbp_0;\bba}\cap\bB^{\beta_0}_{\bbb1;\bba}$,
there is a constant $C_0=C_0(\alpha,\beta_0,p_0,d, \|f_0\|_{\bB^{\beta_0}_{\bbb1;\bba}})>0$ such that if
\begin{align}\label{SMA1}
\|f_0\|_{\bB^{\beta_0}_{\bbp_0;\bba}}\leq C_0,
\end{align}
then there is a unique global smooth solution $f$ to VPFP \eqref{VPFP} such that for any $\beta\geq 0$,
\begin{align}\label{HZBB03}
\sup_{t\in(0,1]}\left(t^{\frac{\beta-\beta_0}\alpha}\big(\|f(t)\|_{\bB^{\beta,1}_{\bbp_0;\bba}}
+\|f(t)\|_{\bB^{\beta,1}_{\bbb1;\bba}}\big)\right)<\infty,
\end{align}
and for any $p\in[2,\infty]$ and $\bbp:=(p,1)$, 
\begin{align}\label{HZBB02}
\sup_{t\geq 1}\left(t^{\frac{d(1+\alpha)}{\alpha}\cdot(1-\frac 1p)}\|f(t)\|_{\bB^{\beta,1}_{\bbp;\bba}}\right)<\infty.
\end{align}
In particular, 
\begin{align}\label{HZAA05}
\sup_{t\geq 1}\left(t^{\frac{(\alpha+1)(d-1)}{\alpha}}\|E(t)\|_\infty\right)<\infty.
\end{align}
%%%%%%%%%%%%%%%%%%%%%%%%%%%%%%%%
\et
\begin{proof}
To apply Theorem \ref{thm230}, we need to first choose suitable parameters so that {\bf (H$_0$)} holds.
Let 
%$\varrho_1>\frac{d}{d-1}$ and
$$
\bbrho_0=(1,\infty), \bbrho_1=(2,\infty), q_b=\infty,\ \beta_b=1+\alpha,\ \bbp_1:=(2,1).
$$
With these choice of parameters, it is easy to see that
$$
\sA_0:=\bba\cdot(\tfrac d{\bbp_0}+\tfrac d{\bbrho_0}-\bbd)=\tfrac{(1+\alpha)d}{p_0},\ \
\sA_1:=\bba\cdot(\tfrac d{\bbp_1}+\tfrac d{\bbrho_1}-\bbd)=0.
%(1+\alpha)d\tfrac{(\bar\varrho_1-p_1)\vee0}{p_1\bar\varrho_1}.
$$
By \eqref{CD7}, one sees that $\beta_b\not=\sA_i-\beta_0$, $i=0,1$  and
$$
0<\sA_0-\beta_0-\beta_b\leq\alpha-1,\ \ \sA_0-\beta_0-\beta_b<\alpha+\beta_0,
$$
which means that  {\bf (H$_0$)} holds.
Moreover, $\bba\cdot\tfrac d{\bbrho_1}=\tfrac{(1+\alpha)d}2>\alpha$. Thus by Lemma \ref{Le51} and 
Theorem \ref{thm230}, there is a constant $C_0>0$
depending on $\|f_0\|_{\bB^{\beta_0}_{\bbb1;\bba}}$ such that if
\begin{align}\label{SMA2}
\|f_0\|_{\bB^{\beta_0}_{\bbp_0;\bba}}+\|f_0\|_{\bB^{\beta_0}_{\bbp_1;\bba}}\leq C_0,
\end{align}
then there is a unique smooth solution $f$ so that \eqref{HZBB03} and
\eqref{HZBB02} for $\bbp=\bbp_0,\bbp_1$ hold. 
Let $\theta\in[0,1]$ be defined by $\frac{1-\theta}{p_0}+\theta=\frac1{2}$.
Note that by the interpolation inequality \eqref{Sob},
$$
\|f_0\|_{\bB^{\beta_0}_{\bbp_1;\bba}}\lesssim \|f_0\|^{1-\theta}_{\bB^{\beta_0}_{\bbp_0;\bba}}\|f_0\|^\theta_{\bB^{\beta_0}_{\bbb1;\bba}}.
$$
By adjusting the small constant $C_0$ in \eqref{SMA1}, one sees that \eqref{SMA1} implies \eqref{SMA2}.

Now we devote to showing \eqref{HZBB02} for all $\bbp=(p,1)$, where $p\in[2,\infty]$.
If $p_0=\infty$, then by the interpolation inequality \eqref{Sob}, we immediately have
\begin{align}\label{HZBB09}
\|f(t)\|_{\bB^{\beta,1}_{\bbp;\bba}}\leq \|f(t)\|^{\frac 2p}_{\bB^{\beta,1}_{(2,1);\bba}}
\|f(t)\|^{1-\frac 2p}_{\bB^{\beta,1}_{(\infty,1);\bba}}
\lesssim t^{-\frac{d(1+\alpha)}{\alpha}\cdot(1-\frac 1p)},\ \ t\geq1.
\end{align}
If $p_0<\infty$, then  for $\beta>\frac d{p_0}$, by the Sobolev embedding \eqref{Sob1} and \eqref{HZBB02} for $\bbp=\bbp_0$, we have
$$
\|f(t)\|_{\bB^{\beta}_{(\infty,1);\bba}}\lesssim\|f(t)\|_{\bB^{\beta,1}_{\bbp_0;\bba}}
\lesssim t^{-\frac{d(1+\alpha)}{\alpha}\cdot(1-\frac 1{p_0})},\ \ t\geq1.
$$
In particular, one can choose $t_0\geq 1$ large enough so that 
$$
\|f(t_0)\|_{\bB^{\beta}_{(\infty,1);\bba}}\leq C_0.
$$
Starting from the time $t_0$ and with the choice of $p_0=\infty$, we then deduce \eqref{HZBB02} by \eqref{HZBB09}.

Finally, for \eqref{HZAA05}, by (ii) of Lemma \ref{lemA02} with $\gamma=d-1$ and $p=1$, we have
for $t\ge1$,
\begin{align*}
\|E(t)\|_\infty&=
\|\nabla U*\<f\>(t)\|_\infty\lesssim \|\<f\>(t)\|_{1}^{\frac{1}d}\|\<f\>(t)\|_\infty^{\frac {d-1}d}\\
&\lesssim\|f(t)\|_{\mL^1}^{\frac{1}d}\|f(t)\|_{\bB^{0,1}_{(\infty,1);\bba}}^{\frac {d-1}d}
\lesssim t^{-\frac{(\alpha+1)(d-1)}{\alpha}}.
\end{align*}
The proof is complete.
\end{proof}
\br\label{Rk53H}\rm
Theorem \ref{thm:VP01} extends the well-posedness theory for VPFP equations to the distribution initial data  
(see \cite{VD90}, \cite{Bou93}, \cite{CS97}, \cite{OS00}).
Notice that when $\alpha=2$, in three dimensional case, Carrillo and Soler \cite{CS97} studied the well-posedness of VPFP equation
with initial data being in Morrey spaces.
More precisely, under more restricted smallness assumption on initial data,
$$
\|f_0\|_{\cM}+\|f_0\|_{M\wt L_{9/4}}\leq\eps,
$$
where  $M\wt L_p$ is some Morrey space (see \cite[(1.4)]{CS97}), they showed global well-posedess and obtained the following decay estimate  for the force field $E(t,x)$,
$$
\|E(t)\|_\infty\lesssim t^{-1/2},\ \ t\ge1.
$$
Note that one can check that for any $p_0\in(9/4,3)$,
$$
M\wt L_{9/4}\cap\cM\subset\bB^{0}_{(9/4,1);\bba}\cap
\bB^{0}_{\bbb1;\bba}\subset\bB^{9/p_0-4}_{(p_0,1);\bba}\cap\bB^{9/p_0-4}_{\bbb1;\bba}.
$$
Later, under smallness condition $\|f_0\|_{L^1}+\|f_0\|_{L^1_v(L^\infty_x)}\leq\eps$, 
 Ono and Strauss \cite{OS00} showed the following optimal rate of large time decay 
\begin{align*}
\|E(t)\|_\infty\lesssim t^{-3(d-1)/2},\ \ t\ge1,
\end{align*}
which coincides with \eqref{HZAA05}.
It is pointed out that requiring small $\|f_0\|_{\cM}$ rules out the application for the associated mean-field SDE since we always have
$\|f_0\|_{\cM}\equiv1$ for SDEs. In fact, we have paid much efforts to drop the smallness of $\|f_0\|_{\cM}$ in Theorem \ref{thm230}.
%%%%%%%%%%%%%%%%%%%%%%%%%%%%%%%%%%%%%
\er

Now we consider the following mean-field SDE:
\begin{align}\label{SDE9}
\left\{
\begin{aligned}
&X_t=X_0+\int_0^t V_s\dif s,\\
&V_t=V_0+\int^t_0(K*\mu_{X_s})(X_s)\dif s+L^{(\alpha)}_t.
\end{aligned}
\right.
\end{align}
By Theorem \ref{Main0} we have
\bt\label{thm:VP02}
Let $\alpha\in(1,2]$, $\beta_0\in(-\alpha,0)$, $p_0\in[2,\infty]$, $\bbp_0:=(p_0,1)$ and $(X_0,V_0)\sim\mu_0\in\cP(\mR^{2d})$.
\begin{enumerate}[{\bf (A)}]
\item Suppose that
$$
\tfrac{(1+\alpha)d}{p_0}<2\alpha+1+\beta_0+\beta_0\wedge(-1).
$$
Then there is a constant $C_0=C_0(\alpha,\beta_0,p_0,d)>0$ such that if
\begin{align}\label{SM2}
\|\mu_0\|_{\bB^{\beta_0}_{\bbp_0;\bba}}\leq C_0,
\end{align}
then there is a unique weak solution $Z_t=(X_t,V_t)$ to SDE \eqref{SDE9} such that 
for each $t>0$, $Z_t$ admits a smooth density $f(t)$  that satisfies \eqref{HZBB03} and \eqref{HZBB02}.
\item Suppose that 
\begin{align*}
\tfrac{(1+\alpha)d}{p_0}<\tfrac 52\alpha-2+\beta_0.
\end{align*}
Then under \eqref{SM2}, there is a unique strong solution to SDE \eqref{SDE9}.
\end{enumerate}
Moreover, when $d=3$ and $\alpha=2$, and if $\mE|V_0|^2<\infty$,  then
the above {\bf (A)}  and {\bf (B)} still hold without smallness condition \eqref{SM2}.
\et
\begin{proof}
Note that by \eqref{Cx3},
$$
\|\mu_0\|_{\bB^0_{\bbb1;\bba}}\leq C(d,\alpha).
$$
As the choice of parameters in Theorem \ref{thm:VP01}, 
{\bf (A)} and {\bf (B)} are direct applications of Theorem \ref{Main0}. 
When $d=3$ and $\alpha=2$, in order to obtain a global solution without smallness condition \eqref{SM2},
we shall use a result of \cite[Theorem 5.1]{OS00}. More precisely, it suffices to find a small time $t_0>0$ so that
$$
\|f(t_0)\|_{L^1_v(L^\infty_x)}:=\int_{\mR^d}\sup_{x\in\mR^d}|f(t_0,x,v)|\dif v<\infty
$$
 and
 $$
 \int_{\mR^{2d}}|v|^2f(t_0,x,v)\dif x\dif v<\infty.
 $$
For any $f_0\in\bB^{\beta_0}_{\bbp_0;\bba}$, by Theorem \ref{Cor35},
one can find a $t_0>0$ small enough so that 
there is a local solution for SDE \eqref{SDE9} up to time $t_0$. Moreover, for any $\beta\geq 0$,
$$
\sup_{t\in(0,t_0]}\left( t^{\beta/\alpha}\|f(t)\|_{\mZ_{\beta+\beta_0}}\right)<\infty.
$$
By \eqref{SDE9} we have
\begin{align*}
(\mE|V_t|^2)^{1/2}&\leq (\mE|V_0|^2)^{1/2}+\left(\mE\Big|\int_0^t (K*\mu_{X_s})(X_s)\dif s\Big|^2\right)^{1/2}+\sqrt2(\mE|W_t|^2)^{1/2}\\
&\leq(\mE|V_0|^2)^{1/2}+\int_0^t\|K*\mu_{X_s}\|_\infty\dif s+\sqrt{2 t},
\end{align*}
which implies by (i) of Lemma \ref{Le45} that for any $t\in[0,t_0]$,
$$
\int_{\mR^{2d}}|v|^2f(t,x,v)\dif x\dif v=\mE|V_t|^2<\infty.
$$
By \eqref{AB2} and \eqref{Sob1}, for any $t\in(0,t_0]$ and $\beta>3(1+\alpha)(1-\frac1{p_0})$, we have
$$
\|f(t)\|_{L^1_v(L^\infty_x)}\leq\sum_{j\geq 0}\|\cR^\bba_jf(t)\|_{L^1_v(L^\infty_x)}=
\|f(t)\|_{\bB^{0,1}_{(\infty,1);\bba}}\lesssim \|f(t)\|_{\bB^{\beta}_{\bbp_0;\bba}}<\infty,
$$
and for $\beta>3(1+\alpha)(1-\frac1{p_0})+3$, 
$$
\|f(t)\|_{L^\infty_{x,v}}\lesssim \|f(t)\|_{\bB^{\beta}_{\bbp_0;\bba}}<\infty.
$$
Now starting from time $t_0/2$, we can use \cite[Theorem 5.1]{OS00} to find a unique global smooth solution $\wt f$ for \eqref{VPFP} with initial value $\wt f(0)=f(t_0/2)$. Since $f$ is smooth on $[\frac{t_0}2,t_0]$, by the uniqueness we have $f(t+t_0/2)=\wt f(t)$ on $[0,\frac{t_0}2]$.
Thus we can patch up the solutions and conclude the proof.
\end{proof}
\br\rm
For $\alpha\in(1,2)$, it seems to be an open question to obtain the global solution for fractional VPFP equation
since there is no finite second-order moment for $\alpha$-stable process.
\er

\iffalse
When $\alpha=2$, we have the following global existence.
Define a Lyapunov function
$$
H(t,x,v):=|v|^2/2+U(t,x)
$$
Note that
$$
\p_t U(t,x)=\mE(V_t\cdot \nabla G(x-X_t)).
$$
By It\^o's formula, we have
\begin{align*}
H(t, Z_t)=H(0,Z_0)+\int^t_0(\p_s U(s,X_s)+V_s\cdot\nabla_x U(s,X_s)+F(s,X_s)\cdot V_s)\dif s
\end{align*}
\fi

\section{Application to fractional Navier-Stokes equations and SQG equations}

In this section we apply our results to the 2D and 3D fractional Navier-Stokes equations:
\begin{align}\label{NSE}
\p_t u=\Delta^{\alpha/2}u+u\cdot\nabla u+\nabla p,\ \ \div u=0, 
\end{align}
where $\alpha\in(1,2]$, $u$ is the velocity field and $p$ is the pressure.
For 2D fractional Navier-Stokes equations, we show the global well-posedness 
when the initial vorticity is in $\bB^{-\frac{1+\alpha}2+}_{\infty}:=\cup_{\eps>0}\bB^{\eps-\frac{1+\alpha}2}_{\infty}$. 
For 3D fractional Navier-Stokes equations, we show the global well-posedness when the 
initial vorticity is in $\bB^{-\frac{1+\alpha}2+}_{\infty}$ and small.
To the best of our knowledge, these results are new for fractional Navier-Stokes equations. 
Moreover, we also consider the SQG equation in Example 2 in the introduction.

We first consider the 2D fractional NSEs with initial vorticity in Besov space $\bB^{\beta_0}_{p_0}$.
Let $w={\rm curl}(u)$ be the vorticity of $u$ that satisfies the following scalar equation in $\mR^2$:
\begin{align}\label{AK1}
\p_t w=\Delta^{\alpha/2} w+u\cdot\nabla w,
\end{align}
where the velocity field $u$ can be recovered from $w$ by the Biot-Savart law (cf. \cite{MB02}):
$$
u(x)=K_2*w(x)=\frac1{2\pi}\int_{\mR^2}\frac{(y_2-x_2, x_1-y_1)}{|x-y|^{2}} w(y)\dif y.
$$
By (iv) of Lemma \ref{lemA02}, one has
\begin{align}\label{AV8}
K_2\in \bB^{\beta}_{\varrho_0}+\bB^{\beta}_{\varrho_1},\ \ \beta=\tfrac 2 {\varrho_0}-1,\ \varrho_0\in[1,2),\ \varrho_1\in(2,\infty].
\end{align}
By Theorem \ref{thm02} we have
\bt\label{Th51}
 Let $\alpha\in(1,2]$, $p_0\in[1,\infty]$  and $\beta_0\in(-\frac{\alpha+1}2,0)$ satisfy
\begin{align}\label{AV6}
\tfrac{2}{p_0}\leq\alpha+\beta_0,\ \beta_0\in(-1,0),\ \ \tfrac{2}{p_0}<\alpha+1+2\beta_0,\ \beta_0\in(-\tfrac{\alpha+1}2,-1].
\end{align}
For any $w_0\in \cup_{q\in[1,\infty)}(\bB^{\beta_0,q}_{p_0}\cap \bB^{\beta_0,q}_{1})$, there is a unique global smooth solution $w
$ to vorticity equation \eqref{AK1} on $\mR_+$ so that for any $T>0$ and $\beta\geq 0$,
$$
\sup_{t\in(0,T]}\left(t^{-(\beta-\beta_0)/\alpha}\big(\|w(t)\|_{\bB^{\beta,1}_{p_0}}+\|w(t)\|_{\bB^{\beta,1}_1}\big)\right)<\infty.
$$
\et
\begin{proof}
We only need to check that  under \eqref{AV6}, the assumptions in Theorem \ref{thm02} hold for $b=K_2$ and $q_b=\infty$. 
Let $\varrho_0=1$, $\varrho_1=\infty$ and $\beta_b=\tfrac 2 {\varrho_0}-1=1$. By \eqref{AV8}, one sees that the assumptions of Theorem \ref{thm02} hold. Moreover, since $\div K_2=0$, we have the global solution.
\end{proof}

\br\rm
If we take $p_0=\infty$ in \eqref{AV6}, then for any $w_0\in\bB^{-\frac{1+\alpha}2+}_\infty\cap \bB^{-\frac{1+\alpha}2+}_{1}$, there is 
a unique global smooth solution  to equation \eqref{AK1} on $\mR_+$. Since the two dimensional Brownian white noise locally
belongs to $\bB^{-1-}_\infty$, one can apply Theorem \ref{Th51} to the 2D fractional Navier-Stokes equations with initial vorticity being the Brownian white noise.
\er

Now we consider the following mean-field SDE in $\mR^2$ with kernel $K_2$:
\begin{align}\label{DDSDE0}
\dif X_t=(K_2*\mu_{X_t})(X_t)\dif t+\dif L^{(\alpha)}_t,\ \ X_0\sim \mu_0\in\cP(\mR^2).
\end{align}
By Theorem \ref{Main1} with $\varrho_1=\infty$, we have
\bt
Let $\alpha\in(1,2]$, $p_0\in[1,\infty]$  and $\beta_0\in(-\frac{\alpha+1}2,0)$ satisfy
$$
\left\{
\begin{aligned}
&\tfrac{2}{p_0}<\alpha+1+\beta_0+\beta_0\wedge(-1),&\ \ {\rm weak\ solution},\\
&\tfrac{2}{p_0}<\tfrac{3\alpha}2-1+\beta_0,&\ \ {\rm strong\ solution}.
\end{aligned}
\right.
$$
For any $\mu_0\in \cP(\mR^2)\cap\bB^{\beta_0}_{p_0}$, there exists a unique weak/strong solution $X_t$ for SDE \eqref{DDSDE0}
so that for each $t>0$, $X_t$ has a smooth density $\mu_{X_t}(\dif x)=\rho_t(x)\dif x$ 
with regularity that for any $T, \beta\geq 0$,
$$
\sup_{t\in(0,T]}\left(t^{\frac{\beta}{\alpha}}\|\rho_t\|_{\bB_{p_0}^{\beta+\beta_0}}\right)<\infty.
$$
\et
%\subsection{3D- fractional NSEs  with initial vorticity in $\bB^{-\beta_0}_{p_0}$}

Next we consider the following mean-field SDE in $\mR^2$ associated with SQG \eqref{SQG1}:
\begin{align}\label{DDSDE30}
\dif X_t=(\sR*\mu_{X_t})(X_t)\dif t+\dif L^{(\alpha)}_t,\ \ X_0\sim \mu_0\in\cP(\mR^2),
\end{align}
where $\sR=(-x_2,x_1)/|x|^3$. Note that by (iii) of Lemma \ref{lemA02},
$$
\sR\in \bB^{\beta}_{\varrho_0},\ \ \beta=\tfrac 2 {\varrho_0}-2,\ \varrho_0\in(1,\infty].
$$
By Theorem \ref{Main1} with $\varrho_1=\varrho_0$, we similarly have
\bt
Let $\alpha\in(1,2]$, $p_0\in[1,\infty]$  and $\beta_0\in(-\alpha,0)$ satisfy
$$
\left\{
\begin{aligned}
&\tfrac{2}{p_0}<\alpha+\beta_0+\beta_0\wedge(-1),&\ \ {\rm weak\ solution},\\
&\tfrac{2}{p_0}<\tfrac{3\alpha}2-2+\beta_0,&\ {\rm strong\ solution}.
\end{aligned}
\right.
$$
For any $\mu_0\in \cP(\mR^2)\cap\bB^{\beta_0}_{p_0}$, there exists a unique weak/strong solution $X_t$ for SDE \eqref{DDSDE30}
so that for each $t>0$, $X_t$ has a smooth density $\mu_{X_t}(\dif x)=\rho_t(x)\dif x$  with regularity that for any $T, \beta\geq 0$,
$$
\sup_{t\in(0, T]}\left(t^{\frac{\beta}{\alpha}}\|\rho_t\|_{\bB_{p_0}^{\beta+\beta_0}}\right)<\infty.
$$
\et
\br\rm
From the condition, one sees that for strong solutions, we need to assume $\alpha\in(\tfrac43,2]$.
\er

Finally we consider the 3D fractional NSEs with initial vorticity in non-homogeneous space $\bB^{\beta_0}_{p_0}$.
Let $w={\rm curl}(u)$ be the vorticity of $u$. Then 
\begin{align}\label{AK33}
\p_t w=\Delta^{\frac\alpha2} w-u\cdot\nabla w+w\cdot\nabla u=\Delta^{\frac\alpha2} w+\div(u\otimes w-w\otimes u),
\end{align}
where for two divergence free vector fields $w_1,w_2$,
$$
\div (w_1\otimes w_2):=w_1\cdot\nabla w_2.
$$ 
By the Biot-Savart law, one knows that (cf. \cite{MB02})
$$
u(x)=K_3*w(x)=\int_{\mR^3}K_3(x-y)w(y)\dif y,
$$
where
$$
K_3(x) h=\frac1{4\pi}\frac{x\times h}{|x|^3},\ h\in\mR^3.
$$
Note that by (iv) of Lemma \ref{lemA02},
\begin{align}\label{AV2}
K_3\in \bB^{\beta}_{\varrho_0}+\bB^{\beta}_{\varrho_1},\ \ \beta=\tfrac 3 {\varrho_0}-2,\ \varrho_0\in[1,\tfrac32),\ \varrho_1\in(\tfrac32,\infty].
\end{align}
By Duhamel's formula, we have
$$
w(t)=P_t w_0+\int^t_0 P_{t-s}\div((K_3*w)\otimes w-w\otimes (K_3*w))\dif s.
$$
The following result can be proved by using similar argument as in proving Theorem \ref{thm230}.
\bt\label{thm010}
 Let $\alpha\in(1,2]$, $p_0\in[2,\infty]$  and $\beta_0\in(-\alpha,0)$ satisfy
\begin{align}\label{AV0}
\tfrac{3}{p_0}\leq\alpha+\beta_0,\ \beta_0\in(-1,0),\ \ \tfrac{3}{p_0}<\alpha+1+2\beta_0,\ \beta_0\in(-\tfrac{\alpha+1}2,-1].
\end{align}
For any initial vorticity $w_0\in \bB^{\beta_0}_{p_0}\cap \bB^{\beta_0}_{1}$,
there is a constant $C_0=C_0(\alpha, \beta_0,p_0, \|w_0\|_{\bB^{\beta_0}_1})>0$ such that if
\begin{align}\label{AV9}
\|w_0\|_{\bB^{\beta_0}_{p_0}}\leq C_0,
\end{align}
then there is a unique global smooth solution $w$  to equation \eqref{AK33} on $\mR_+$
so that for any  $\beta\geq 0$,
\begin{align}\label{HZAA01}
\sup_{t>0}\left(t^{(\beta-\beta_0)/\alpha}\big(\|w(t)\|_{\bB^{\beta,1}_{p_0}}+\|w(t)\|_{\bB^{\beta,1}_1}\big)\right)<\infty
\end{align}%%%%%%%%%%%%%%%%%%%%%%%%%%%%%%%
and for any $p\in[2,\infty]$ and $\beta\geq 0$, there is a constant $C=C(\alpha,p,\beta,\beta_0,p_0,\|w_0\|_{\bB^{\beta_0}_1})>0$ such that
\begin{align}\label{HZAA02}
\sup_{t\geq 1}\left(t^{(3-\frac3{p})/\alpha}\|w(t)\|_{\bB^{\beta,1}_{p}}\right)<\infty.
\end{align}
%%%%%%%%%%%%%%%%%%%%%%%%%%%%%%%%%%%%%%%
\et
\begin{proof}
For given $\beta_0,p_0$ satisfying \eqref{AV0}, by \eqref{AV2},
 one can choose $\varrho_0=1$, $\varrho_1=2$, $\beta_b=\tfrac 3 {\varrho_0}-2=1$, $q_b=\infty$  so that
 $b=K_3\in\bB^{\beta_b}_{\varrho_0}+\bB^{\beta_b}_{\varrho_1}$, and \eqref{AD606} holds and
$$
\tfrac1{p_0}+\tfrac1{\varrho_0}\geq 1,\ \ \ \tfrac 3{\varrho_1}>\alpha-1.
$$
Define
$$
\cH_t(u,w):=\int^t_0 P_{t-s}\div((K_3*u)\otimes w-w\otimes (K_3*u))\dif s.
$$
As in the proof of Theorem \ref{thm230}, one can show that there is a small constant $C_0>0$ such that under condition \eqref{AV9}, there is a unique global smooth solution $w$ satisfying
$$
w(t)=P_tw_0+\cH_t(w,w),
$$
and  \eqref{HZAA01}  and \eqref{HZAA02} for $p=2,p_0$ hold.

Now we devote to showing \eqref{HZAA02} for all $p\in[2,\infty]$. 
Fix $\beta\geq 0$. Note that it is completely the same as proving Theorem \ref{thm230}
that if $p_0=\infty$, then 
\begin{align*}
\sup_{t\geq 1}\Big(t^{3/\alpha} \|w(t)\|_{\bB^{\beta,1}_\infty}+t^{3/(2\alpha)}\|w(t)\|_{\bB^{\beta,1}_2}\Big)<\infty,
\end{align*}
which implies by interpolation inequality \eqref{Sob},
\begin{align}\label{HZBB06}
\|w(t)\|_{\bB^{\beta,1}_{p}}\leq \|w(t)\|^{\frac 2p}_{\bB^{\beta,1}_2}\|w(t)\|^{1-\frac 2p}_{\bB^{\beta,1}_\infty}
\lesssim t^{-(3-\frac3{p})/\alpha},\ \ t\geq1.
\end{align}
If $p_0<\infty$,  by \eqref{HZAA01} with $\beta=\frac3{p_0}$ and the Sobolev embedding inequality \eqref{Sob1}, we have
\begin{align*}
\|w(t)\|_\infty\lesssim \|w(t)\|_{\bB^{3/p_0,1}_{p_0}}\lesssim t^{-(3-3/p_0)/\alpha},\ \ t\geq 1,
\end{align*}
and
$$
\sup_{t\geq 1}\|w(t)\|_{\bB^{0,1}_1}<\infty.
$$
In particular, we can choose $t_0$ large enough so that for the small constant $C_0$ in \eqref{AV9},
$$
\|w(t_0)\|_{\bB^{\beta_0}_{\infty;\bba}}\leq C\|w(t_0)\|_\infty\leq C_0.
$$
Now starting from $t_0$, we can use \eqref{HZBB06} to conclude the proof of \eqref{HZAA02}.
%%%%%%%%%%%%%%%%%%%%%%%%%%%%%%%%%%%%
\end{proof}
\br\rm
When $\alpha=2$, it is by no means that Theorem \ref{thm010} is new. For fractional Navier-Stokes equation with $\alpha\in(1,2)$, 
Wu \cite[Theorem 6.1]{Wu05} showed the global well-posedness  for small initial velocity in homogeneous Besov spaces.
As one knows, there are vast literatures about 3D Navier-Stokes equations. Here we only want to emphasize that
we present a general framework that contains Navier-Stokes equations as one application.
\er

\section{Application to fractional kinetic porous medium equation with viscosity}
In this section we consider the following fractional kinetic porous medium equation with viscosity
\begin{align}\label{KD1}
\p_t f=\Delta_v^{\alpha/2}f-v\cdot\nabla_xf+\div_v((\nabla_v\Delta^{-s}_v f) f ),
\end{align}
where $s\in(0,1]$ and $d\geq 3$. Note that
$$
(\nabla_v\Delta^{-s}_v f)(x,v)= K*f(x,v)
$$
with
$$
K(x,v):=c_{d,s}\delta_0(\dif x)\nabla_v |v|^{2s-d},
$$
where $c_{d,s}>0$ is a normalized  constant and $\delta_0$ is the Dirac measure at point $0$.

The following lemma provides the necessary regularity estimate for $K$.
\bl
Fix $p_0\in[1,\infty]$. Let $\bbp=(p_x,p_v)$ with $p_x\in[1,\infty]$ and $p_v\in[p_0,\infty]$.
For any $\beta_0<\frac d {p_0}-\frac d{p_v}$, there is a constant $C=C(\alpha, d,p_0,\beta_0,\bbp)>0$ such that for any 
$U(v)\in\bB^{\beta_0}_{p_0}$,
\begin{align}\label{aK1}
\|\delta_0\otimes U\|_{\bB^{\beta_0+d-\frac d {p_0}-\bba\cdot(\bbd-\frac{d}{\bbp})}_{\bbp;\bba}}
\lesssim_C \|U\|_{\bB^{\beta_0}_{p_0}}.
\end{align}
In particular, for $s\in(0,1]$ and $\bbrho_i=(1,\varrho_i)$, $i=0,1$, where $\varrho_0\in[1,\infty]$ and $\varrho_1\in(\frac d{d-2s},\infty]$, 
\begin{align}\label{HA1}
\delta_0(\dif x)\nabla|v|^{2s-d}\in\bB^{\beta_b}_{\bbrho_0;\bba}+\bB^{\beta_b}_{\bbrho_1;\bba},
\end{align}
where $\beta_b<(2s-1-d+\frac{d}{\varrho_0})\wedge 0$.
\el
\begin{proof}
By definition \eqref{Ph0} we have
\begin{align*}
\cR^\bba_j (\delta_0\otimes U)(x,v)=(\check{\phi}^\bba_j*_vU)(x,v)
=\sum_{k\geq 0}(\check{\phi}^\bba_j*_v\cR^v_kU)(x,v),
\end{align*}
where $*_v$ stands for the convolution with respect to the variable $v$ and $\cR^v_k U=\check\phi_k*U$.
Noting that
$$
\widehat{(\check{\phi}^\bba_j*_v\cR^v_kU)}(\xi,\eta)=\phi^\bba_j(\xi,\eta) \phi_k(\eta) \hat U(\eta),
$$
and by \eqref{Cx8} and \eqref{Cx9},
$$
\phi^\bba_j(\xi,\eta)\phi_k(\eta)=0,\ \ k>j+2,
$$
by Young's inequality, scaling of $\check{\phi}^\bba_j$ and Bernstein's inequality \eqref{Ber}, we have for any $j\geq 0$,
\begin{align*}
\|\cR^\bba_j (\delta_0\otimes U)\|_{\bbp}
&\leq\sum_{k\leq j+2}\|\check{\phi}^\bba_j*_v\cR^v_kU\|_{(p_x,p_v)}
\leq\sum_{k\leq j+2}\|\check{\phi}^\bba_j\|_{(p_x,1)}\|\cR^v_kU\|_{p_v}\\
&=\sum_{k\leq j+2}2^{j(1+\alpha)(d-\frac d{p_x})}\|\check{\phi}^\bba_1\|_{(p_x,1)}\|\cR^v_kU\|_{p_v}\\
&\lesssim2^{j(1+\alpha)(d-\frac d{p_x})}\sum_{k\leq j+2}2^{k(\frac dp-\frac d{p_v}-\beta)}\|U\|_{\bB^{\beta_0}_{p_0}}\\
&\lesssim2^{j(\bba\cdot(\bbd-\frac{d}{\bbp})+\frac d{p_0}-d-\beta_0)}\|U\|_{\bB^{\beta_0}_{p_0}},
\end{align*}
where the last step is due to $\beta_0<\frac d{p_0}-\frac d{p_v}$.
Thus we get \eqref{aK1}.
Now, fix $s\in(0,1]$.
%Let $U(x,v):=c_{d,s}\delta_0(\dif x)|v|^{2s-d}$. Since $K=\nabla_v U$, by Bernstein's inequality, it suffices to show $U\in \bB^{\beta+1}_{\bbp;\bba}$. 
Let $U(v):=\nabla_v|v|^{2s-d}$.  For any 
$\varrho_1\in(\frac{d}{d-2s},\infty]$, by (iv) of Lemma \ref{lemA02}, we have
$$
U\in \bB^{2s-1}_1+\bB^{2s-1}_{\varrho_1}.
$$
Now for any $\varrho_0\in[1,\infty]$,
if we take $\beta_0<(2s-1)\wedge(d-\frac d{\varrho_0})$, $p_0=1$ and $\bbp=(1,\varrho_0)$  in \eqref{aK1},
then for $U_0\in \bB^{2s-1}_1\subset\bB^{\beta_0}_1$,
$$
\delta_0\otimes U_0\in \bB^{\beta_0+d/\varrho_0-d}_{\bbrho_0;\bba}.
$$
Similarly,  for any $\varrho_1\in(\frac{d}{d-2s},\infty]$, if we take $\beta_0=\beta_b<0$,
$p_0=\varrho_1$ and $\bbp=(1,\varrho_1)$ in \eqref{aK1}, 
then for $U_1\in \bB^{2s-1}_{\varrho_1}\subset\bB^{\beta_b}_{\varrho_1}$,
$$
\delta_0\otimes U_1\in \bB^{\beta_b}_{\bbrho_1;\bba}.
$$
Thus we obtain \eqref{HA1} by noting $\beta_0+\frac{d}{\varrho_0}-d<(2s-1)\wedge(d-\frac d{\varrho_0})+\frac{d}{\varrho_0}-d
=(2s-1-d+\frac{d}{\varrho_0})\wedge 0$.
%$U\in \bB^{\beta+1}_{\bbp;\bba}$. The proof is complete.
\end{proof}

The associated density-distribution dependent SDE to \eqref{KD1} is taken as the form:
\begin{align}\label{SDE69}
\left\{
\begin{aligned}
&X_t=X_0+\int_0^t V_s\dif s,\\
&V_t=V_0+\int^t_0(K*f_{Z_s})(X_s,V_s)\dif s+L^{(\alpha)}_t,
\end{aligned}
\right.
\end{align}
where $f_{Z_s}$ stands for the distribution density of $Z_s=(X_s,V_s)$.
By Theorem \ref{Main0} we have
\bt\label{thm:PME}
Let $s\in(0,1]$, $\beta_0\in(-\alpha,0)$, $\bbp_0=(p_{0x},p_{0v})\in[1,\infty]^2$
and $(X_0,V_0)\sim\mu_0\in\cP(\mR^{2d})$.
Suppose that
\begin{align}\label{CC8}
(1+\alpha)\tfrac{d}{p_{0x}}<\alpha+(2s-1-\tfrac{d}{p_{0v}})\wedge 0+\beta_0+\beta_0\wedge(-1).
\end{align}
If $\mu_0\in\bB^{\beta_0}_{\bbp_0;\bba}$ satisfies the following smallness condition: 
for some $C_0=C_0(s,\alpha,d,\bbp_0,\beta_0)>0$,
\begin{align}\label{aK2}
\|\mu_0\|_{\bB^{\beta_0}_{\bbp_0;\bba}}\leq C_0,
\end{align}
then there is a unique weak solution $Z_t=(X_t,V_t)$ to SDE \eqref{SDE69} in the class that
for each $t>0$, $Z_t$ admits a smooth density $f(t)$ so that for any $T>0$ and $\beta\geq 0$,
\begin{align}\label{SM44}
\sup_{t\in(0,T]}\left(t^{(\beta-\beta_0)/\alpha}\|f(t)\|_{\bB^{\beta,1}_{\bbp_0;\bba}}\right)<\infty.
\end{align}
Moreover, if $\bbp_0=(p_0,p_0)$, then we can drop \eqref{aK2}.
\et
\begin{proof}
Under \eqref{CC8},  one can choose $\beta_b<(2s-1-\frac{d}{p_{0v}})\wedge 0$ and $\varrho_0=(1,\frac{p_{0v}}{p_{0v}-1})$ so that
$$
\bba\cdot(\tfrac d{\bbp_0}+\tfrac d{\bbrho_0}-\bbd)-\beta_0-\beta_b<\alpha+\beta_0\wedge(-1).
$$
Thus by \eqref{HA1} and Theorem \ref{Main0},  as in the proof of Theorem \ref{thm:VP01}, we conclude the proof.

Without \eqref{aK2}, by Theorem \ref{Cor35}, there is a local solution to SDE \eqref{SDE69} on some small time interval
$[0,T_0]$, where $T_0$ depends on the norm  $\|\mu_0\|_{\bB^{\beta_0}_{\bbp_0;\bba}}$.
Next we want to extend the solution to arbitrary large time when $\bbp_0=(p_0,p_0)$. 
Fix $t_0\in(0,T_0)$. Let us start from time $t_0$ to make some estimates.
Suppose that $\psi:[0,\infty]\to[0,\infty)$ is a $C^2$-convex function. 
Since $f\in C([t_0,T_0]; C^2_b(\mR^{2d}))$, by \eqref{AC11} we have
$$
\int\psi(f(t))
\leq\int\psi(f(t_0)) -\int^t_{t_0}\int \nabla_v\Delta^{-s}_v f\cdot\nabla_v f\,\psi''(f) f,
$$
where the integral $\int$ is taken over $\mR^{2d}$.
Define
$$
G(r):=\int^r_0 s\psi''(s)\dif s.
$$
It is easy to see that
$$
\int\psi(f(t))\leq\int\psi(f(t_0)) -\int^t_{t_0}\int \nabla_v\Delta^{-s}_v f\cdot\nabla_v G(f)
=\int\psi(f(t_0)) +\int^t_{t_0}\int \Delta^{1-s}_v f\, G(f).
$$
By \eqref{Mb3} and the change of variable, one sees that
\begin{align*}
&\int \Delta^{1-s}_v f\, G(f)
=\lim_{\eps\downarrow 0}\int_{\mR^d}\int_{|w|\geq \eps}G(f(v))\frac{f(v+w)-f(v)}{|w|^{d+2(1-s)}}\dif w\dif v\\
&\qquad=\lim_{\eps\downarrow 0}
\frac12\int_{\mR^d}\int_{|w|\geq\eps}(G(f(v))-G(f(v+w)))\frac{f(v+w)-f(v)}{|w|^{d+2(1-s)}}\dif w\dif v\\
&\qquad=-\lim_{\eps\downarrow 0}\frac12\int_{\mR^d}\int_{|w|\geq\eps}\int^1_0G'(s f(v)+(1-s)f(v+w))\dif s\frac{(f(v+w)-f(v))^2}{|w|^{d+2(1-s)}}\dif w\dif v.
\end{align*}
Since $G'(r)=r\psi''(r)\geq 0$ for $r\geq 0$, we get
$$
\int \Delta^{1-s}_v f\, G(f)\leq0.
$$
Therefore,
$$
\int\psi(f(t))\leq\int\psi(f(t_0)).
$$
In particular,
$$
\|f(t)\|_{p_0}\leq \|f(t_0)\|_{p_0}.
$$
As in the proof of Theorem \ref{TH1}, we can extend the solution to arbitrary time.
\end{proof}

\begin{appendix}
\renewcommand{\thetable}{A\arabic{table}}
\numberwithin{equation}{section}

\section {Several technical lemmas}

The following two elementary lemmas are used to show the Schauder estimate in time-weighted spaces for kinetic semigroups.
\bl\label{LA1}
For any $\gamma>\beta>0$, there is a constant $C=C(\gamma,\beta)\geq 1$ such that for all $\lambda>0$,
$$
\sum_{j\geq 1}2^{j\beta}((2^{-j\gamma}\lambda)\wedge 1)\asymp_C\lambda^{\frac\beta\gamma}
\int^\lambda_0 t^{-\frac\beta\gamma-1}(t\wedge 1)\dif t.
$$
\el
\begin{proof}
Noting that for $j\in\mN$ and $s\in[j-1,j]$,
$$
2^{s\beta}((2^{-(s+1)\gamma}\lambda)\wedge 1)
\leq 2^{j\beta}((2^{-j\gamma}\lambda)\wedge 1)\leq 2^{(s+1)\beta}((2^{-s\gamma}\lambda)\wedge 1),
$$
we have
\begin{align*}
\sI_0:=\int^\infty_02^{s\beta}((2^{-(s+1)\gamma}\lambda)\wedge 1)\dif s\leq
\sum_{j\geq 1}2^{j\beta}((2^{-j\gamma}\lambda)\wedge 1)
\leq \int^\infty_02^{(s+1)\beta}((2^{-s\gamma}\lambda)\wedge 1)\dif s=:\sI_1.
\end{align*}
By the change of variables, we have
$$
\sI_0=\frac{\lambda^{\frac\beta\gamma} 2^{-\beta}}{\gamma\ln 2}\int^\lambda_0 t^{-\frac\beta\gamma-1}(t\wedge 1)\dif t,\ \ 
\sI_1=\frac{\lambda^{\frac\beta\gamma} 2^\beta}{\gamma\ln 2}\int^\lambda_0 t^{-\frac\beta\gamma-1}(t\wedge 1)\dif t.
$$
Thus we obtain the desired estimate.
\end{proof}
\bl\label{lemA2}
For any $\vartheta\in[0,\infty)$ and $\gamma_0,\gamma_1\in(-\infty,1)$,
there is a constant $C=C(\vartheta, \gamma_0,\gamma_1)>0$ such that for all $\lambda,t>0$,
\begin{align*}
\int_0^t \left[(\lambda (t-s))^{-\vartheta}\wedge1\right](1\wedge s)^{-\gamma_0}(1\vee s)^{-\gamma_1}\dif s
\leq C  
\ell_\vartheta(\lambda t)(1\wedge t)^{1-\gamma_0}(1\vee t)^{1-\gamma_1},
\end{align*}
where 
$$
\ell_\vartheta(\lambda):=(\lambda\vee 1)^{-\vartheta\wedge 1}(1+\b1_{\vartheta=1}\ln(1\vee\lambda)).
$$
\el
\begin{proof}
We first show the following estimate: for any $\vartheta\in[0,\infty)$ and $\gamma<1$,
\begin{align}\label{SX2}
\sI(\lambda,t):=\int_0^t \left[(\lambda (t-s))^{-\vartheta}\wedge1\right]s^{-\gamma}\dif s\lesssim \ell_\vartheta(\lambda t)t^{1-\gamma},\ \ \forall\lambda, t>0.
\end{align}
By the change of variable,  we have
$$
\sI(\lambda,t)=t^{1-\gamma}\sI(t\lambda,1).
$$
So, it suffices to prove \eqref{SX2} for $t=1$. Note that
\begin{align*}
\sI(\lambda,1)
&\leq 2^{\gamma\vee 0}\int_0^{\frac12} \left[(\lambda s)^{-\vartheta}\wedge1\right]\dif s+\left[(\tfrac\lambda2)^{-\vartheta}\wedge1\right]\int_{\frac12}^1(1-s)^{-\gamma}\dif s\\
&=2^{\gamma\vee0}\lambda^{-1}\int_0^{\frac\lambda2} \left[s^{-\vartheta}\wedge1\right]\dif s+\left[(\tfrac\lambda2)^{-\vartheta}\wedge1\right]\int_0^{\frac12}s^{-\gamma}\dif s.
\end{align*}
If $0\leq\vartheta\not=1$, then
\begin{align*}
\sI(\lambda,1)\lesssim \lambda^{-1}\wedge 1+\lambda^{-\vartheta}\wedge 1\leq 2(\lambda^{-(\vartheta\wedge 1)}\wedge 1)
=2(\lambda\vee 1)^{-\vartheta\wedge 1}.
\end{align*}
If $\vartheta=1$, then
$$
\sI(\lambda,1)\lesssim (\lambda^{-1}\wedge 1)\ln(1\vee\lambda)+\lambda^{-1}\wedge 1\lesssim
(\lambda\vee1)^{-1}\ln(1\vee\lambda).
$$
Hence,
$$
\sI(\lambda,1)\lesssim\ell_\vartheta(\lambda).
$$
%$$x^{-(\vartheta\wedge 1)}\wedge 1\le  x^{\gamma-\beta-1},\ \ x>0.$$
For $t\in(0,2]$, we have
\begin{align*}
&\int_0^t \left[(\lambda (t-s))^{-\vartheta}\wedge1\right](1\wedge s)^{-\gamma_0}(1\vee s)^{-\gamma_1}\dif s\\
&\qquad\lesssim\int_0^t \left[(\lambda (t-s))^{-\vartheta}\wedge1\right]s^{-\gamma_0}\dif s
\lesssim_C \ell_\vartheta(\lambda t)t^{1-\gamma_0}.
\end{align*}
For $t>2$, we have
\begin{align*}
&\int_0^t \left[(\lambda (t-s))^{-\vartheta}\wedge1\right](1\wedge s)^{-\gamma_0}(1\vee s)^{-\gamma_1}\dif s\\
&\qquad\leq \int_0^1 \left[(\lambda (t-s))^{-\vartheta}\wedge1\right] s^{-\gamma_0}\dif s
+\int_1^t \left[(\lambda (t-s))^{-\vartheta}\wedge1\right] s^{-\gamma_1}\dif s\\
&\qquad\leq \left[(\lambda (t-1))^{-\vartheta}\wedge1\right] \int_0^1 s^{-\gamma_0}\dif s
+\int_0^t \left[(\lambda (t-s))^{-\vartheta}\wedge1\right]  s^{-\gamma_1}\dif s\\
&\qquad\lesssim \left[(\lambda t)^{-\vartheta}\wedge1\right]+\ell_\vartheta(\lambda t)t^{1-\gamma_1}
\lesssim \ell_\vartheta(\lambda t)(1\vee t)^{1-\gamma_1}.
\end{align*}
The proof is complete.
\end{proof}
Below we recall some well-known results about singular kernels mentioned  in the introduction. 
For the readers' convenience, we give full proofs. 
\bl\label{lemA02}
Let  $\gamma\in(0,d)$ and $K:\mR^d\to \mR^N$ be a measurable function with that for some $C>0$,
\begin{align*}
|K(x)|\le C|x|^{\gamma-d},\ \ x\in\mR^d.
\end{align*}
\begin{enumerate}[(i)]
\item Let $1<q<p<\infty$ satisfy 
\begin{align*}
\tfrac1p+\tfrac\gamma{d}=\tfrac1q.
\end{align*}
Then there is a constant $C=C(p,q,d,\gamma)>0$ such that for any $f\in L^q(\mR^d)$,
\begin{align}\label{AW1}
\|K*f\|_p\le C\|f\|_q.
\end{align}
\item Let $p\in[1,d/\gamma)$ and $\theta:=p\gamma/d$. Then there is a constant $C=C(p,d,\gamma)>0$ 
such that for any $f\in L^p(\mR^d)\cap L^\infty(\mR^d)$,
\begin{align}\label{AW2}
\|K*f\|_\infty\le C\|f\|_p^\theta\|f\|_\infty^{1-\theta}.
\end{align}
\item For any $p\in(\frac{d}{d-\gamma},\infty]$,  it holds that 
\begin{align}\label{HH1}
K\in\bB^{d/p-d+\gamma,\infty}_{p},\ \ \nabla K\in\bB^{d/p-d+\gamma-1,\infty}_{p}.
\end{align}
\item If  $K$ is a $(d-\gamma)$-order smooth homogenous function on $\mR^d\setminus\{0\}$, i.e.,
$K(\lambda x)=\lambda^{\gamma-d}K(x)$, $\lambda>0$,
then for any $p_0\in[1,\frac{d}{d-\gamma})$ and $p_1\in(\frac{d}{d-\gamma},\infty]$,  it holds that
\begin{align}\label{HH2}
K\in\bB^{d/p_0-d+\gamma,\infty}_{p_0}+\bB^{d/p_0-d+\gamma,\infty}_{p_1}.
\end{align}
\item Let $K(x)=\nabla|x|^{\gamma-d}$ and $K_\eps(x)=\nabla(|x|\vee\eps)^{\gamma-d}$ for $\eps>0$. Then 
for any $p\in(\frac{d}{d-\gamma},\infty]$ and $r\in[1,\frac{d}{d-\gamma})$, there is constant $C=C(r,d,p,\gamma)>0$
such that for all $\eps>0$,
\begin{align}\label{AW3}
\|K-K_\eps\|_{\bB^{\beta,\infty}_{p}}\lesssim \eps^{\frac dr-d+\gamma},\ \ \beta:=\tfrac dp-\tfrac dr.
\end{align}
\end{enumerate}
\el
\begin{proof}
(i) Let $r=d/(d-\gamma)$. By \cite[Theorem 1.5]{BCD11}, we have
\begin{align*}
\|K*f\|_p\lesssim \|K\|_{r,\infty}\|f\|_q,
\end{align*}
where
\begin{align*}
\|K\|_{r,\infty}:=\sup_{\lambda>0} \lambda^rm\big\{x: |K(x)|>\lambda\big\}\lesssim \sup_{\lambda>0} \lambda^rm\big\{x\in \mR^d: C|x|^{\gamma-d}>\lambda\big\}<\infty.
\end{align*}
\\
(ii) Let $p'=p/(p-1)$. For any $\lambda>0$, by H\"older's inequality we have
\begin{align*}
|K*f(x)|&\leq \int_{|x-y|\lesssim \lambda}\frac{|f(y)|}{|x-y|^{d-\gamma}}\dif y
+\int_{|x-y|> \lambda}\frac{|f(y)|}{|x-y|^{d-\gamma}}\dif y\\
&\lesssim \|f\|_\infty \lambda^{\gamma}+\|f\|_p\left(\int_{|x-y|> \lambda}\frac{1}{|x-y|^{(d-\gamma) p'}}\dif y\right)^{1/p'}\\
&\lesssim \|f\|_\infty \lambda^{\gamma}+\|f\|_p \lambda^{d/p'-d+\gamma},
\end{align*}
In particular, taking $\lambda=\|f\|_p^{p/d}\|f\|_\infty^{-p/d}$, we get
\begin{align*}
|K*f(x)|\lesssim \|f\|_p^\theta\|f\|_\infty^{1-\theta}.
\end{align*}
(iii) For $p\in(\frac{d}{d-\gamma},\infty)$, letting $\tfrac1p+\tfrac\gamma{d}=\tfrac1q$,
by \eqref{AW1} we have for any $j\geq 1$,
$$
\|\cR_jK\|_{p}\lesssim \|\check\phi_j\|_q\lesssim 2^{j(d-\frac d{q})}\|\check{\phi}_1\|_q
=2^{j(d-\frac d{p}-\gamma)}\|\check{\phi}_1\|_q.
$$
For $p=\infty$, by \eqref{AW2} we have
$$
\|\cR_jK\|_\infty\leq \|\check\phi_j\|^{\frac\gamma d}_1\|\check\phi_j\|^{1-\frac\gamma d}_\infty
=\|\check\phi_1\|^{\frac\gamma d}_1\|\check\phi_1\|^{1-\frac\gamma d}_\infty 2^{j(d-\gamma)}.
$$
Hence, $K\in\bB^{d/p-d+\gamma,\infty}_{p}$. By Bernstein's inequality, we have
$\nabla K\in\bB^{d/p-d+\gamma-1,\infty}_{p}$. 
\\
(iv) Now, let $\chi$ be a smooth function with $\chi=1$ on $B_1$
and $\chi=0$ on $B^c_2$. Define
$$
K_1(x):=\chi(x)K(x),\ \ K_2(x):=(1-\chi(x)) K(x).
$$
Fix $p_0\in[1,\frac{d}{d-\gamma})$ and $p_1\in(\frac{d}{d-\gamma},\infty]$.
It is easy to see that for any $m\in\mN_0$,
$$
\nabla^ mK_2\in L^{p_1}\Rightarrow K_2\in \bB^{d/p_0-d+\gamma,\infty}_{p_1}.
$$ 
%Moreover, by \cite[page 108, Proposition 2.93]{BCD11}, we have
%$$
%K_1\in \dot\bB^{d/p_0-\gamma,\infty}_{p_0}\subset\bB^{d/p_0-\gamma,\infty}_{p_0}.
%$$
Next we show $K_1\in\bB^{d/p_0-d+\gamma,\infty}_{p_0}$. Note that
$$
\|\cR_0K_1\|_{p_0}\leq C\|K_1\|_{p_0}<\infty.
$$
For $j\geq 1$, since $K_1=K-K_2$, by  \cite[Proposition 2.21]{BCD11} we have
$$
\|\cR_j K_1\|_{p_0}\leq \|\cR_j K\|_{p_0}+\|\cR_j K_2\|_{p_0}\leq C2^{j(d-\gamma-\frac d {p_0})}+\|\cR_j K_2\|_{p_0},
$$
and
\begin{align*}
\|\cR_j K_2\|_{p_0}&=\|\Delta^{-d}\cR_j\Delta^d K_2\|_{p_0}\lesssim 2^{-2d j}\|\Delta^d K_2\|_{p_0}\\
&\lesssim 2^{j(d-\gamma-\frac d {p_0})}\Big(\|1_{1\leq |x|\leq 2} G\|_{p_0}+\|1_{|x|\geq 1} \nabla^{2d}K\|_{p_0}\Big),
\end{align*}
where $G$ is a bounded function on $\{x: 1\leq |x|\leq 2\}$. Hence,
$$
\|K_1\|_{\bB^{d/p_0-d+\gamma,\infty}_{p_0}}= \sup_{j\geq 0}2^{j(\frac d {p_0}-d+\gamma)}\|\cR_j K_1\|_{p_0}<\infty.
$$
(v) For $p\in(\frac{d}{d-\gamma},\infty)$, letting $1+\tfrac1p=\tfrac 1{r}+\tfrac1q$, by Young's inequality we have
$$
\|\cR_j(K-K_\eps)\|_{p}\leq \|\nabla\check\phi_j\|_q\||\cdot|^{\gamma-d}-(|\cdot|\vee\eps)^{\gamma-d}\|_r
\lesssim 2^{j(d-\frac dq)} \eps^{\frac dr-d+\gamma}.
$$
This completes the proof.
\end{proof}

We also need the following simple lifting proposition.
\bp\label{Pro1}
Suppose $K\in\bB^{s,q}_{\varrho}$ for some $s\in\mR$ and $q,\varrho\in[1,\infty]$. Let 
$$
K_1(x,v)=K(x),\ K_2(x,v)=K(v),\ \bbrho_1=(\varrho,\infty),\ \bbrho_2=(\infty, \varrho).
$$ 
Then we have
$$
\|K_1\|_{\bB^{(1+\alpha)s,q}_{\bbrho_1;\bba}}\asymp \|K\|_{\bB^{s,q}_{\varrho}}\asymp\|K_2\|_{\bB^{s,q}_{\bbrho_2;\bba}}.
$$ 
\ep
\begin{proof}
We only prove $\|K_1\|_{\bB^{(1+\alpha)s,q}_{\bbrho_1;\bba}}\asymp \|K\|_{\bB^{s,q}_{\varrho}}$. By definition we have for $j\ge0$,
$$
\cR^\bba_j K_1(x,v)=\int_{\mR^{2d}}\check{\phi}^\bba_j(x-x',v-v')K(x')\dif x'\dif v'=h_j*K(x),
$$
where
$$
h_j(x):=\int_{\mR^{d}}\check{\phi}^\bba_j(x,v)\dif v.
$$
Note that for $j\geq 1$,
$$
\hat h_j(\xi)=2^{j(1+\alpha)d}(2\pi)^{-d}\int_{\mR^{2d}}\e^{ix\cdot\xi}\check{\phi}^\bba_1 (2^{j(1+\alpha)}x,v)\dif x\dif v=\phi^\bba_1(2^{-(1+\alpha)j}\xi,0),
$$
where the support of $\phi^\bba_1(2^{-(1+\alpha)j}\cdot,0)$ is contained in
\begin{align*}
\Big\{ 2^{j-1}<|\xi|^{1/(1+\alpha)}<2^{j+1}\Big\}=\Big\{2^{(1+\alpha)(j-1)}<|\xi|<2^{(1+\alpha)(j+1)}\Big\}.
\end{align*}
In particular, an equivalent norm of the Besov space $\bB^{s,q}_{\varrho}$ is then given by
$$
\|f\|_{\bB^{s,q}_{\varrho}}\asymp \left(\sum_{j\geq0}\big(2^{ j(1+\alpha)s}\|h_j*f\|_{\varrho}\big)^q\right)^{1/q}.
$$
Hence,
$$
\|K_1\|_{\bB^{s,q}_{\bbrho_1;\bba}}^q=\sum_{j\geq0}2^{j(1+\alpha)s q}\|\cR^\bba_jK_1\|_{\bbrho_1}^q
= \sum_{j\geq0}\big(2^{ j(1+\alpha)s}\|h_j*K\|_{\varrho}\big)^q\asymp\|K\|_{\bB^{s,q}_{\varrho}}^q.
$$
The proof is complete.
\end{proof}

\section{Basic properties of Besov spaces}

\begin{proof}[Proof of Lemma \ref{lemB1}]

(i) By definition we have
$$
\|f\|_\bbp\leq\sum_{j\geq 0}\|\cR^\bba_jf\|_\bbp=\|f\|_{\bB^{0,1}_{\bbp;\bba}}
$$
and
$$
\|f\|_{\bB^{0,\infty}_{\bbp;\bba}}=\sup_{j\geq 0}\|\cR^\bba_jf\|_\bbp\lesssim \|f\|_\bbp.
$$

(ii) For $j\geq 0$, by \eqref{KJ2} and Young's inequality, we have
$$
\|\cR^{\bba}_j (f*g)\|_{\bbp}=\|(\cR^{\bba}_jf)*(\wt\cR^{\bba}_jg)\|_{\bbp}\leq\|\cR^{\bba}_jf\|_{\bbp_1}\|\wt\cR^{\bba}_jg\|_{\bbp_2},
$$
where $\wt\cR^{\bba}_jg=\sum_{|i-j|\leq 2}\cR^\bba_ig$.
Hence,
\begin{align*}
\|f*g\|_{\bB^{\beta,q}_{\bbp;\bba}}&=\left(\sum_{j}(2^{j\beta}\|\cR^{\bba}_j(f*g)\|_{\bbp})^q\right)^{1/q}\\
&\leq \left(\sum_{j}(2^{j\beta_1}\|\cR^{\bba}_jf\|_{\bbp_1}2^{j\beta_2}\|\wt\cR^{\bba}_jg\|_{\bbp_2})^q\right)^{1/q}\\
&\leq 5\|f\|_{\bB^{\beta_1,q_1}_{\bbp_1;\bba}}\|g\|_{\bB^{\beta_2,q_2}_{\bbp_2;\bba}}.
\end{align*}

(iii) For $j\geq 0$, we define $S_j=\sum_{k\leq j-1}\cR^\bba_j$ with convention $\cR^\bba_{-1}=0$ and $S_{-1}=0$. Recall the following Bony decomposition:
\begin{align*}
fg=f\prec g+ f\circ g+g\prec f,
\end{align*}
where
\begin{align*}
f\prec g:=\sum_{j\ge 0}S_{j-1}f\cR^{\bba}_j g\quad \text{and}\quad f\circ g:=\sum_{j\geq 0}\cR^\bba_jf\wt\cR^{\bba}_jg.
\end{align*}
Here $\wt\cR^\bba_j$ is defined in \eqref{KJ2}.
Noting that
\begin{align}\label{AP00}
\cR_k\(S_{j-1}f\cR^{\bba}_j g\)=0,\quad \forall |k-j|>3,
\end{align}
we have for $k\geq 0$,
\begin{align*}
&\|\cR^\bba_k(f\prec g)\|_{\bbp}\leq\sum_{j\ge 0}\|\cR_k\(S_{j-1}f\cR^{\bba}_j g\)\|_\bbp
=\sum_{|j-k|\le 3}\|\cR_k^\bba\(S_{j-1}f\cR^{\bba}_j g\)\|_\bbp\\
&\quad\lesssim
\sum_{|j-k|\le 3}\|S_{j-1}f\cR^{\bba}_j g\|_\bbp\leq \sum_{|j-k|\le 3}\|S_{j-1}f\|_{\bbp_1}\|\cR^{\bba}_j g\|_{\bbp_2}
\lesssim\|f\|_{\bbp_1} \sum_{|j-k|\le 3}\|\cR^{\bba}_j g\|_{\bbp_2}.
\end{align*}
Hence,
\begin{align*}
\|f\prec g\|_{\bB^{s,q}_{\bbp;\bba}}&\lesssim \|f\|_{\bbp_1}\left(\sum_{k\ge-1}2^{skq}\sum_{|j-k|\le 3}\|\cR^{\bba}_j g\|^q_{\bbp_2}\right)^{1/q}\lesssim \|f\|_{\bbp_1}\|g\|_{\bB^{s,q}_{\bbp_2;\bba}}.
\end{align*}
By symmetry, we also have
$$
\|f\prec g\|_{\bB^{s}_{\bbp;\bba}}\lesssim\|f\|_{\bB^{s,q}_{p_1;\bba}}\|g\|_{\bbp_2}.
$$
Moreover, noting that
$$
\cR^\bba_k( \cR^{\bba}_jf\wt\cR^{\bba}_j g)=0,\ \ k>j+5,
$$
for $k>5$ we have
\begin{align*}
\|\cR^{\bba}_k( f\circ g)\|_\bbp&\le\sum_{j\geq k-5}\|\cR^{\bba}_k( \cR^{\bba}_jf\wt\cR^{\bba}_j g)\|_\bbp
\lesssim\sum_{j\geq k-5}\| \cR^{\bba}_jf\wt\cR^{\bba}_j g\|_\bbp\\
&\leq \sum_{j\geq k-5}\|\cR^{\bba}_jf\|_{\bbp_1}\|\wt\cR^{\bba}_j g\|_{\bbp_2}
\lesssim  \sum_{j\geq k-5}\|\cR^{\bba}_jf\|_{\bbp_1}\| g\|_{\bbp_2}.
\end{align*}
Since $s>0$, by Young's inequality for discrete summation, we get
$$
\|f\circ g\|_{\bB^{s,q}_{\bbp;\bba}}
\lesssim\left(\sum_{k\geq 0}2^{skq}\left(\sum_{j\geq k-5}\|\cR^{\bba}_jf\|_{\bbp_1}\right)^q\right)^{1/q}\| g\|_{\bbp_2}
\lesssim \|f\|_{\bB^{s,q}_{\bbp_1;\bba}}\| g\|_{\bbp_2}.
$$
Combining the above estimates, we complete the proof of  (iii).

(iv)
For $\theta=0,1$, \eqref{Sob} is direct by \eqref{Ber}. Below we assume $\theta\in(0,1)$.
Let $\bbp_2\leq\bbw\in[1,\infty]^N$ and $s_2\in\mR$ be defined by 
\begin{align}\label{DD02}
\tfrac{1}{\bbp}=\tfrac{1-\theta}{\bbp_1}+\tfrac{\theta}{\bbw},\ \ s_1(1-\theta)+s_3\theta=s.
\end{align}
By H\"older's inequality, we have
\begin{align*}
\|\cR^\bba_j f\|_{\bbp}&=\|(\cR^\bba_j f)^{1-\theta}(\cR^\bba_j f)^{\theta}\|_{\bbp}
\leq\|\cR^\bba_j f\|^{1-\theta}_{\bbp_1}\|\cR^\bba_j f\|^\theta_{\bbw}.
\end{align*}
Hence,
\begin{align}
2^{js}\|\cR^\bba_j f\|_{\bbp}&\leq\Big(2^{js_1}\|\cR^\bba_j f\|_{\bbp_1}\Big)^{1-\theta}
\Big(2^{js_3}\|\cR^\bba_j f\|_{\bbw}\Big)^\theta.\label{DD03}
\end{align}
Note that by \eqref{DD01} and \eqref{DD02},
$$
s_2=s_3+\bba\cdot(\tfrac d{\bbp_2}-\tfrac d{\bbw}),
$$ 
and by \eqref{Ber},
$$
\|\cR^\bba_j  f\|_{\bbw}\lesssim 2^{j\bba\cdot(\frac d{\bbp_2}-\frac d{\bbw})}\|\cR^\bba_j  f\|_{\bbp_2},\ j\geq 0.
$$
Substituting it into \eqref{DD03}, we obtain
$$
2^{js}\|\cR^\bba_j f\|_{\bbp}\lesssim\Big(2^{js_1}\|\cR^\bba_j f\|_{\bbq}\Big)^{1-\theta}
\Big(2^{js_2}\|\cR^\bba_j f\|_{\bbq}\Big)^\theta.
$$
Thus we obtain (iv)  by H\"older's inequality with respect to $j$.

(v)  Note that for $\xi=(\xi_1,\xi_2)$,
$$
\widehat{\cR^x_k\cR^\bba_jf}(\xi)=\phi_k(\xi_1)\phi^\bba_j(\xi)\hat f(\xi).
$$
It is easy to see that 
$$
\phi_k(\xi_1)\phi^\bba_j(\xi)=0,\ \ k>(1+\alpha)(j+2).
$$
Hence,
\begin{align*}
\sup_{k,j\geq 1}2^{\frac{ks_0}{1+\alpha}}2^{ js_1}\|\cR^x_k\cR^\bba_{j}f\|_\bbp
&=\sup_{k\leq (1+\alpha)(j+2)}2^{\frac{ks_0}{1+\alpha}}2^{ js_1}\|\cR^x_k\cR^\bba_{j}f\|_\bbp\\
&\leq\sup_{k\leq (1+\alpha)(j+2)}2^{\frac{k(s_0-\eps)}{1+\alpha}}2^{ j(s_1+\eps)}\|\cR^x_k\cR^\bba_{j}f\|_\bbp.
\end{align*}

(vi) By definition and Bernstein's inequality we have
\begin{align*}
\|f\|_{\bB^{\theta \gamma, (1-\theta)\beta}_{\bbp; x,\bba}}
&=\sup_{k,j\geq 0}\left(2^{\frac{k \gamma}{1+\alpha}}\|\cR^x_k\cR^\bba_{j}f\|_\bbp
\right)^\theta\left(2^{ j\beta}\|\cR^x_k\cR^\bba_{j}f\|_\bbp\right)^{1-\theta}\\
&\lesssim\sup_{k,j\geq 0}\left(2^{\frac{k\gamma}{1+\alpha}}\|\cR^x_kf\|_\bbp
\right)^\theta\left(2^{ j\beta}\|\cR^\bba_{j}f\|_\bbp\right)^{1-\theta}=
\|f\|^\theta_{\bB^{ \gamma/(1+\alpha)}_{\bbp;x}}
\|f\|^{1-\theta}_{\bB^\beta_{\bbp;\bba}}.
\end{align*}
The proof is complete.
\end{proof}

\section{A priori $L^1$-estimates for kinetic equation}

In this appendix we use the probabilistic representation to show some a priori estimates for linear kinetic Fokker-Planck equation
with smooth coefficients, which has independent interest. 

\bl\label{LEC1}
Let $T>0$ and $\alpha\in(0,2]$. Let  $\DD\in L^1([0,T]; C^1_b(\mR^{2d}))$ and
$u\in C([0,T];C^2_b(\mR^{2d}))$ be a smooth solution of the linear Fokker-Planck equation
$$
\p_t u=\Delta^{\alpha/2}_vu-v\cdot\nabla_x u-\div_v(\DD u).
$$
For any $t\in[0,T]$, it holds that
\begin{align}\label{AM2}
\|u(t)\|_{\mL^\infty}\leq\|u(0)\|_{\mL^\infty}\e^{\int^t_0\|(\div_v \DD_s)^-\|_{\mL^\infty}\dif s}
\end{align}
and
\begin{align}\label{AM1}
\|u(t)\|_{\mL^1}\leq\|u(0)\|_{\mL^1}.
\end{align}
\el
\begin{proof}
(i)  Fix $t_0\in(0,T]$. Let $\bD_t(z):=\bD_t(x,v):=(-v,-\DD_t(x,v))$ and $\sigma:=(0,\mI)^*$.
Let $Z_t(z)=Z_t$ solve the following SDE 
$$
\dif Z_t=\bD_{t_0-t}(Z_t)\dif t+\sigma\dif L^{(\alpha)}_t, \ \ Z_0=z,\ \ t\in[0,t_0].
$$
Since we are considering the additive noise and $\DD\in L^1([0,T]; C^1_b(\mR^{2d}))$, it is easy to see that
$$
z\mapsto Z_t(z)\mbox{ is a $C^1$-diffeomorphism,}
$$
and
$$
\nabla Z_t(z)=\mI_{2d\times 2d}+\int^t_0\nabla\bD_{t_0-s}(Z_s(z))\nabla Z_s(z)\dif s.
$$
Let $G_t$ be the inverse of matrix $\nabla Z_t(z)$, which solves the following ODE
\begin{align}\label{Mx7}
G_t(z)=\mI_{2d\times 2d}-\int^t_0G_s(z)\cdot\nabla\bD_{t_0-s}(Z_s(z))\dif s.
\end{align}
Note that
$$
\p_t u=\Delta^{\alpha/2}_vu-v\cdot\nabla_x u-\DD\cdot \nabla_v u-\div_v \DD u.
$$
By applying It\^o's formula to the function $(s,z)\mapsto u(t_0-s,z)$, we have
$$
u(t_0-t, Z_t(z))=u(t_0,z)-\int^t_0(\div_v \DD u)(t_0-s, Z_s(z))\dif s+\mbox{a martingale}.
$$
Using product It\^o's formula we get
$$
\mE \left(u(t_0-t, Z_t(z))\e^{-\int^t_0\div_v \DD_{t_0-s}(Z_s(z))\dif s}\right)=u(t_0,z).
$$
Hence,
$$
|u(t_0,z)|\leq\|u(0)\|_{\mL^\infty}\e^{\int^{t_0}_0\|(\div_v \DD_{t_0-s})^-\|_{\mL^\infty}\dif s}.
$$
Replacing $t_0$ by $t$, we obtain \eqref{AM2}. Moreover, by the change of variables, we also have
\begin{align}
\|u(t_0)\|_{\mL^1}&\leq\e^{\int^{t_0}_0\|(\div_v \DD_{t_0-s})^-\|_{\mL^\infty}\dif s}\mE\left(\int_{\mR^{2d}} |u(0, Z_{t_0}(z))|\dif z\right)\no\\
&\leq\e^{\int^{t_0}_0\|(\div_v \DD_{t_0-s})^-\|_{\mL^\infty}\dif s}\mE\left(\int_{\mR^{2d}} |u(0, z)|\det(\nabla Z^{-1}_{t_0}(z))\dif z\right)\no\\
&\leq\e^{\int^{t_0}_0\|(\div_v \DD_{t_0-s})^-\|_{\mL^\infty}\dif s}\sup_{z}\mE(\det G_{t_0}(z))\|u(0)\|_{\mL^1}\no\\
&\leq\e^{\int^{t_0}_0(\|(\div_v \DD_{t_0-s})^-\|_{\mL^\infty}+\|(\div_v \DD_{t_0-s})^+\|_{\mL^\infty})\dif s}\|u(0)\|_{\mL^1}\no\\
&=\e^{\int^{t_0}_0\|\div_v \DD_s\|_{\mL^\infty}\dif s}\|u(0)\|_{\mL^1},\label{Mx6}
\end{align}
where we have used that $\nabla Z^{-1}_{t_0}(z)=G_{t_0}\circ Z^{-1}_{t_0}(z)$ in the third inequality, and the last inequality is due to 
(cf. \cite{MB02})
$$
\det G_t(z)=1-\int^t_0\det G_s(z)\cdot\div_v H_{t_0-s}(Z_s(z))\dif s.
$$

(ii) %For proving \eqref{AM1}, we priorly do not know whether $u(t)\in\mL^1$. So we have to use some cutoff technique.
Let $\psi: \mR\to\mR$ be a smooth convex function with 
\begin{align}\label{Mx66}
|\psi(r)|+|\psi'(r)|\leq Cr,\ \ |\psi''(r)|\leq C.
\end{align}
Let $\chi\in C^\infty_c(\mR^{d})$ be a cutoff function in $\mR^d$ with $\chi(v)=1$ for $|v|\leq 1$ and
$\chi(v)=0$ for $|v|>2$. For $R>0$, define
$$
\chi_R(v):=\chi(v/R).
$$
By the chain rule we have
$$
\p_t \psi(u)\chi_R=\psi'(u)\Delta^{\alpha/2}_vu \chi_R-v\cdot\nabla_x\psi(u)\chi_R-\div_v(\DD u)\psi'(u)\chi_R.
$$
Integrating both sides over $[0,t]\times\mR^{2d}$, we get
$$
\int\psi(u(t))\chi_R=\int\psi(u(0))\chi_R+\int^t_0\int\Delta^{\alpha/2}_vu \psi'(u)\chi_R
-\int^t_0\int \div_v(\DD u)\psi'(u)\chi_R,
$$
where we have used that
$$
\int v\cdot\nabla_x\psi(u)\chi_R=\int \div_x((v\chi_R) \psi(u))=0.
$$
By \eqref{Mx6} and \eqref{Mx66}, letting $R\to\infty$ and by the dominated convergence theorem, we get
$$
\int\psi(u(t))=\int\psi(u(0))+\int^t_0\int\Delta^{\alpha/2}_vu \psi'(u)-\int^t_0\int \div_v(\DD u)\psi'(u).
$$
We consider the second term in the right hand side. For $\alpha=2$, we directly have
$$
\int\Delta_vu\psi'(u)=-\int|\nabla_vu|^2\psi''(u)\leq 0.
$$
For $\alpha\in(0,2)$, by \eqref{Mb3} and symmetry we have
\begin{align*}
&\int_{\mR^d}\Delta^{\alpha/2}_vu(v)\psi'(u(v))\dif v
=\lim_{\eps\downarrow 0}\int_{\mR^d}\int_{|w|\geq \eps}\psi'(u(v))\frac{u(v+w)-u(v)}{|w|^{d+\alpha}}\dif w\dif v\\
&\qquad=\lim_{\eps\downarrow 0}
\frac12\int_{\mR^d}\int_{|w|\geq\eps}(\psi'(u(v))-\psi'(u(v+w)))\frac{u(v+w)-u(v)}{|w|^{d+\alpha}}\dif w\dif v\\
&\qquad=-\lim_{\eps\downarrow 0}\frac12\int_{\mR^d}\int_{|w|\geq\eps}\int^1_0\psi''(s u(v)+(1-s)u(v+w))\dif s\frac{(u(v+w)-u(v))^2}{|w|^{d+\alpha}}\dif w\dif v\leq 0.
\end{align*}
Hence,
\begin{align}\label{AC11}
\int\psi(u(t))\leq\int\psi(u(0))+\int^t_0\int (\DD u\cdot\nabla_v u)\psi''(u).
\end{align}
Now we take
$$
\psi(r)=\psi_\delta(r):=\sqrt{r^2+\delta}-\sqrt\delta,\ \ \delta>0.
$$
Clearly,
$$
\lim_{\delta\downarrow 0}\psi_\delta(r)=|r|,\ |\psi_\delta(r)|\leq|r|,\ \psi'_\delta(r)=\frac{r}{\sqrt{r^2+\delta}},\ \ \psi''_\delta(r)=\frac{\delta}{(r^2+\delta)^{3/2}}.
$$
By \eqref{AC11}, we have
$$
\int\psi_\delta(u(t))\leq\int\psi_\delta(u(0))+\int^t_0\int (\DD u\cdot\nabla_v u)\psi''_\delta(u).
$$
Noting that
$$
|\psi''_\delta(u)|\cdot |\DD|\cdot|u|\cdot|\nabla_v u|\leq |\DD|\cdot|\nabla_v u|,
$$
by taking limits $\delta\to 0$ and the dominated convergence theorem, we get
\begin{align*}
\int|u(t)|
\leq\int|u(0)|+\lim_{\delta\to 0}\int^t_0\int |\psi''_\delta(u)|\cdot |\DD|\cdot|u|\cdot|\nabla_v u|=\int|u(0)|.
\end{align*}
The proof is complete.
\end{proof}

\end{appendix}

\end{document}